\documentclass[11pt]{article}
\usepackage{amsfonts,amssymb,stmaryrd,amscd,amsmath,latexsym,amsbsy,graphicx}

\newcommand{\beq}{\begin{eqnarray*}}
\newcommand{\eeq}{\end{eqnarray*}}
\newcommand{\beqn}{\begin{eqnarray}}
\newcommand{\eeqn}{\end{eqnarray}}
\newcommand{\ub}[2]{\stackrel{\; #2}{#1}}

\title{The eigenvalue equation on the Eguchi-Hanson space}
\date{}
\author{Andreas Malmendier\footnote{email: malmendi@mit.edu} \\ \\
\em \small Department of Mathematics,\\ 
\em \small Massachusetts Institute of Technology,\\
\em \small Cambridge, MA 02139-4307}

\begin{document}
\maketitle
\begin{abstract}
We consider the eigenvalue equation for the Laplace-Beltrami operator
acting on scalar functions on the non-compact Eguchi-Hanson space. The
corresponding differential equation is reducible to a confluent Heun
equation with Ince symbol $[0,2,1_2]$. We construct approximations for
the eigenfunctions and their asymptotic scattering phases with the
help of the Liouville-Green approximation (WKB). Furthermore, for
specific discrete eigenvalues obtained by a continued T-fraction we
construct the solution by the Frobenius method and determine its
scattering phase by a monodromy computation.
\end{abstract}

\section{Introduction}
The Eguchi-Hanson metric \cite{Eguchi:1979gw} is a self-dual, asymptotically
locally Euclidean (ALE) metric on the cotangent bundle of
the 2-sphere $\mbox{T}^*S^2$.
Geometrically, this corresponds to a Ricci-flat metric on the smooth
resolution space
of an $A_1$-singularity -- such a singularity looks like the origin of
$\mathbb{C}^2/\mathbb{Z}_2$ where the $\mathbb{Z}_2$-group
acts by point reflection at the origin.

Apart from Euclidean quantum gravity this metric has physical applications
in string compactification. In fact, it is well
known that the only nontrivial two dimensional Calabi-Yau manifold -- this
is a K3-surface -- can be
obtained from the $\mathbb{Z}_2$-orbifold limit
$\mathbb{T}^2_\mathbb{C}/\mathbb{Z}_2$ by blowing up its 16
$A_1$-singularities
(where $\mathbb{T}^2_\mathbb{C}$ is the complex two-dimensional torus).
By gluing the Eguchi-Hanson metric together with the flat torus metric one
can explicitly construct an almost Ricci-flat
metric on a K3-surface \cite{Bozhkov:1988} which is related to a
Ricci-flat metric on K3 by a gauge transformation
\cite{Taubes:1982}.

In this paper we shall examine the eigenvalue equation for the
Laplace-Beltrami operator on the Eguchi-Hanson space. The
problem is interesting from a mathematical point of view since the
differential equation is separable and reduces to an
ordinary differential equation which due to its singularities can be
identified as confluent Heun equation with corresponding Ince symbol 
$[0,2,1_2]$ (see \cite{DMR:1978} for definitions). The Heun equation is an 
ordinary differential equation with four regular singularities on the punctured Riemann sphere. By 
coalescing two of the regular singularities to one irregular singularity 
one obtains the confluent Heun equation 
(analogously to the procedure by which one obtains the confluent
hyper-geometric differential equation from the hyper-geometric one).

The problem is also interesting from a physical point of view since the
functions that are obtained by gluing the eigenfunctions on the
Eguchi-Hanson space together with the well known
eigenfunctions on the flat torus describe the quantum mechanical limit
of string fields on $\mbox{K}3$.

The plan of the paper is as follows: In Sect.~\ref{eigenvalue_equation} we
introduce the ordinary differential equation that
describes the radial part of the eigenvalue equation of the Laplace operator
on the Eguchi-Hanson space. In Sect.~\ref{WKB} we
construct its solutions and their corresponding scattering phases by the
Liouville-Green approximation (WKB). This extends and corrects a
result in \cite{Mign:1991}. In \cite{Mign:1991} the author used the ad-hoc version
of the Liouville-Green approximation. This version did not give approximations for 
the wave functions which are valid over the whole range. 
Furthermore, we give error bounds for the constructed
solutions and their scattering phases. The explicit calculations can be found in the
Appendix (cf. App.~\ref{WKB_calculation} and \ref{error_bounds}).
In Sect.~\ref{Con_Frac} we construct the exact solutions for specific discrete
eigenvalues by the Frobenius method. These special
values are given by the vanishing condition for a continued fraction. This
approach is similar to the treatment of the generalized
spheroidal wave equation in \cite{Wilson:1928}, \cite{Leav:1986}.
In addition, we can determine the exact scattering phase by
a monodromy computation. Finally, we will show how this information can be used
to compute the asymptotic scattering phase for this discrete set of eigenvalues.
In Sect.~\ref{results} we present the numerical results obtained by the method of 
Sect.~\ref{WKB} and Sect.~\ref{Con_Frac} and show that they match up to a high
accuracy. In Sect.~\ref{conclusions} we give the conclusions of this article and
a brief outlook.

\section{The Eigenvalue Equation}
\label{eigenvalue_equation}
In this section we introduce the eigenvalue equation for the
Laplace-Beltrami operator of the Eguchi-Hanson space.
Using the $SU(2)\times U(1)$ symmetry of the Eguchi-Hanson space this 
eigenvalue equation will reduce
to an ordinary differential equation
of second order. A complete derivation of these and similar results can be
found in \cite{Gibbons:1986df},\cite{Perr:1978}.

The Riemannian metric on $\mbox{T}^* S^2$ that is known as Eguchi-Hanson
metric is Hyperk\"ahler. Moreover, a complex structure which is compatible
with the Hyperk\"ahler structure can be introduced by identifying
$\mbox{T}^*S^2$ with the complex manifold $\mbox{T}^* \mathbb{CP}^1$. In
particular, the latter can be covered by two coordinate charts $U \cong 
U^\prime \cong \mathbb{C}^2$ with coordinates $(u,\xi), \, (u^\prime,
\xi^\prime)$, respectively. Here, $u, \, u^\prime$ denote Euclidean
coordinates on the base $\mathbb{C}P^1$, and $\xi, \, \xi^\prime$ parameterize
the fiber of the bundle $\mbox{T}^*\mathbb{CP}^1 \longrightarrow 
\mathbb{CP}^1$. This means that
\beqn
\label{charts}
U \cap U^\prime \cong \mathbb{C}^\ast \times \mathbb{C}\,, \quad
(u^\prime, \xi^\prime) \in U\cap U^\prime: \; (u^\prime, \xi^\prime)=
(u^{-1}, u^2\xi)\,.
\eeqn
Note that $\mbox{T}^* \mathbb{CP}^1$ describes the minimal resolution
$\widetilde{\mathbb{C}^2/\mathbb{Z}_2}$ of $\mathbb{C}^2/\mathbb{Z}_2$,
where $\mathbb{Z}_2$ acts on $(z^1,z^2) \in \mathbb{C}^2$ by
$(z^1,z^2) \mapsto (-z^1,-z^2)$. To see this, on $U_0:= \left\lbrace (u,\xi)
\in U \mid \xi\not = 0\right\rbrace$ use $(u,\xi) = (\frac{z^1}{z^2}, 
(z^2)^2)$, and analogously for the other chart. Then (\ref{charts}) is
given by $z^1 \leftrightarrow z^2$, and by the above we identify 
$\mathbb{C}^2/\mathbb{Z}_2$ with $\mbox{T}^* \mathbb{CP}^1$ with the zero 
section removed. 
The exceptional divisor of $\widetilde{\mathbb{C}^2/\mathbb{Z}_2}$ therefore
corresponds to the zero section of $\mbox{T}^*\mathbb{CP}^1$.
In these coordinates, i.e. outside the exceptional divisor, 
the Eguchi-Hanson metric takes the form

\beq
g^{\mbox{\scriptsize EH}}_{ \imath \bar\jmath} & = & 
\frac{\sqrt{c^4+R^4}}{R^2}
\left\lbrace \delta_{\imath \bar\jmath} - \frac{c^4 \, z_\imath
z_{\bar\jmath}}{R^2 \, (c^4+R^4)} \right\rbrace \,, \\
\mbox{where} \quad R^2 & = & |z^1|^2 + |z^2|^2 \;,
\eeq
and $c>0$ is the parameter of the Eguchi-Hanson metric. Introducing the
Euler angles $(\theta, \phi, \psi)$ on $S^3$ by 
\beq
\left(\begin{array}{c} z^1 \\ z^2 \end{array} \right) = R \,
                         \left(\begin{array}{c}  \cos{\frac{\theta}{2}} \,
e^{\frac{i}{2}(\psi + \phi)} \\
                         \sin{\frac{\theta}{2}} \, e^{\frac{i}{2}(\psi -
\phi)} \end{array} \right) \,,
\eeq
we can describe coordinates (in the coordinate patch with $z_2\not = 0$) by a radial part $R \in \mathbb{R}_{>0}$
and the angles $\theta \in (0,\pi)$, $\phi \in [0,2\pi)$, $\psi \in [0,4\pi)$.
In the following, it will be more convenient to use $r:=\sqrt[4]{R^4 + c^4}$
instead of $R$. 

The singularity of $\mathbb{C}^2/\mathbb{Z}_2$ then sits
at $r=c$. On the exceptional divisor, the Eguchi-Hanson metric takes the
form of the Fubini-Study metric on $\mathbb{CP}^1$, and therefore assigns
volume $\pi$ to the exceptional divisor \cite{Eguchi:1979gw}.

The $\mathbb{Z}_2$-operation mentioned above takes the form
\beq
\mathbb{Z}_2: \qquad \quad \left(\begin{array}{c}  z^1 \\ z^2 \end{array} 
\right) &
\mapsto &
       - \left(\begin{array}{c}  z^1 \\ z^2 \end{array} \right) \\
       \big( \theta, \psi , \phi \big) & \mapsto & \big( \theta, \psi + 2
\pi, \phi \big) \;.
\eeq
In these new coordinates the metric takes the form originally found by
Eguchi and Hanson.
The Lie group $SO(3)$ acts on itself by multiplication from the left and from 
the
right. Let us call the vector fields
that generate the  right-multiplication $(\xi_1, \xi_2, \xi_3)$, and the 
ones
that generate the
left-multiplication $(\bar\xi_1, \bar\xi_2, \bar\xi_3)$. Since
$SO(3) \cong SU(2)/\mathbb{Z}_2 \cong S^3/\mathbb{Z}_2$
we have an action of these vector fields on the Eguchi-Hanson space. In
\cite{Gibbons:1986df} the authors derive
that $\forall i,j: [\xi_i, \bar\xi_j]=0$ and that $(\xi_1, \xi_2, \xi_3,
\bar\xi_1)$ define Killing vector fields
that generate an $\mathfrak{su}(2)_L \oplus \mathfrak{u}(1)_R$ symmetry
algebra of the Eguchi-Hanson space.

Let $\Delta$ denote the Laplace-Beltrami operator on functions that is
associated to the Eguchi-Hanson metric. Furthermore, let $\Psi$ denote
a smooth function on the Eguchi-Hanson space. We consider the eigenvalue 
equation for the real positive eigenvalue\footnote{We are considering 
only the scattering case since we want to glue 
together the eigenfunctions on the Eguchi-Hanson space and the flat torus.
However, eigenfunctions for $E<0$ would correspond to bounded states.}, i.e.~ 
$(E+\Delta) \Psi = 0$ for $E>0$.
The operator $E+\Delta$ can be expressed merely  in terms of 
$\frac{d}{dr}, \, \xi^2=\sum_{i=1}^3 \xi_i^2, \, \xi_3^2$, i.e.
\beq
 E + \Delta
= E + \left(1-\frac{c^4}{r^4} \right) \frac{d^2}{d r^2} + \left( \frac{3}{r}
 + \frac{c^4}{r^5} \right) \frac{d}{d r}
      + \frac{4 \xi^2}{r^2} + \frac{4c^4 \, \xi_3^2}{r^2(r^4-c^4)}\;.
\eeq
Due to their commutation relations we can diagonalize the operators $\xi_3,
\bar{\xi}_3, \xi^2$ simultaneously and expand the eigenfunctions in terms 
of Wigner functions $D^{j}_{qm}(\theta,\psi,\phi)$ (see~\cite{Abra:1955} for definitions), i.e.
\beq
\Psi(r,\theta,\psi,\phi) = \sum_{j=0}^\infty \sum_{q,m=-j}^j \alpha^j_{qm} 
\; A(j,q,\beta|z) \; D^j_{qm}(\theta,\psi,\phi) \,,
\eeq
where $z:=\frac{r^2}{c^2}$, $\beta:=\frac{c^2\,E}{4}$, and $\alpha^j_{qm}$ are 
complex coefficients, and the $A(j,q,\beta|z)$ are functions that depend 
only on the coordinate $z$. The Wigner functions fulfill
\beq
\bar{\xi}_3 \, D^{j}_{qm} & = & - i \, m \, D^{j}_{qm} \;,\\
       \xi_3 \, D^{j}_{qm} & = &   i \, q \, D^{j}_{qm} \;,\\
\Big( \xi_1^2 + \xi_2^2 + \xi_3^2 \Big) \, D^{j}_{qm} & = &   -j(j+1) \,
D^{j}_{qm} \;,\\
\mbox{where} \quad j \in \mathbb{N}\,,&& \; q,m \in \lbrack -j , j 
\rbrack_{\mathbb{N}} \,.
\eeq
Here, $j,q$ are the quantum numbers that label the $SU(2)$ representation.
The eigenvalue equation for $\Psi$ then reduces to an ordinary differential
equation for $A(j,q,\beta|.)$,
\beqn
\label{ODE}
0 & = & \left[ \frac{d^2}{d z^2} + \frac{2\,z}{z^2-1} \frac{d}{d z}
      + \frac{(\beta z - j(j+1))(z^2-1)- q^2}{(z^2-1)^2} \right] \;
A(j,q,\beta | z) \;.
\eeqn
The differential equation (\ref{ODE}) has three singularities 
which we have summarized in the following table
\begin{center}
\begin{tabular}{c|c|l|c}
$z$ & $r$ & singularity & roots\\
\hline
  1   & $c$           & regular   & $\pm \frac{q}{2}$\\
-1   & $i \, c$      & regular   & $\pm \frac{q}{2}$\\
$\infty$ & $\infty $ & irregular & --
\end{tabular}\\
(for $q=0$ the regular singularities are logarithmic).
\end{center}
One should notice that our differential equation has an irregular 
singularity at infinity, thus
is {\em not} of Fuchsian type.
More precisely, it is a confluent Heun equation with Ince symbol $[0,2,1_2]$.

To construct a continuous solution valid on the whole Eguchi-Hanson sphere
we have to choose the regular boundary condition at $z=1$ (then we can extend
this solution to $z=1$, i.e. the blown up $A_1$-singularity). This means that 
we are looking for the \textsl{recessive} solutions that behave like
$(z-1)^\frac{q}{2}$ for $z\to 1$. This means that all the eigenfunctions 
obtain a constant value on the entire $\mathbb{CP}^1$ which gives the
exceptional divisor of the blow-up. By approaching the singular point
of $\mathbb{C}^2/\mathbb{Z}_2$ with different slopes, one reaches different
points in the exceptional divisor. But for a well-defined solution of
the differential equation, the limit does not depend on the chosen slope. 

We remark that the differential equation (\ref{ODE}) does not depend on $q$ but only on $q^2$.
Therefore, it suffices to restrict ourselves to $q \ge 0$. Then, the function
$\Psi$ takes the following form
\beq
\Psi(r,\theta,\psi,\phi) = \sum_{j=0}^\infty \sum_{q=0}^j A(j,q,\beta|z)
\, \sum_{m=-j}^j \left( \alpha^j_{qm} \, D^j_{qm}(\theta,\psi,\phi) +
 \alpha^j_{-qm} \, D^j_{-qm}(\theta,\psi,\phi) \right) \,.
\eeq
On these functions the differential operator of Eq.~(\ref{ODE}) is
self-adjoint: Since the differential operator in Eq.~(\ref{ODE})
is already formally self-adjoint the statement follows from the
application of the results in \cite[Sect. 3.9]{Codd:1955} and a
detailed study of the singularity at $z=1$ for $q=0$, where the 
differential equation is of a \textsl{limit-circle} type,
and for $q\ge 1$, where it is of \textsl{limit-point} type 
-- for more details see \cite{Codd:1955}.

The specific problems that arise in the treatment of the differential
equation (\ref{ODE}) are due to the structure of singularities. In
particular, the major issue is to treat {\em three} singularities
-- one being an \textsl{irregular} singularity -- at the same time.
The problem is to derive the connection between the different bases 
which provide expansions
of the solutions in the neighborhood of the singular points. Another
important question is how given a system of solutions the solutions
transform into each other
when passing through a cycle around the corresponding singularity.

For the thrice-punctured Riemann sphere with only regular singularities
(this is a generalized hyper-geometric equation) this can be done in terms
of the Meijer transcendental functions -- as they were recently applied in 
\cite{Greene:2000ci}. Here, the key technique lies in a representation
of the solutions in terms of Mellin-Barnes integrals.

For only two singularities -- one being regular, one being singular --
techniques can be applied which are familiar from the treatment of the
Bessel differential equation. This corresponds to the confluent 
hyper-geometric equation, i.e. a hyper-geometric equation where two
regular singularities are coalescing and forming one irregular singularity. 
Here, the key technique is a generalized Borel transformation 
\cite{Gura:1994b}, \cite{Gura:1994a} which relates a cycle around the
regular singularity at $z=1$ to the irregular one at $z=\infty$.

In the case of two regular and one irregular singularities we cannot apply 
Mellin-Barnes integrals. Due to the irregular singularity solutions are
oscillating at large real values. Thus, the conditions for convergence for
the Mellin-Barnes integrals when continuing into the complex plane to close 
the path of integration to a cycle are not satisfied any more. 

On the other hand, due to the third singularity at $z=-1$ a cycle around 
infinity is also not homologous to a cycle around $z=1$. Thus, an 
asymptotic expansion of the regular solution cannot be derived by 
a Borel transformation. However, as we will show in Sect.~\ref{monodromy},
in some particular cases one can still use a similar argument.

\section{WKB type solutions}
\label{WKB}
In this section we give approximations to the recessive solutions of
differential equation (\ref{ODE}) that can be obtained by the 
Liouville-Green
approximation (WKB). The explicit derivation of the results can be found in
App.~\ref{case_I} to \ref{case_IV}.

From now on we suppress the labels $\beta,j,q$ in $A(\beta,j,q|z)$.
Substituting $A(z) = (z^2-1)^{-\frac{1}{2}} \, w(u,z)$ in Eq.~(\ref{ODE})
we find
\beqn
\label{WKB-ODE}
\frac{d^2}{dz^2} w(u,z)   & = & \left[ - u^2 f(z) + g(z) \right] w(u,z) \,,
\eeqn
where $u^2 = \beta$, $g(z) = - (z^2-1)^{-2}$, $a=\frac{j(j+1)}{\beta}$,
$b=\frac{q^2}{\beta}$ and
\beq
f(z) & = & \frac{(z-a)(z^2-1)-b}{(z^2-1)^2} \;.
\eeq
We remark that due to the possible values for $j,q,\beta$ we will always 
have $1\le b+1<a$ and $0\le a$.

Eq.~(\ref{WKB-ODE}) is the standard form of a differential equation 
of second order considered for Liouville-Green approximation (WKB).  
Following the
discussion in \cite{Olve:1974}, the specific approximation depends on
the order of the pole at $z=1$ as well as the number of \textsl{simple
transition points} (tps), i.e. simple zeros, in the region $z>1$ which is of
interest in connection with the geometry discussed in 
Sect.~\ref{eigenvalue_equation}.
A simple analysis of the function $f$ shows that we have to deal with the 
four different cases (I to IV) which are summarized in the table below.

\begin{center}
\begin{tabular}{c|l|l|l}
  case & condition           & order of pole at $z=1$  & tp for $z>1$\\
\hline 
  I    & $1 < b+1< a$        & 2                       & yes\\
  II   & $ b=0, 1<a $        & 1                       & yes\\
  III  & $ b=0, a=1 $        & 0                       & no\\
  IV   & $ b=0, 0 \le a < 1$ & 1                       & no\\
\end{tabular}
\end{center}
For the different cases one can then apply the Liouville-Green approximation: 
First, one performs a Liouville-transformation of the variable $z$ to the new
variable $\zeta$. The transformation is given by an integral equation of the form 
\beq
 \forall z \ge z_0: \, G(\zeta)= \int_{z_0}^z \sqrt{|f(t)|} \, dt
\eeq
such that $\zeta$ and $z$ are analytic functions of each other, and where the function $G$ and the 
point $z_0$ depend on the case (I to IV) we are dealing with. Simultaneously, we replace the
function $w$ by a function $W$ according to $w(u,z)= \sqrt{|\frac{dz}{d\zeta}|} \, W(u,\zeta)$.
The choices are as follows:
\begin{center}
\begin{tabular}{c|l|l}
  case & $G(\zeta)$ & $z_0$\\
\hline
  I    & $\frac{2}{3}(-\zeta)^\frac{2}{3}$ & $>1$ such that $f(z_0)=0$\\
  II   & $\zeta^\frac{1}{2}(\zeta-\alpha)^\frac{1}{2} 
          - \frac{\alpha}{2} \, \ln{\left(\frac{2\zeta-\alpha+2\zeta^\frac{1}{2}(\zeta-\alpha)^\frac{1}{2}}{\alpha}\right)}$ & $a$\\
  III  & $\zeta$ & $1$\\
  IV   & $(-\zeta)^\frac{1}{2}$ & $1$\\
\end{tabular}
\end{center}
and
\beq
\alpha := \frac{2}{\pi} \int_1^a dt \, \sqrt{-f(t)} \;.
\eeq
The aim of this transformation is to transform the differential equation (\ref{WKB-ODE})
to a differential equation of the form 
\beq
 \frac{d^2}{d\zeta^2} W(u,\zeta) & = & \Big( T_1(u,\zeta) + T_2(\zeta) \Big) \, W(u,\zeta)\,,
\eeq
where $T_1$ and $T_2$ are real-valued functions such that the
approximating differential equation obtained by omitting
$T_2$ has solutions which are functions of a single variable. However,
this transformation has been done in a way that the approximate
solution will still reflect the right behavior of the solution at the
singularity at $z=1$ and for $z \to \infty$. From now on we suppress
the dependence of $u$ in $w$ and $W$.

We have found the following approximations for the solutions of
differential equation (\ref{ODE}) that are regular at $z=1$
\begin{center}
\begin{tabular}{c|l|l}
  case & approx. function & form\\
\hline
  I    & Airy       & $w(z) = \sqrt[4]{- \frac{\zeta}{f(z)}} \; \mbox{Ai}(u^\frac{2}{3} \zeta)$\\
  II   & Whittaker  & $w(z) = \sqrt[4]{\frac{\zeta-\alpha}{2u \zeta \,f(z)}} \;
 e^{-\frac{i\pi}{4}} \, \mbox{M}_{\frac{i u\alpha}{2},0} \left(2i u\zeta\right)$\\
  III  & Bessel     & $w(z) = \frac{\zeta^\frac{1}{2}}{\sqrt[4]{f(z)}} \, \mbox{J}_0(u\zeta)$\\
  IV   & Bessel     & $w(z) = \frac{|\zeta|^\frac{1}{2}}{\sqrt[4]{4|\zeta|f(z)}} \,
  \mbox{J}_0(u |\zeta|^\frac{1}{2})$\\
\end{tabular}
\end{center}
Asymptotically, the Eguchi-Hanson metric becomes the flat metric. Therefore, 
the solutions of the differential equation ($\ref{ODE}$) must have the following 
behavior
\beqn
\label{k1}
A(z) \ub{\sim}{z\to \infty} \frac{1}{z^\frac{3}{4}} \sin\left( 2 \sqrt{\beta 
z} + \Delta_{j,q} \right)\,,
\eeqn
where $\Delta_{j,q}$ is called scattering phase. Based on the results of
Olver et al. (see App.~\ref{case_I} to \ref{case_IV} for details), for the scattering 
phase we obtain the results listed below, where in the last column we give
the equation number for the error bounds that are determined in App.~\ref{WKB_calculation}:
\begin{center}
\begin{tabular}{c|l|l}
 case  & $\Delta_{j,q}(\beta)$ & error\\
\hline
  I    & $\sqrt{\beta} \; \lim_{z \to \infty} \left( \int_{z_0}^z
\sqrt{f(t)} \, dt - 2 \sqrt{z} \right) + \frac{\pi}{4} $
  & (\ref{errorI}) \\
  II   & $- 2 \sqrt{\beta} \, \sqrt{1+a} \; \;
\mbox{E}\left[\sqrt{\frac{2}{1+a}} \right]+ \frac{\alpha\sqrt{\beta}}{2}
- \frac{\alpha\sqrt{\beta}}{2} \ln{(\frac{\alpha\sqrt{\beta}}{2})}+
\mbox{arg} \; \Gamma\left( \frac{1}{2} + \frac{i\alpha\sqrt{\beta}}{2}
\right) +  \frac{\pi}{4}$ & (\ref{errorII})\\
  III  & $- \sqrt{8 \beta} + \frac{\pi}{4} $ & (\ref{errorIII}) \\
  IV   & $ (1-a) \, \sqrt{2\beta} \;\;
\mbox{K}\left[\sqrt{\frac{1+a}{2}} \right] -
2 \sqrt{2 \beta} \; \; \mbox{E}\left[\sqrt{\frac{1+a}{2}} \right]
+ \frac{\pi}{4} $ & (\ref{errorIV}) 
\end{tabular}
\end{center}
Here, $\mbox{F}(\phi,m), \mbox{E}(\phi,m)$ denote elliptic integrals of the
first and second kind, and
$\mbox{K}(m) = \mbox{F}(\frac{\pi}{2},m)$ and $\mbox{E}(m) =
\mbox{E}(\frac{\pi}{2},m)$ are the corresponding complete
elliptic integrals (c.f.~\cite[Sect. 17]{Abra:1955}).
This corrects a result in \cite{Mign:1991}.

\section{Solutions related to continued fractions}
\label{Con_Frac}
In this section we determine for which values of $\beta,j,q$ the exact
solution of the differential equation (\ref{ODE}) can be obtained by
a formal power series expansion (Frobenius method) around the singularity
at $z=1$.

For an expansion around the regular singularity at $z=1$ a transformation
according to $\zeta = \frac{1}{2}(z-1)$ is suitable. Substitution of
\beqn
 \label{k2}
A(z) = (z^2 -1 )^\frac{q}{2} \, u(\zeta)
\eeqn
in Eq.~(\ref{ODE}) yields
\beqn
\label{ODE2}
0 & = & \zeta(\zeta+1) \; u''(\zeta) + (q+1)(2\zeta+1) \; u'(\zeta) + \Big(
\beta (2\zeta+1) + \mu \Big) \; u(\zeta)\,,
\eeqn
where $\mu = q(q+1) - j(j+1)$. In a neighborhood of $z=1$ the two linearly
independent solutions can be represented by the series
\beqn
\label{FPS_reg}
u_{\mbox{\scriptsize reg} }(\zeta) & = & \sum_{k=0}^\infty a_k(\beta,j,q) \;
\zeta^k \,,\\
\nonumber
u_{\mbox{\scriptsize sing}}(\zeta) & = & u_{\mbox{\scriptsize reg} }(\zeta)
\; \Big( \ln{\zeta} + B_u \Big)
+ \frac{1}{\zeta^q} \sum_{k=0}^\infty b_k(\beta,j,q) \; \zeta^k \;,
\eeqn
where $B_u$ is a real number. Since the singular solution is not unique (one can 
always add a multiple of the regular solution), the parameter $B_u$ is not uniquely
determined. However, the parameter can be fixed by fixing the asymptotic scattering 
phase of the singular solution.

By the Frobenius method (cf. \cite[Chapt. 3.6]{Rabe:1972}) we obtain the 
coefficients $(a_k)_{k \ge -1}$ as solutions of the following three-term 
recurrence relation
\beqn
\label{REC_REL}
\nonumber
a_{-1} & = & 0 \,,\\
\forall k \ge 0: \qquad a_{k+1} & = & - \frac{k(k+2q+1)+\mu+\beta}{(k+1)(k+q+1)} \, a_k
- \frac{2\beta}{(k+1)(k+q+1)} \,a_{k-1} \;.
\eeqn
Notice that a rescaling of the parameter $a_0$ results in a general rescaling of all the coefficients
$(a_k)_{k \ge -1}$ since Eq.~(\ref{REC_REL}) is linear. However, the crucial information, i.e. the ratio
of $a_1$ and $a_0$ is fixed by the recurrence relation (\ref{REC_REL}).
By standard methods one can show that there are only two types of solutions for $(a_k)_{k\in \mathbb{N}}$ in (\ref{REC_REL}) depending on the radius of convergence $r_c$
of the series $\sum_{k=0}^\infty a_k \, \zeta^k$
\beq
\lim_{k \to \infty} \frac{a_{k+1}}{a_{k}} = \left\{ \begin{array}{ccl} -1
&\quad & r_c = 1 \,,
  \\ 0 &\quad & r_c = \infty \,.\end{array} \right.
\eeq
A solution with $r_c=1$ corresponds to a solution of the recurrence
relation for generic values of $(\beta,j,q)$
whereas the solution with $r_c=\infty$ is the \textsl{minimal
solution}.\footnote{A minimal solution
$(g_n)_{n \in \mathbb{N}}$ is defined by the universal property that for any
other solution
$(h_n)_{n \in \mathbb{N}}$ one obtains $\lim_{n \to \infty} \frac{g_n}{h_n} =
0$.}
The question for which values of the parameters $(\beta,j,q)$ this minimal
solution exists will be considered in Chap.~\ref{minimal}. 

One should mention that it is possible to derive an explicit
representation of the coefficients $(a_k)_{k \in \mathbb{N}}$ by the
use of Babister's inhomogeneous hyper-geometric functions
\cite{Exton:1991b}. However, this approach did not enable us to derive
the asymptotic behavior in the case of a minimal solution. To solve
the differential equation (\ref{ODE2}) by a Laplace transformation 
as in \cite{Exton:1991a} for similar differential equations
or an Euler transformation as in \cite{Kaz:1998} is not possible due to the structure 
of the coefficients in Eq.~(\ref{ODE2}).

To gain a better understanding of solutions to (\ref{ODE2}), let us take a look
at the solution $u_{\mbox{\scriptsize reg}}$ which is regular at $z=1$ ($\zeta =0$) 
in terms of the singularity at $z=-1$ ($\zeta=-1$). Let $v_{\mbox{\scriptsize reg}}(\zeta+1)$ 
be the solution regular at $z=-1$ ($\zeta=-1$). Its continuation to $\zeta \to \infty$ shall have an
asymptotic scattering phase of $\gamma$, which by (\ref{k1}) and (\ref{k2}) gives an asymptotic expansion
\beq
v_{\mbox{\scriptsize reg}}(\zeta+1) & \ub{\sim}{\zeta \to \infty} &
\frac{\sin{(2\sqrt{2\beta\zeta}+\gamma)}}{\zeta^{q+\frac{3}{4}}} \; \left( 1 + O\left(\frac{1}{\zeta}\right) \right) \,.
\eeq
We fix a solution singular at $z=-1$ ($\zeta=-1$) by the requirement that
its asymptotic scattering phase differs from the regular solution
by a phase of $\frac{\pi}{2}$, giving an asymptotic expansion
\beq
v_{\mbox{\scriptsize sing}}(\zeta+1) & \ub{\sim}{\zeta \to \infty} &
 \frac{\cos{(2\sqrt{2\beta\zeta}+\gamma)}}{\zeta^{q+\frac{3}{4}}} \; \left( 1 + O\left(\frac{1}{\zeta}\right) \right)\,.
\eeq
To resume, looking for the nontrivial solutions of
the differential equation (\ref{ODE}) we have already found two
contributions: for a generic set of parameters $(\beta,j,q)$ the
solution regular at $z=1$ ($\zeta=0$) has a regular singularity at
$z=-1$ ($\zeta=-1$). This means that it is a linear combination of
$v_{\mbox{\scriptsize reg}}$ and $v_{\mbox{\scriptsize sing}}$, without loss
of generality:
\beq
u_{\mbox{\scriptsize reg}}(\zeta) & = & \cos{(\alpha)} \;
v_{\mbox{\scriptsize reg}}(\zeta+1) + \sin{(\alpha)} \; v_{\mbox{\scriptsize
sing}}(\zeta+1) \\
& \ub{\sim}{\zeta \to \infty} & \frac{1}{\zeta^{q+\frac{3}{4}}} \,
\sin{(2\sqrt{2\beta\zeta}+\underbrace{\gamma +  \alpha}_{=
\Delta_{j,q}(\beta)})} \; \left( 1 + O\left(\frac{1}{\zeta}\right) \right)\;.
\eeq
Because of the singularity at $z=-1$ ($\zeta=-1$) the
representation of $u_{\mbox{\scriptsize reg}}$ by a power series
expansion must break down for all values $|\zeta|>1$. Therefore, the
asymptotic scattering phase of $u_{\mbox{\scriptsize reg}}$ can in
general not be obtained from (\ref{FPS_reg}).

Only for those values $(\beta_n)_{n \in \mathbb{N}}$ that lead to a
minimal solution of the recurrence relation (\ref{REC_REL}) the
solution regular at $z=1$ is regular at $z=-1$ as well, and
therefore of class  $C^\infty(\mathbb{R})$. This means 
$\sin{\alpha} =0$ and therefore
\beqn
\label{phase}
\Delta_{j,q}(\beta_n) & = & \gamma \mod{\pi}\;.
\eeqn
Since we can always change the sign of the solution by an overall factor
the scattering phase is only determined up to an integer multiple of  $\pi$.
In Chapt.~\ref{monodromy} we will eventually determine $\gamma$.

If we had made an expansion around the regular singularity
at $z=-1$ a transformation according to $\chi = \frac{1}{2} (z+1)$ would
have been suitable. Substitution of 
\beq
A(z) = (z^2 -1 )^\frac{q}{2} \, v(\chi)
\eeq
in Eq.~(\ref{ODE}) then yields
\beqn
\label{IsoCHE} 
0 & = & v''(\chi) + \left( \frac{q+1}{\chi} + \frac{q+1}{\chi-1} \right) v'(\chi) + \left(
\frac{\mu - \beta}{\chi(\chi-1)} + \frac{ 2\beta}{\chi -1} \right) \; v(\chi)\,.
\eeqn
Eq.~(\ref{IsoCHE}) can also be understood as the isomonodromic deformation of the confluent Heun equation
with an additional apparent regular singularity moving towards $\chi=1$ \cite{Slav:1999}. By a theorem 
\cite[Sect.~2, Th.]{Slav:1999} Eq.~(\ref{IsoCHE}) is then related to a Painlev\'e equation, in our case $P^{VI}$: The 
Painlev\'e equation is the Newtonian equation of motion corresponding to the quantum Hamiltonian given by means of the 
Heun equation.

This is clear since the $SU(2)$-invariant self-dual metrics are specified exactly by a solution of
the Painlev\'e equation $P^{VI}$ \cite{Okum:2003}.

\subsection{Minimal solutions and continued fractions}
\label{minimal}
According to Pincherle's Theorem (see \cite[Sect. 5.3]{Jone:1980}) a
three-term recurrence relation
\beq
\forall_{n\ge1} \qquad y_{n+1} & = & - \delta_n \, y_{n} + \gamma_n \, y_{n-1}
\eeq
has a minimal solution if and only if the following continued fraction
converges:
\beq
\sum_{k=1}^\infty \frac{\gamma_k|}{|\delta_k} =
\frac{\gamma_1}{\displaystyle \delta_1 + \frac{\gamma_2}{\displaystyle
\delta_2 + \frac{\gamma_3}{\displaystyle \delta_3 + \dots}}} \;.
\eeq
In particular, if $(a_n)_{n\in\mathbb{N}}$ is the minimal solution it
follows that $a_0 \not = 0$ and
\beq
-\frac{a_1}{a_0} + \sum_{k=1}^\infty \frac{\gamma_k|}{|\delta_k} = 0 \;.
\eeq
For any set of parameters $(j,q)$, where as before $q\ge 0$, we now define the function $\widetilde{\mathcal{M}}(j,q
\,|x)$ by the continued fraction of {\em Thron} type (or T-fraction) 
\beq
\widetilde{\mathcal{M}}(j,q \,|x) & := & \delta_0 + \sum_{k=1}^\infty \frac{\gamma_k|}{|\delta_k}\,,\\
\mbox{where} \quad  \forall_{k\ge0} \quad \delta_k & := & \frac{k(k+2q+1)+\mu-x}{(k+1)(k+q+1)} \,,\\
           \forall_{k\ge1} \quad \gamma_k & := & \frac{2x}{(k+1)(k+q+1)} \;.
\eeq
Using the Umordnungssatz \cite[Kap.~6.42, Satz 2]{Perr:1977} and 
\cite[Sect. 7.3, Th. 7.23]{Jone:1980}, we can conclude that the 
function $\widetilde{\mathcal{M}}(j,q \,|x)$ is a meromorphic function on
$\mathbb{C}$. Moreover, by the above its real zeros determine the values of
$x=-\beta$ for which the recurrence relation (\ref{REC_REL}) has a minimal
solution, i.e. the radius of convergence of the power series expansion
(\ref{FPS_reg}) becomes infinite.

The corresponding solution $u_{\mbox{\scriptsize reg}}$ is then of class
$C^\infty(\mathbb{R})$. Since the relation $A(\beta,j,q|z) = A(-\beta,j,q|-z)$ 
holds for the differential equation (\ref{ODE}), it follows that any such 
smooth solution for $\beta,j,q$ in $z$ is simultaneously a smooth solution 
for $-\beta,j,q$ in the variable $-z$. 
Therefore, $x=\beta$ must be another zero of $\widetilde{\mathcal{M}}(j,q \,|x)$. This is
\beq
\widetilde{\mathcal{M}}(j,q\,|\beta)=0   \quad \Leftrightarrow \quad
\widetilde{\mathcal{M}}(j,q\,|-\beta)=0\;.
\eeq

Using the Umordnungssatz \cite[Kap.~6.42, Satz 2]{Perr:1977}, one can see that one can cancel
the factor $(k+1)(k+q+1)$ in $\gamma_k, \delta_k$. Namely, it is
equivalent to calculate the zeros of the meromorphic function 
$\mathcal{M}(j,q\,|x)$ instead of $\widetilde{\mathcal{M}}(j,q\,|x)$ where
$\mathcal{M}(j,q\,|x)$ is defined by the continued fraction 
\beqn
 \label{T-frac}
\mathcal{M}(j,q \,|x) & := & d_0 + \sum_{k=1}^\infty \frac{c_k|}{|d_k}\,,\\
\nonumber
\mbox{where} \quad  \forall_{k\ge0} \quad d_k & := & k(k+2q+1)+\mu-x \,,\\
\nonumber
           \forall_{k\ge1} \quad c_k & := & 2k(k+q)x \;.
\eeqn
For this continued fraction one can even prove separate convergence (see \cite{Thron:1991} for definitions).
To see this, recall that $\mu=(q-j)(j+q+1)$, and rewrite the continued fraction in the following way:
\beq
 \forall_{j-q\ge 1} \quad \mathcal{M}(j,q \,|x) & = & d_0 + \sum_{k=1}^{j-q-1} \frac{c_k|}{|d_k} 
 \begin{array}{l} \vspace*{0.1cm} \\ + \; \frac{\textstyle 2j(j-q)}  
 {\textstyle -1 +  \frac{m(j,q|x)}{x}} \end{array} \;,\\
 \forall_{j\ge 0} \quad \mathcal{M}(j,j \,|x) & = & x \left( -1 + \frac{m(j,j|x)}{x} \right) \;.
\eeq 
Here, $m(j,q\,|x) := \sum_{k=j-q+1}^\infty \frac{c_k|}{|d_k}$ can be written
as a T-fraction using again the Umordnungssatz \cite[Kap.~6.42, Satz 2]{Perr:1977}, i.e.
\beqn
 \label{T-frac_red}
 m(j,q\,|x) & = & \sum_{k=j-q+1}^{\infty} \frac{F_k\, x|}{|1 + G_k \, x} \,,\\
 \nonumber
\mbox{where} \quad  \forall_{k\ge j-q+1}  \quad   G_k & := & -\frac{1}{(k-j+q)(k+j+q+1)} \,,\\
 \nonumber
                    \forall_{k\ge j-q+2}  \quad   F_k & := & \frac{2k(k+q)}{(k-j+q-1)(k-j+q)(k+j+q)(k+j+q+1)} \,,\\
\nonumber
                                            F_{j-q+1} & := & j-q+1 \,. 
\eeqn
Since we have both $\sum_{k=j-q+1}^\infty |F_k| < \infty$ and $\sum_{k=j-q+1}^\infty |G_k| < \infty$ we can
apply \cite[Sect. 3, Th. 3.1]{Thron:1991} to the T-fraction (\ref{T-frac_red}): let $A_k(x)$ and $B_k(x)$ be 
the numerators and denominators, respectively, of the n$^{th}$ approximant of the T-fraction  (\ref{T-frac_red}). Then
the sequences $(A_k(x))_{k\ge j-q+1}$ and $(B_k(x))_{k\ge j-q+1}$ converge, uniformly on compact subsets of
$\mathbb{C}$, to entire functions $A(x)$ and $B(x)$ of order at most one. 
Further $B(0)=1$, $A(0)=0$, $A^\prime(0)=F_{j-q+1}$ so that neither function is identically zero, 
and $\frac{1}{x} m(j,q|x)$ is well defined at $x=0$.

\subsection{The determination of the scattering phase and the monodromy}
\label{monodromy}
In this section we calculate the scattering phase of the regular solution
$u_{\mbox{\scriptsize reg}}$ in the case that it is also regular at
$\zeta=-1$. 

Let us first look at the asymptotic expansion of the solutions of
Eq.~(\ref{ODE2}). For an asymptotic expansion a Fabry transformation
\cite{Olve:1974} is suitable, i.e.~a change of the variable according to $x^2=\zeta$. With
$u(\zeta)=U(x)$ Eq.~(\ref{ODE2}) becomes
\beq
(x^2+1) U''(x) + \left( - \frac{x^2+1}{x} + \frac{2(q+1)(2x^2+1)}{x}  \right) \; U'(x)
 + \big(4\beta (2x^2+1) + 4\mu \big) \;  U(x) = 0 \;.
\eeq
By standard methods (c.f. \cite{Codd:1955}) one finds that this equation has two linearly independent
solutions $H_1, H_2$, such that for some $\delta>0$ the following asymptotic expansions hold:
\beqn 
 \label{asymp}
 \begin{array}{crccllcrcl}
 \forall_x &  -2\pi + \delta & \le & \mbox{arg} \, x & \le &  \pi - \delta &: \quad &   H_1(x) 
 & \sim & \frac{1}{x^{2q+\frac{3}{2}}} e^{2i\sqrt{2\beta}x}
\sum\limits_{k=0}^\infty \frac{(-1)^k \, p_k}{(4i\sqrt{2\beta}x)^k} \,,\\ \vspace*{-0.4cm} \\
 \forall_x & -\pi + \delta & \le & \mbox{arg} \, x & \le &  2 \pi - \delta &: \quad &   H_2(x) 
 & \sim & \frac{1}{x^{2q+\frac{3}{2}}} e^{-2i\sqrt{2\beta}x}
\sum\limits_{k=0}^\infty \frac{p_k}{(4i\sqrt{2\beta}x)^k} \,,
\end{array}
\eeqn
where the coefficients $(p_k)_{k \in \mathbb{N}}$ fulfill the following four-term
recurrence relation
\beqn
 \nonumber
 p_{-2} & = & p_{-1} =  0 \;,\\
 \nonumber
 \forall k \ge 0: \qquad p_{k+1} 
 & = & \frac{k(k+1)-4\beta-4\,j(j+1)-\frac{3}{4}}{k+1} \, p_k
 + 32\beta\,\frac{k+q}{k+1} \, p_{k-1} \\
 \label{REC_REL2}
 & - & 32\beta \, \frac{k^2+(2q-1)k-q+\frac{1}{4}}{k+1} \,
p_{k-2} \;.
\eeqn
It follows from this recurrence relation that for $p_0 \in \mathbb{R}$ the 
coefficients $(p_k)_{k \in \mathbb{N}}$ are real. In particular, we have 
$H_1(x) = \overline{H_2(\bar{x})}$. This property remains true not only asymptotically, 
but also for the actual solutions: Since all the coefficients in the 
differential equation are real, the complex conjugate of any solution is again a solution.

The Stokes phenomenon will lead to a non-trivial monodromy of the solutions if we pass 
through a cycle around infinity $\zeta \mapsto \zeta e^{2\pi i}$. This is because for $x \mapsto x e^{\pi i}$ we
will leave the sector in which the asymptotic expansion (\ref{REC_REL2}) for $H_2$ holds \cite{Olve:1995}. However, since 
the solutions $H_1, H_2$ are a complete system of solutions the new solutions can be expressed as linear 
combinations of them. These are the connection formulae. In the case of a second order
differential equation they had been explicitly calculated, namely in \cite{Olve:1994}.
We obtain for $|x| \gg 0$
\beqn
\label{con}
 \lim_{\theta \rightarrow 1} H_1(x e^{\pi i \theta}) & = & - 2\pi P \, H_1(x) + i   \, H_2(x) \,,\\
\nonumber
 \lim_{\theta \rightarrow 1} H_2(x e^{\pi i \theta}) & = & i  \, H_1(x) \,.
\eeqn
The parameter $P$ can be determined by a generalized
Borel transformation of the asymptotic solution as pointed out
in \cite{Gura:1994a}, \cite{Gura:1994b}. For ordinary
differential equations of second order this has been done explicitly
by Daalhuis and Olver \cite{Olve:1994}. Applying these results we can
determine the parameter $P$ from the coefficients in the
asymptotic expansion (\ref{REC_REL2}), i.e.
\beqn
 \label{P_formula}
 P = \lim_{k \to \infty} \frac{1}{(k-1)!} \;
\frac{p_k}{p_0} \;.
\eeqn
In particular, the parameter $P$ is real. We subsume the solutions $H_1, H_2$ in 
the vector $\overrightarrow{H}:= \,^t(H_1, H_2 )$. 

We need to calculate the linear combination of the regular solution $u_{\mbox{\scriptsize reg} }(\zeta)$ 
in terms of the solutions $\overrightarrow{H}(x)$ derived earlier. This amounts to determining a complex 
parameter $\lambda \in U(1)$ such that
\beq
 u_{\mbox{\scriptsize reg} }(\zeta) & = & \left( \; \lambda  \quad  \overline{\lambda} \; \right) \cdot \overrightarrow{H}(x) \;.
\eeq
The occurrence of $\lambda$ and $\overline{\lambda}$ is due to the fact that $u_{\mbox{\scriptsize reg} }(\zeta)$
is real, i.e. $\overline{u_{\mbox{\scriptsize reg} }(\overline{\zeta})}= u_{\mbox{\scriptsize reg} }(\zeta)$, and
$\overline{H_1(\overline{x})} = H_2(x)$.

Now, the representation of $u_{\mbox{\scriptsize reg} }(\zeta)$ by the power series (\ref{FPS_reg}) is valid for $|\zeta|>1$
iff the coefficients $(a_k)_{k \ge 1}$ constitute a minimal solution of the recurrence relation (\ref{REC_REL}). Therefore,
$u_{\mbox{\scriptsize reg} }(\zeta)$ has a trivial monodromy around both of the regular singularities iff the 
$(a_k)_{k \ge 1}$ constitute a minimal solution of Eq.~(\ref{REC_REL}). Hence, iff for any set of parameters $(\beta,j,q)$ the
equation $\mathcal{M}(j,q|\beta)=0$ holds, then the following equation must hold: 
\beq
 u_{\mbox{\scriptsize reg} }(\zeta) 
 & = & \lim_{\theta \rightarrow 1}  u_{\mbox{\scriptsize reg} }(\zeta e^{2\pi i \theta} ) 
 = \left( \; \lambda  \quad \overline{\lambda} \; \right) \cdot  \lim_{\theta \rightarrow 1} \overrightarrow{H}(x e^{\pi i \theta}) \\
 & = & \left( \; \lambda  \quad \overline{\lambda} \; \right) \cdot \left(\begin{array}{cc} -2\pi P & i
\\ i  & 0 \end{array} \right)  \cdot \overrightarrow{H}(x) \;.
\eeq
Thus, the equation $\mathcal{M}(j,q|\beta)=0$ holds iff  
\beq
 \left( \; \lambda  \quad \overline{\lambda} \; \right) = \left( \; \lambda  \quad \overline{\lambda} \; \right) \cdot  \left(\begin{array}{cc} -2\pi P & i
\\ i  & 0 \end{array} \right)\;.
\eeq
If we set $\lambda = i \, e^{-i\gamma}$ with $\gamma \in \mathbb{R}$ we obtain (since $P$ is real)
\beqn
 \label{monodromy_result}
 \gamma = \frac{3\pi}{4} \mod{\pi} \,, & \quad & P=-\frac{1}{\pi} \;.
\eeqn
Thus, $\lambda$ encodes the crucial information for the 
asymptotic expansion of the regular solution, i.e.
\beq
 u_{\mbox{\scriptsize reg}}(\zeta) & =  & \frac{1}{i} \; 
 \Big( e^{i\gamma} H_2(x) - e^{-i\gamma} H_1(x) \Big) \\
 & \sim & \frac{2\,p_0}{\zeta^{q+\frac{3}{4}}} \Big( \sin{(2\sqrt{2\beta\zeta}+\gamma)} +
O\left(\frac{1}{\sqrt{\zeta}}\right) \Big)\;.
\eeq 
This is the desired formula for the scattering phase in the case that the set of parameters 
$(\beta,j,q)$ induces a minimal solution of the recurrence relation (\ref{REC_REL}).

\subsection{Computation}
\label{computation}
The numerical determination of the scattering phase consists of two steps:
First, we have to determine the successive zeros of the function $\mathcal{M}(j,q \,|x)$. 
However, one should mention that for the necessary evaluation of the continued fraction one has to use 
a backward algorithm since any forward algorithm must be numerically instable as shown in \cite{Gaut:1967}.

We have determined the positive zeros $(\beta_n)_{n \in \mathbb{N}}$ of $\mathcal{M}(j,q \,|x)$ by a 
simple bisection algorithm where the evaluation of the continued fraction was accomplished by the Gautschi 
algorithm. It is essential to use the function $\mathcal{M}(j,q\,|x)$ instead of the earlier defined 
$\widetilde{\mathcal{M}}(j,q\,|x)$: Because of the structure of the coefficients in the continued fraction 
the typical values of $\widetilde{\mathcal{M}}(j,q\,|x)$ become very small, and the determination of its 
zeros unstable. With $(\beta_n,j,q,\gamma)$ we then have the data needed for a numerical interpolation 
of the scattering phase over the whole range of $\beta$.

\section{Numerical Results}
\label{results}

\subsection{Numerical Results for the WKB approximation}
\label{WKB_res}

The figures below show the scattering phases for different quantum numbers
$j$ and $q$ with varying parameter $\beta$ obtained by the WKB approximation
from Sect.~\ref{WKB}. 
\begin{figure}[ht]
 \begin{minipage}[t]{0.5\linewidth}
  \centering
  \includegraphics[width=4cm,angle=-90]{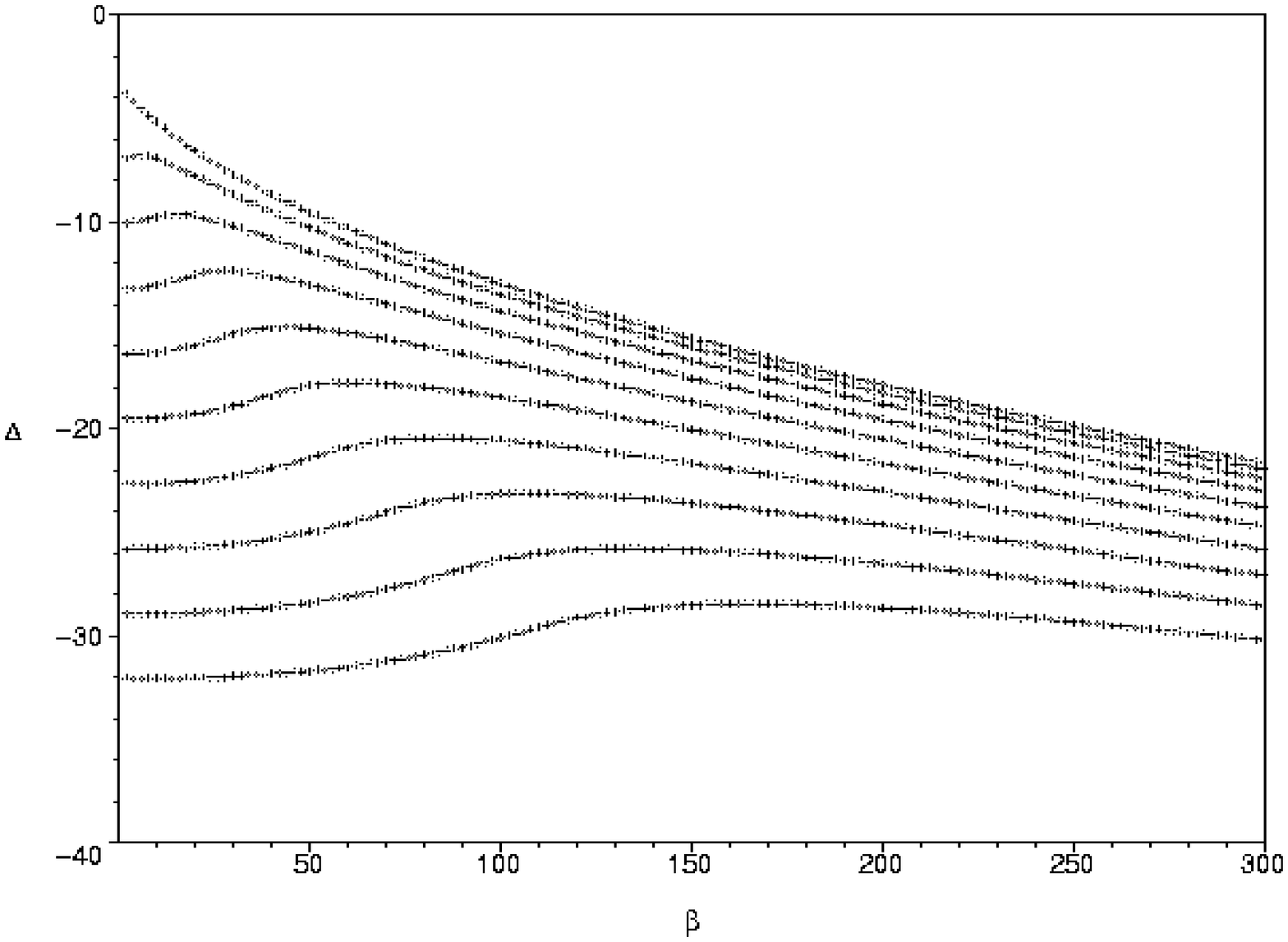}
  \caption{$\Delta_{j,q}$ for $j=1,\dots,10$, $q=1$}
  \label{WKB12345678910_1}
 \end{minipage}%
 \begin{minipage}[t]{0.5\linewidth}
  \centering
  \includegraphics[width=4cm,angle=-90]{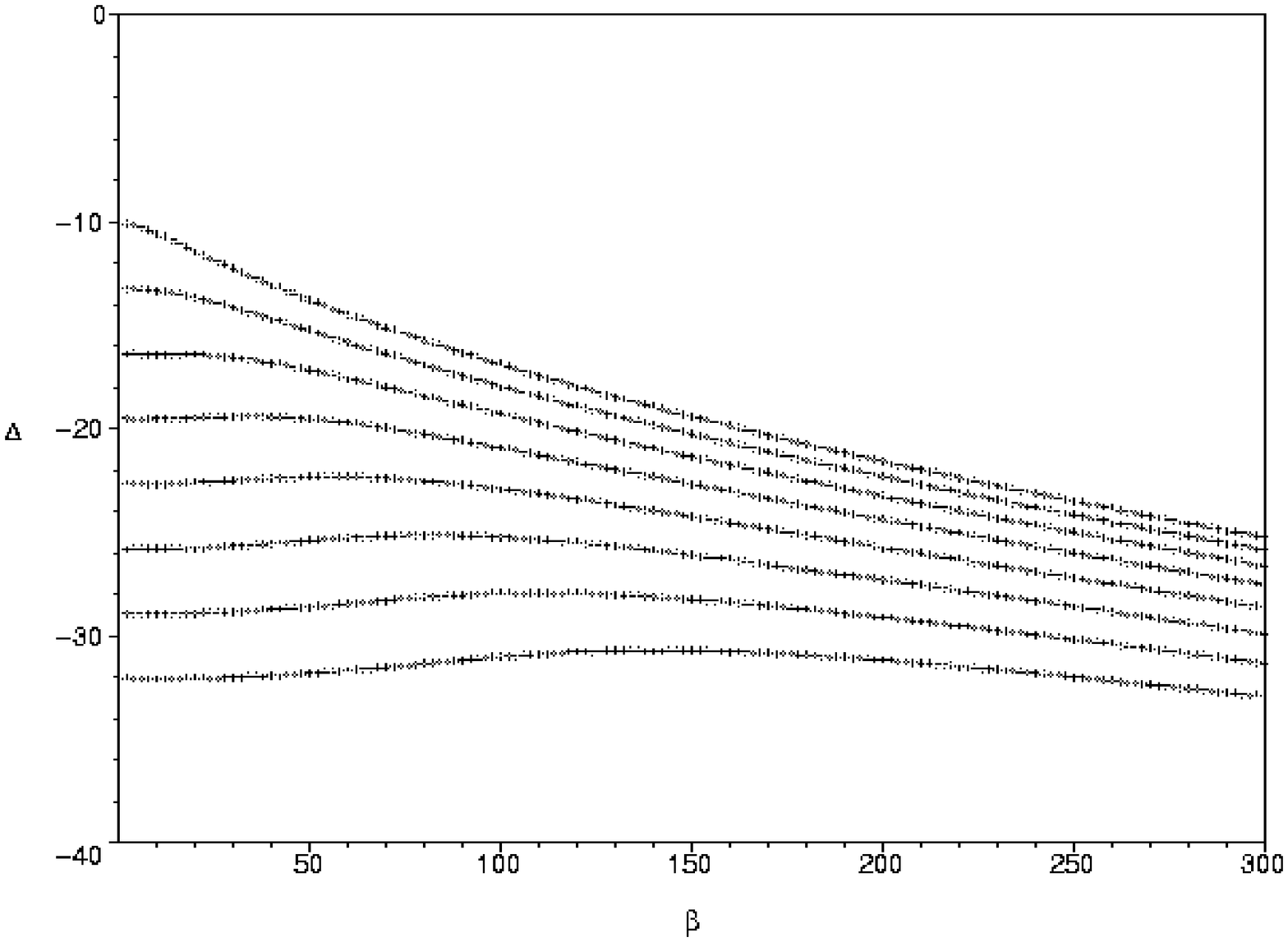}
  \caption{$\Delta_{j,q}$ for $j=3,\dots,10$, $q=3$}
  \label{WKB345678910_3}
 \end{minipage}\\[5pt]
%\end{figure}
%\begin{figure}[ht]
 \begin{minipage}[t]{0.5\linewidth}
  \centering
  \includegraphics[width=4cm,angle=-90]{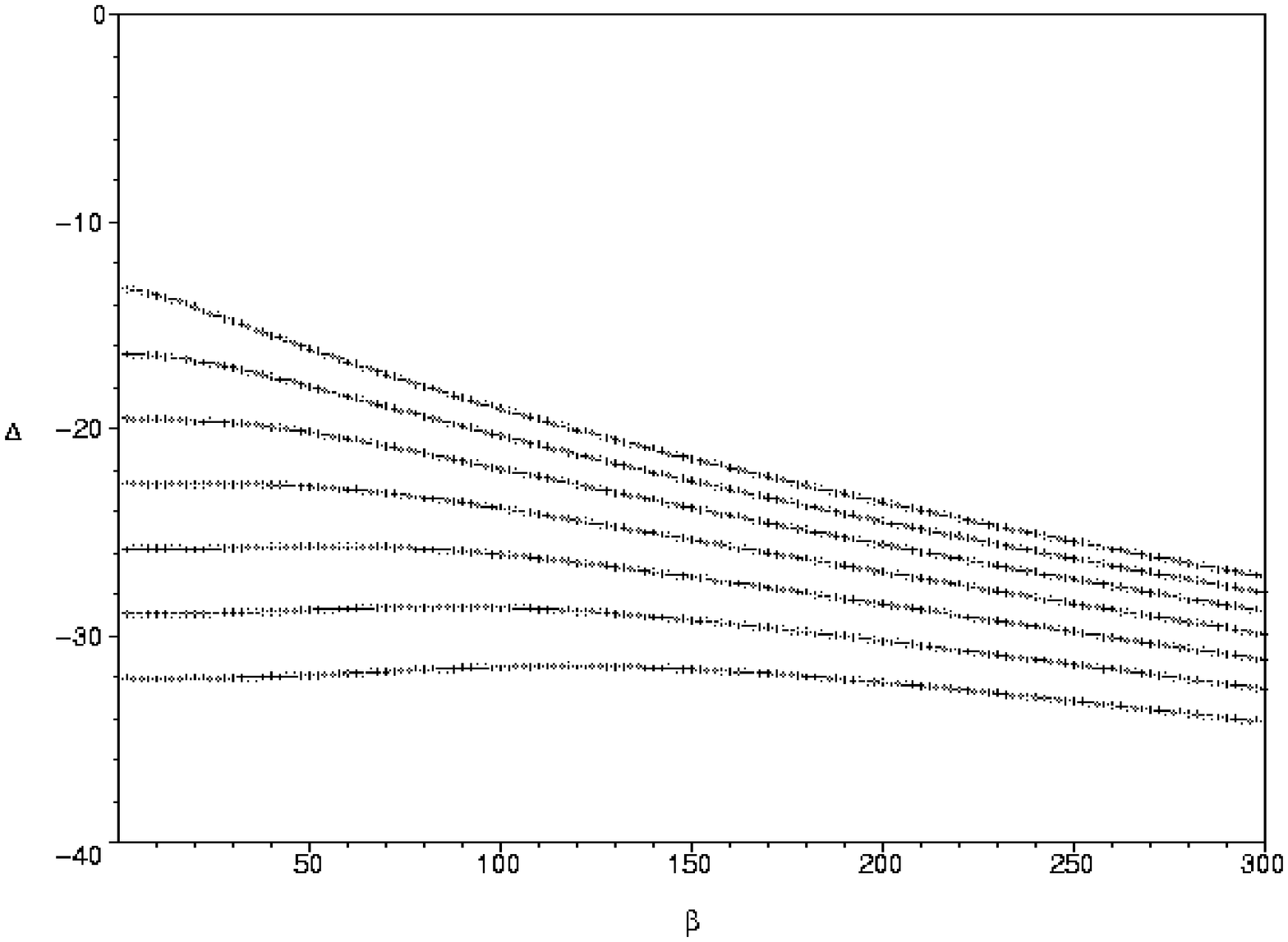}
  \caption{$\Delta_{j,q}$ for $j=4,\dots,10$, $q=4$}
  \label{WKB45678910_4}
 \end{minipage}%
 \begin{minipage}[t]{0.5\linewidth}
  \centering
  \includegraphics[width=4cm,angle=-90]{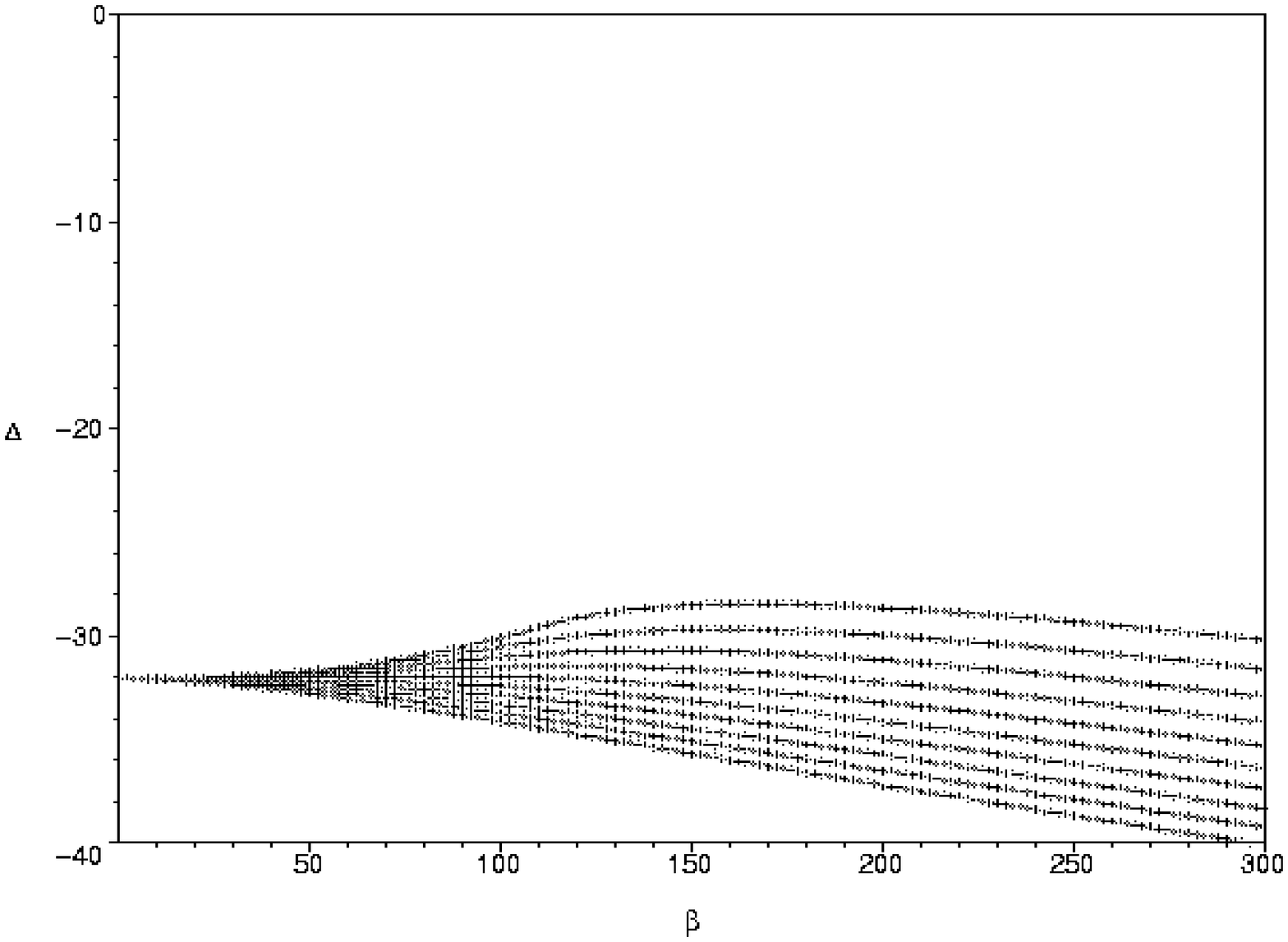}
  \caption{$\Delta_{j,q}$ for $j=10$, $q=1,\dots,10$}
  \label{WKB10_12345678910}
 \end{minipage} \\ [5pt]
%\end{figure}
%\begin{figure}[hb]
 \begin{minipage}[t]{0.5\linewidth}
  \centering
  \includegraphics[width=4cm,angle=-90]{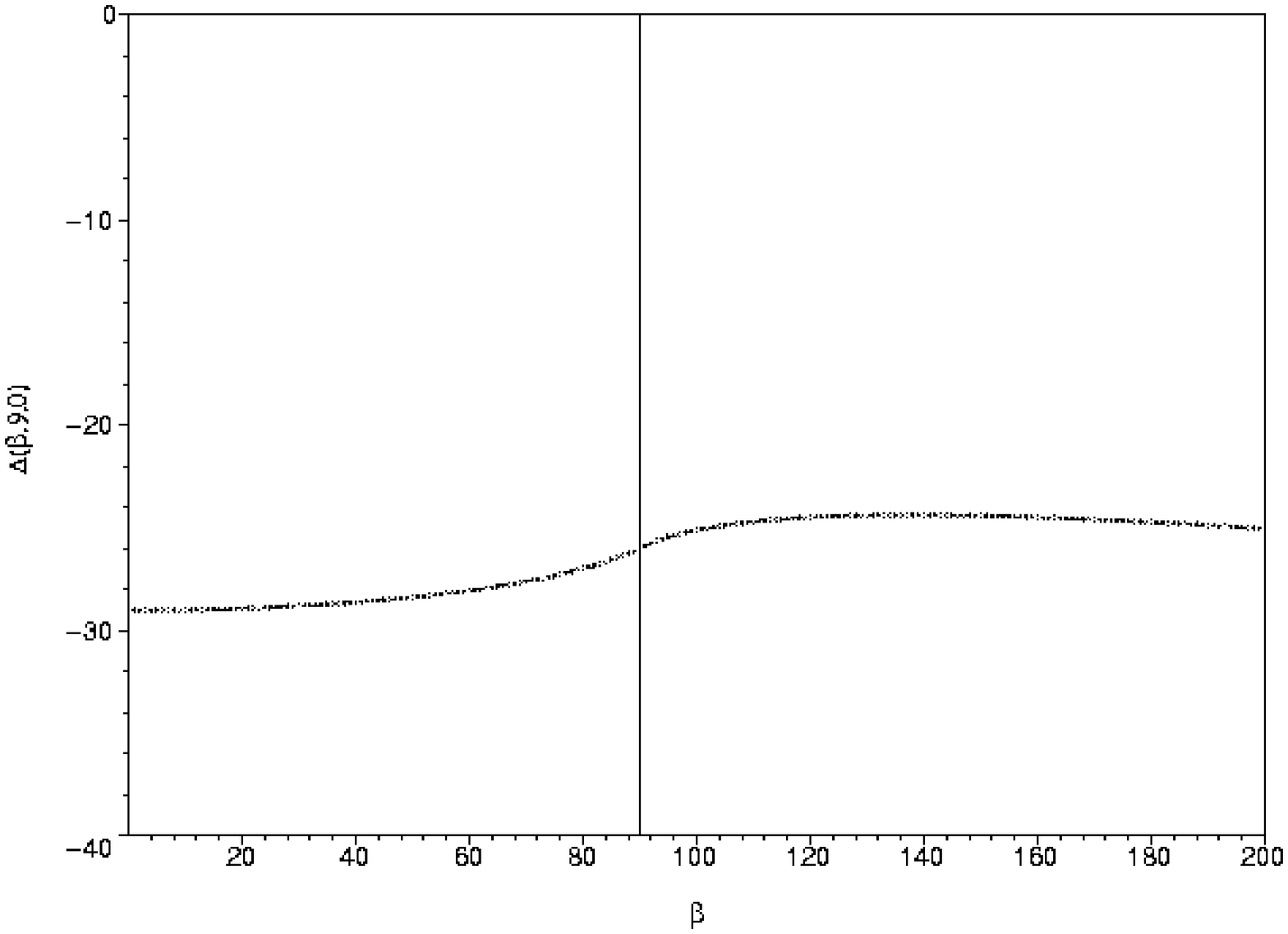}
  \caption{$\Delta_{9,0}$ for case II, III, IV}
  \label{WKB9_0}
 \end{minipage}%
 \begin{minipage}[t]{0.5\linewidth}
  \centering
  \includegraphics[width=4cm,angle=-90]{WKB9_0.ps}
  \caption{$\Delta_{9,0}$ for case I}
  \label{WKB9_0int}
 \end{minipage} \\[5pt]
 \begin{minipage}[t]{0.5\linewidth}
  \centering
  \includegraphics[width=4cm,angle=-90]{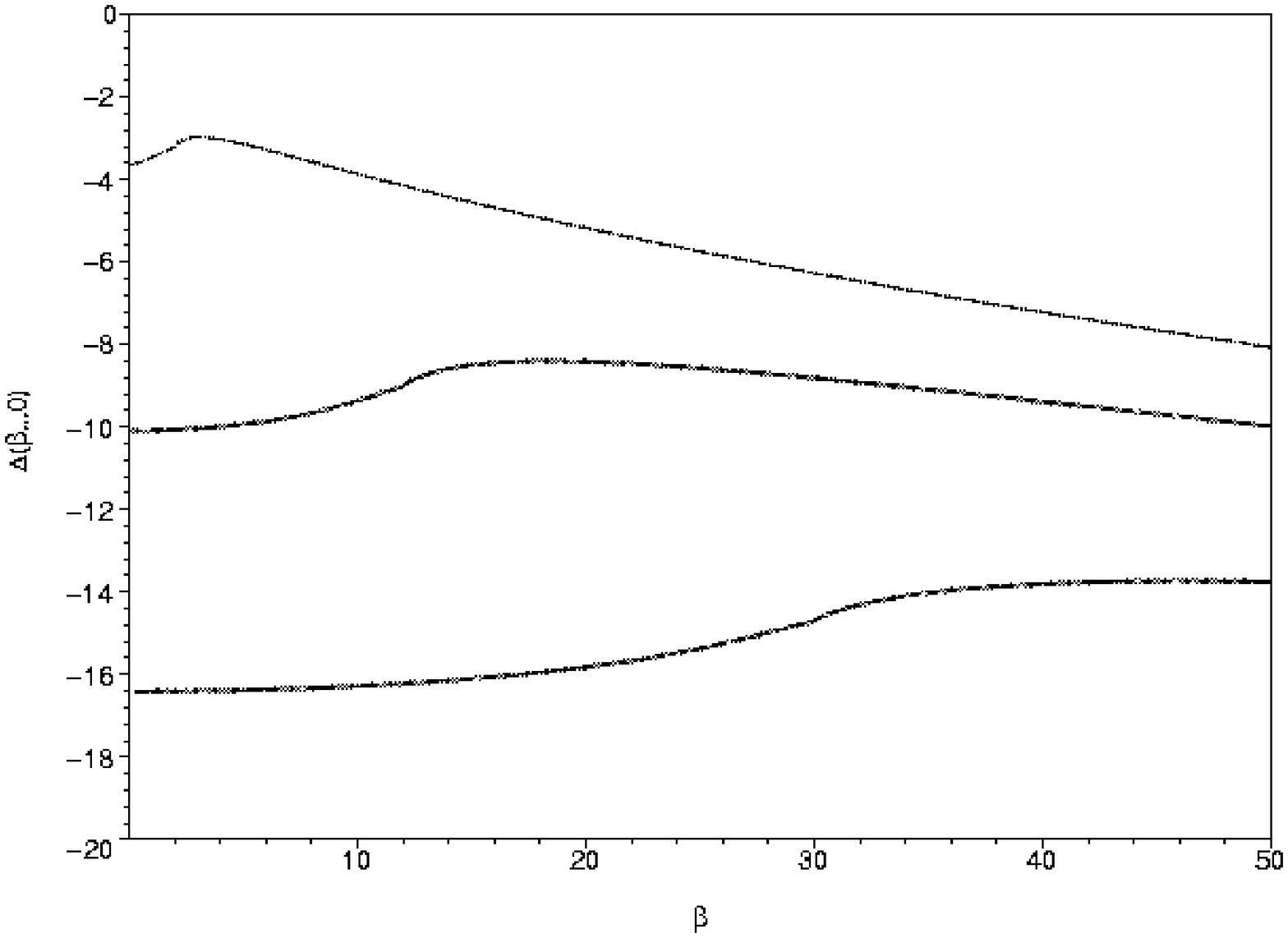}
  \caption{$\Delta_{j,0}$ for $j=1,3,5$}
  \label{WKB135_0}
 \end{minipage}%
 \begin{minipage}[t]{0.5\linewidth}
  \centering
  \includegraphics[width=4cm,angle=-90]{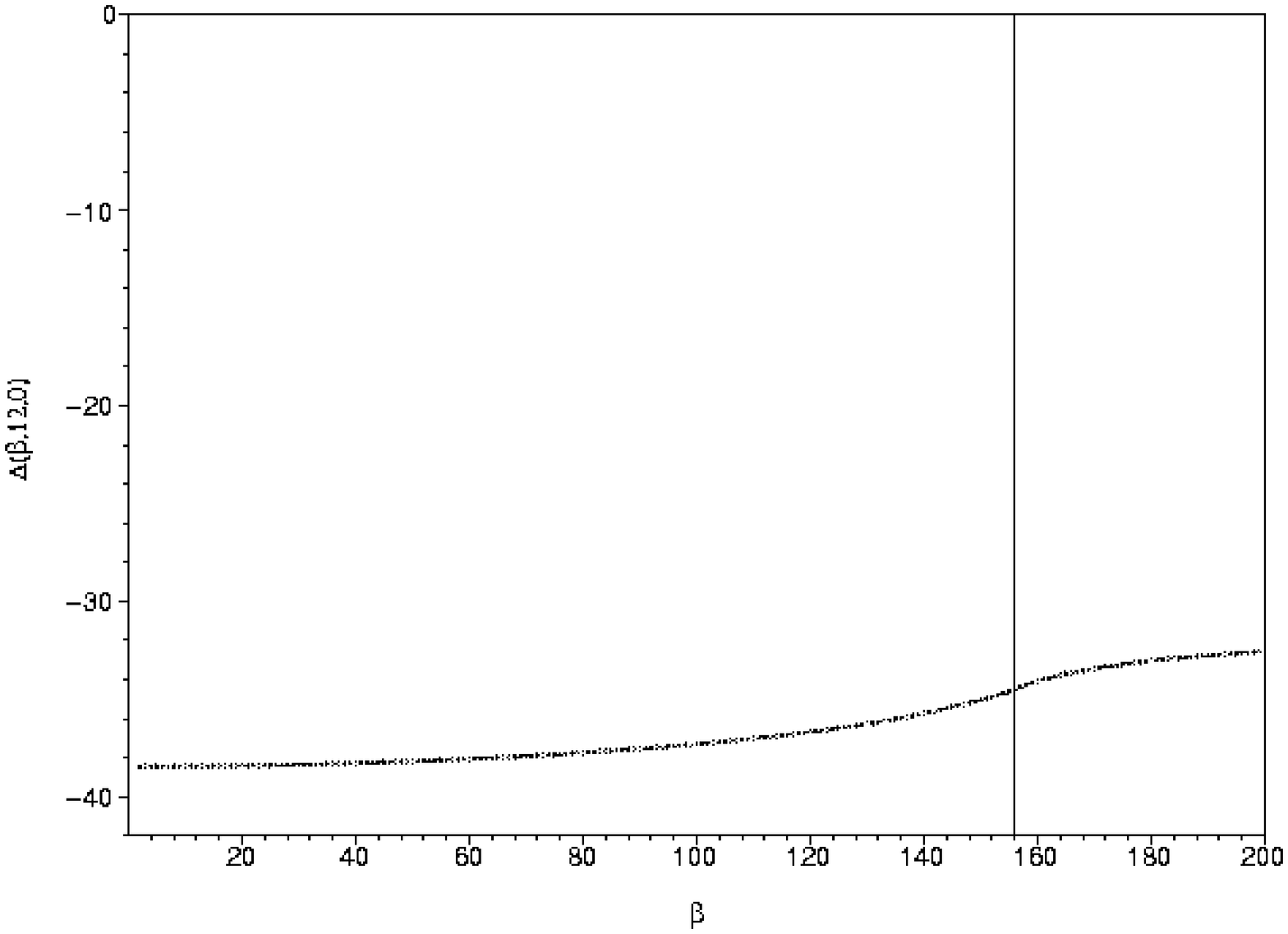}
  \caption{$\Delta_{12,0}$}
  \label{WKB12_0}
 \end{minipage}
\end{figure}

In Fig.~\ref{WKB12345678910_1} to Fig.~\ref{WKB10_12345678910} 
we have applied the WKB approximation of case~I. 

Fig.~\ref{WKB12345678910_1} shows (from the top to the bottom) the graphs for 
$\Delta_{j,q}$ for $j=1,\dots,10$ and $q=1$.

Fig.~\ref{WKB345678910_3} shows (from the top to the bottom) the graphs for 
$\Delta_{j,q}$ for $j=3,\dots,10$ and $q=3$.

Fig.~\ref{WKB45678910_4} shows (from the top to the bottom) the graphs for 
$\Delta_{j,q}$ for $j=4,\dots,10$ and $q=4$.

Fig.~\ref{WKB10_12345678910} shows (from the top to the bottom) the graphs for 
$\Delta_{j,q}$ for $j=5$ and $q=1,3,5$.

In Fig.~\ref{WKB9_0} and Fig.~\ref{WKB9_0int} we have compared the WKB approximation of case~II to IV
(left) with the WKB approximation of case I (right) for the values $j=9, q=0$. The vertical line 
in Fig.~\ref{WKB9_0} is indicating the transition from case~II to IV.

Fig.~\ref{WKB135_0} shows (from the top to the bottom) the graphs for 
$\Delta_{j,0}$ for $j=1,3,5$.

Fig.~\ref{WKB12_0} shows the graph for $\Delta_{12,0}$.

\subsection{Comparison between WKB and Frobenius method}

In this section we compare the numerical results for the asymptotic scattering phase
obtained by the WKB approximation (already shown in Sect.~\ref{WKB_res}) with the
numerical data obtained by the Frobenius/continued-fraction method as described in Sect.~\ref{computation}.

The figures below show the scattering phases for different quantum numbers $j$ and $q$
with varying parameter $\beta$. Each figure will show the graph $(\beta, \Delta_{j,q}^{\mbox{\scriptsize WKB}}(\beta))$
of the corresponding WKB approximation, and crosses will mark the points 
$(\beta_n, \Delta^{\mbox{\scriptsize Fr}}_{j,q}(\beta_n))_{n \ge 1}$ obtained by the Frobenius/continued-fraction method
where $\beta_n$ is the n$^{th}$ zero of $\mathcal{M}(j,q|\beta)$ and 
\beqn
 \label{data}
 \Delta^{\mbox{\scriptsize Fr}}_{j,q}(\beta_n) = - \frac{\pi}{4} - l_{j,q}(n) \, \pi - \max{\left(0,j-1\right)} \, \pi 
\eeqn
with $\forall_{n\ge 1}: l_{j,q}(n) \in \mathbb{N}_0$. 
Notice that the scattering phase is only determined up to a multiple of $\pi$ since we can always 
change the sign of a solution by an overall factor. Therefore, Eq.~(\ref{data}) is equivalent to 
Eq.~(\ref{monodromy_result}). We will also present an additional and more significant 
diagram of $\Big( \Delta_{j,q}(\beta) \mod \pi \Big)$ and varying parameter $\beta$.

\medskip

\begin{figure}[ht]
 \begin{minipage}[t]{0.5\linewidth}
  \centering
  \includegraphics[width=4cm,angle=-90]{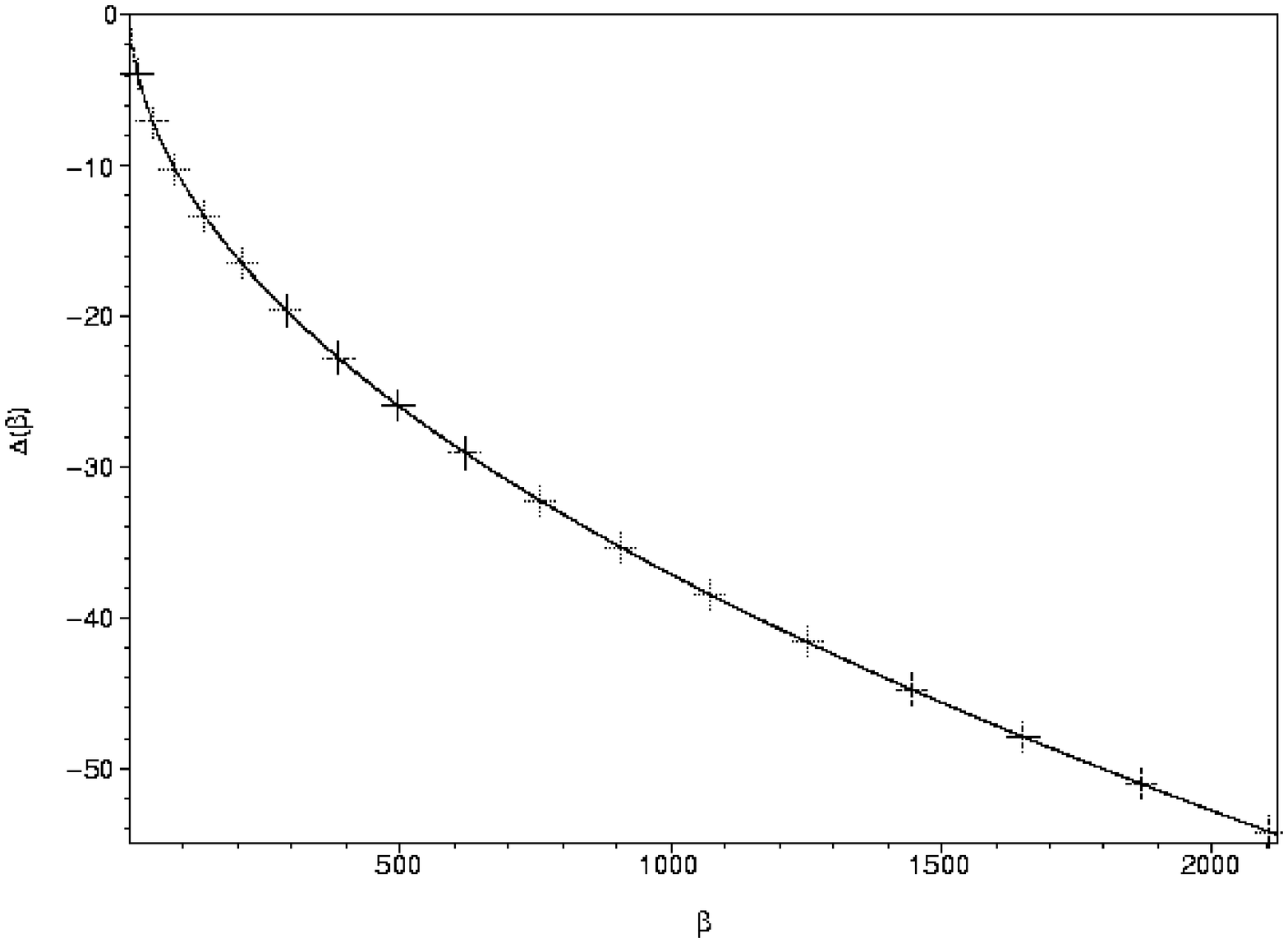}
  \caption{$\Delta^{\mbox{\scriptsize WKB}}_{0,0}$ and $\Delta^{\mbox{\scriptsize Fr}}_{0,0}$}
  \label{WKB_CF_0_0}
 \end{minipage}%
 \begin{minipage}[t]{0.5\linewidth}
  \centering
  \includegraphics[width=4cm,angle=-90]{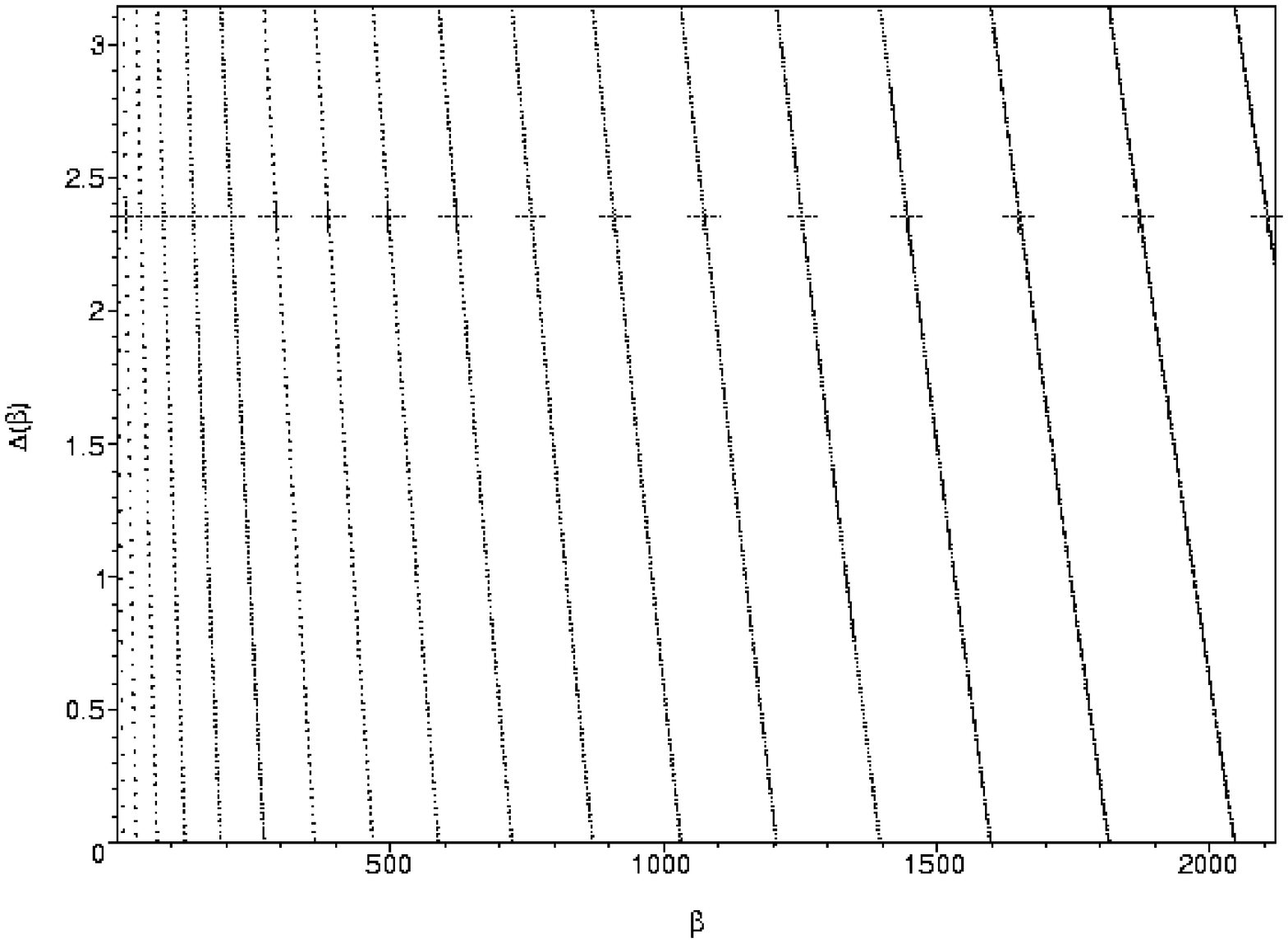}
  \caption{$\Delta^{\mbox{\scriptsize WKB}}_{0,0}$ and $\Delta^{\mbox{\scriptsize Fr}}_{0,0}$  modulo $\pi$}
  \label{WKB_CF_0_0modPi}
 \end{minipage}\\[5pt]
%\end{figure}
%\begin{figure}[ht]
 \begin{minipage}[t]{0.5\linewidth}
  \centering
  \includegraphics[width=4cm,angle=-90]{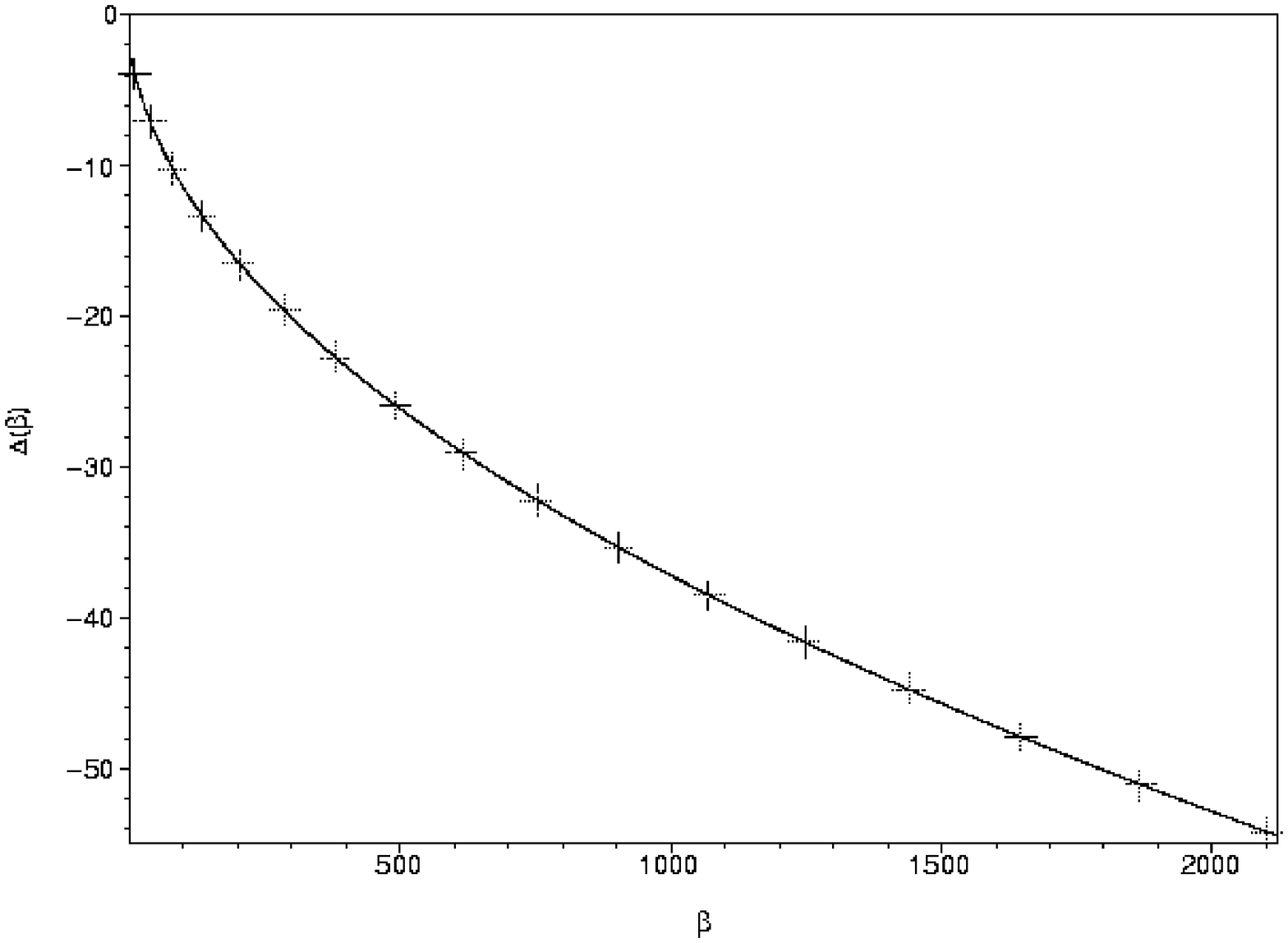}
  \caption{$\Delta^{\mbox{\scriptsize WKB}}_{1,0}$ and $\Delta^{\mbox{\scriptsize Fr}}_{1,0}$}
  \label{WKB_CF_1_0}
 \end{minipage}%
 \begin{minipage}[t]{0.5\linewidth}
  \centering
  \includegraphics[width=4cm,angle=-90]{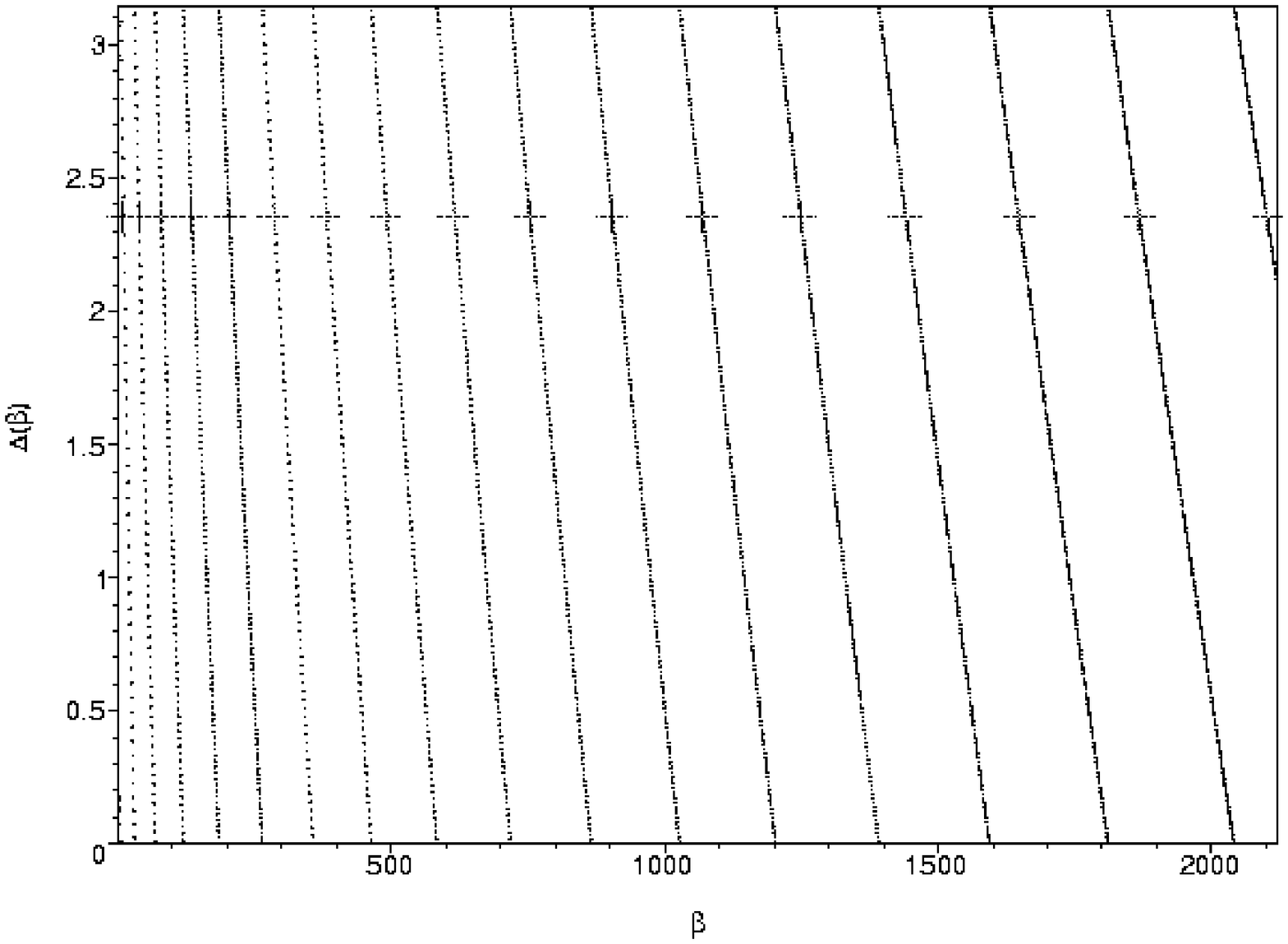}
  \caption{$\Delta^{\mbox{\scriptsize WKB}}_{1,0}$ and $\Delta^{\mbox{\scriptsize Fr}}_{1,0}$ modulo $\pi$}
  \label{WKB_CF_1_0modPi}
 \end{minipage} \\ [5pt]
%\end{figure}
%\begin{figure}[hb]
 \begin{minipage}[t]{0.5\linewidth}
  \centering
  \includegraphics[width=4cm,angle=-90]{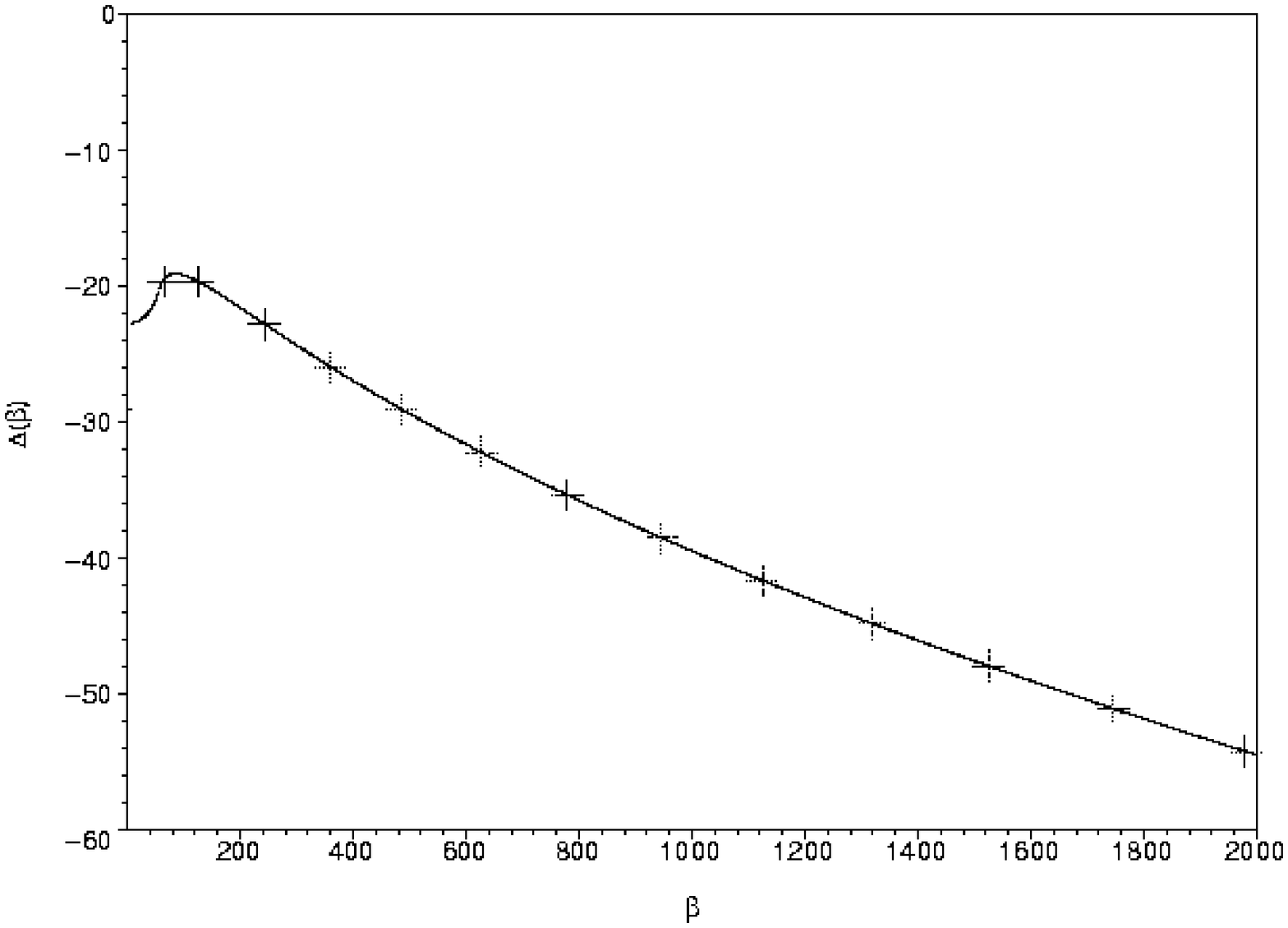}
  \caption{$\Delta^{\mbox{\scriptsize WKB}}_{7,0}$ and $\Delta^{\mbox{\scriptsize Fr}}_{7,0}$}
  \label{WKB_CF_7_0}
 \end{minipage}%
 \begin{minipage}[t]{0.5\linewidth}
  \centering
  \includegraphics[width=4cm,angle=-90]{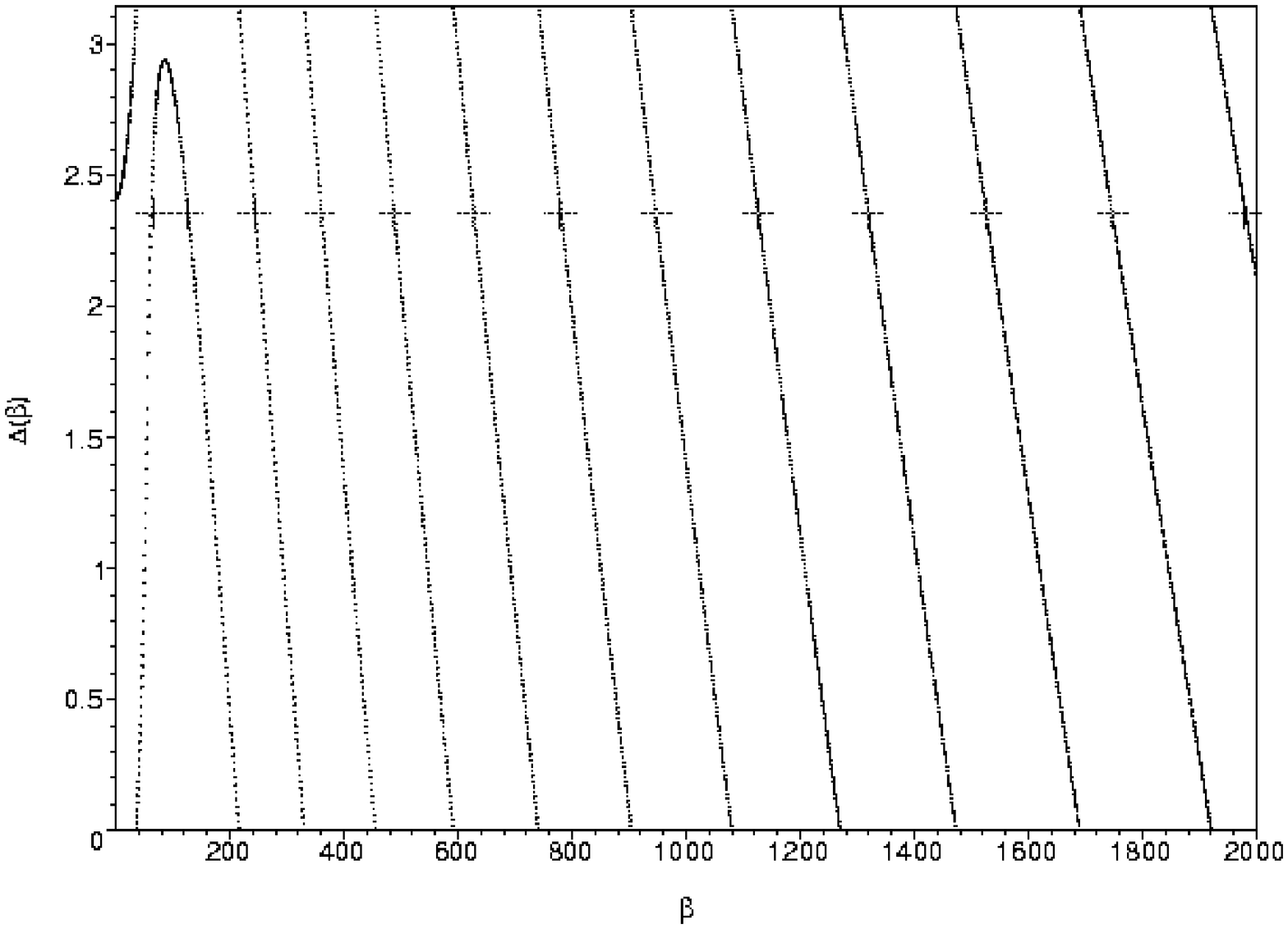}
  \caption{$\Delta^{\mbox{\scriptsize WKB}}_{7,0}$ and $\Delta^{\mbox{\scriptsize Fr}}_{7,0}$ modulo $\pi$}
  \label{WKB_CF_7_0modPi}
 \end{minipage} \\[5pt]
 \begin{minipage}[t]{0.5\linewidth}
  \centering
  \includegraphics[width=4cm,angle=-90]{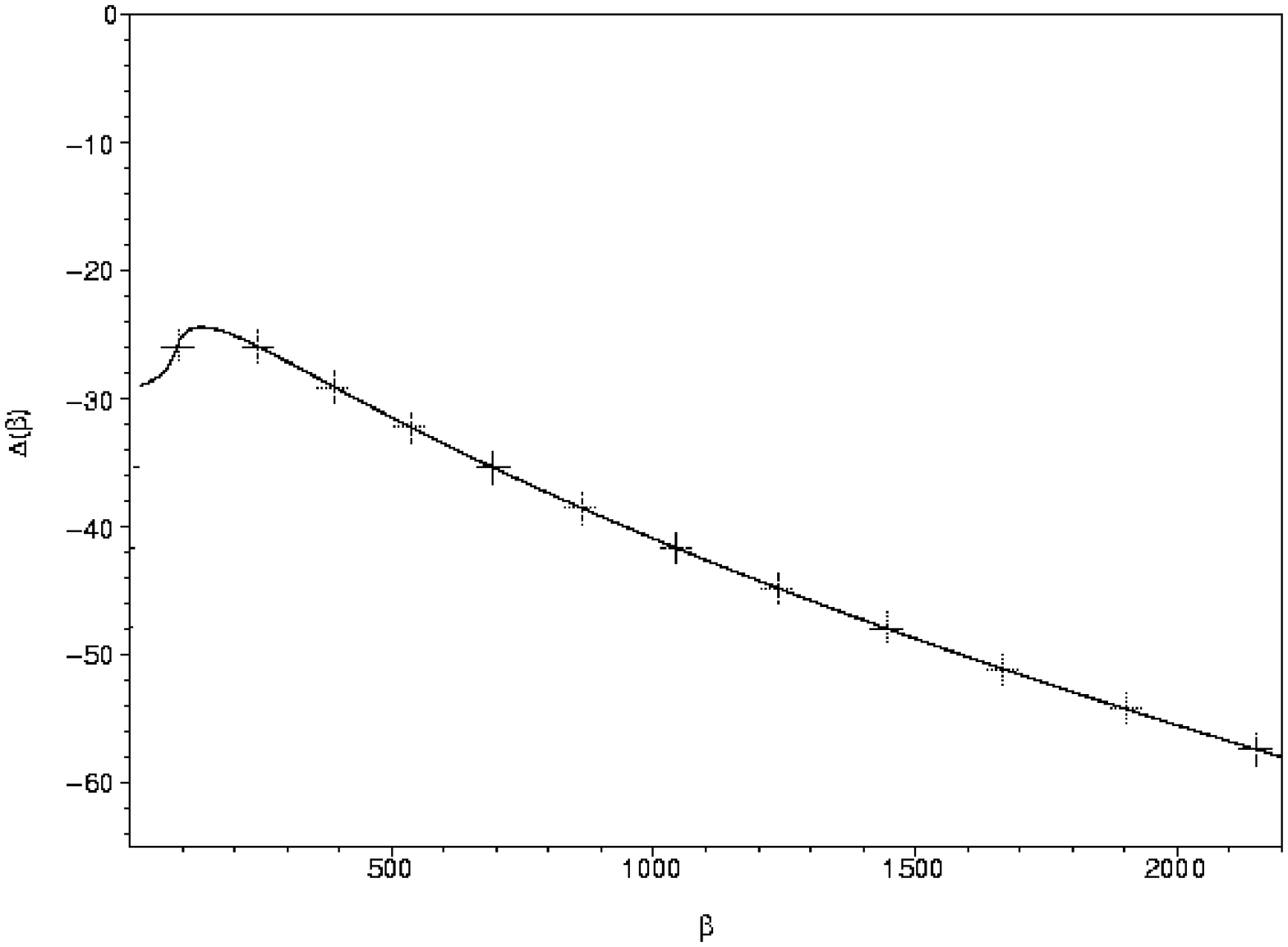}
  \caption{$\Delta^{\mbox{\scriptsize WKB}}_{9,0}$ and $\Delta^{\mbox{\scriptsize Fr}}_{9,0}$}
  \label{WKB_CF_9_0}
 \end{minipage}%
 \begin{minipage}[t]{0.5\linewidth}
  \centering
  \includegraphics[width=4cm,angle=-90]{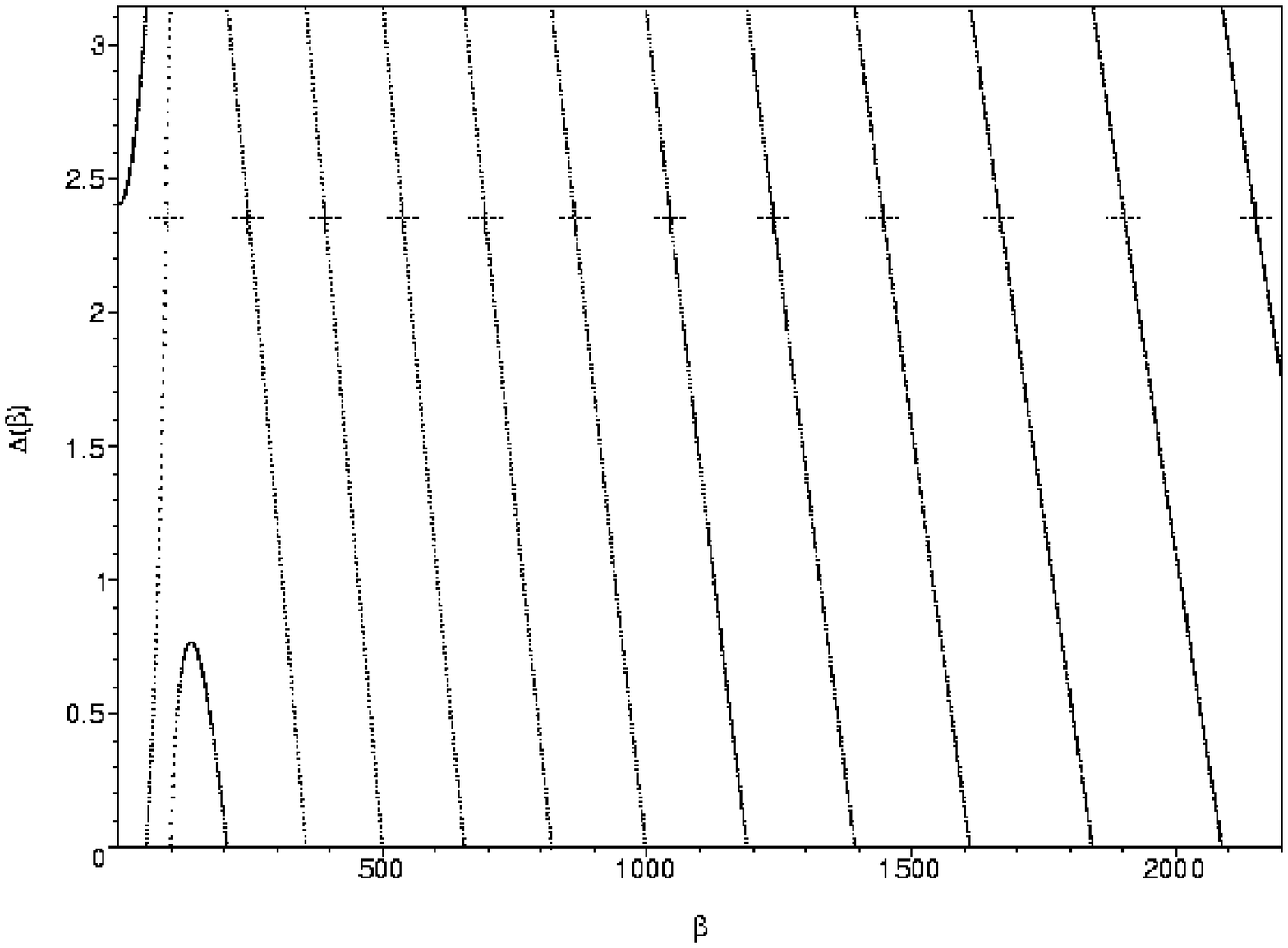}
  \caption{$\Delta^{\mbox{\scriptsize WKB}}_{9,0}$ and $\Delta^{\mbox{\scriptsize Fr}}_{9,0}$ modulo $\pi$}
  \label{WKB_CF_9_0modPi}
 \end{minipage}
\end{figure}

Figs.~\ref{WKB_CF_0_0}, \ref{WKB_CF_1_0}, \ref{WKB_CF_7_0}, \ref{WKB_CF_9_0} show 
$\Delta_{j,0}(\beta)$ with varying $\beta$. For $j \le 6$ it is $l(n):=n$ and all $\beta_n>j(j+1)$, i.e. $1>a\ge 0$ and we have to 
compare with case IV of the WKB approximation. $j=7$ is the lowest value where for $q=0$ we have 
$\beta_1 < j(j+1)$, i.e. $1<a$, and we have to compare with case II of the WKB approximation: 
we have to relabel the $\beta_n$. This can be done by setting $\forall_{j \ge 7}:\; l(n):=|n-\frac{3}{2}|-\frac{1}{2}$.

However, since the phase can be obtained only up to multiples of $\pi$ the function $l(n)$ is 
irrelevant for the comparison of $\Delta_{j,q}^{\mbox{\scriptsize WKB}}(\beta)$ and $\Delta^{\mbox{\scriptsize Fr}}_{j,0}(\beta_n)$.
For this purpose one can look at Figs.~\ref{WKB_CF_0_0modPi}, \ref{WKB_CF_1_0modPi}, \ref{WKB_CF_7_0modPi}, 
\ref{WKB_CF_9_0modPi}, which are independent of $l(n)$. These figures show an exact match
of the data sets obtained by the WKB approximation and the 
Frobenius/continued-fraction method.

\medskip

\begin{figure}[ht]
 \begin{minipage}[t]{0.5\linewidth}
  \centering
  \includegraphics[width=4cm,angle=-90]{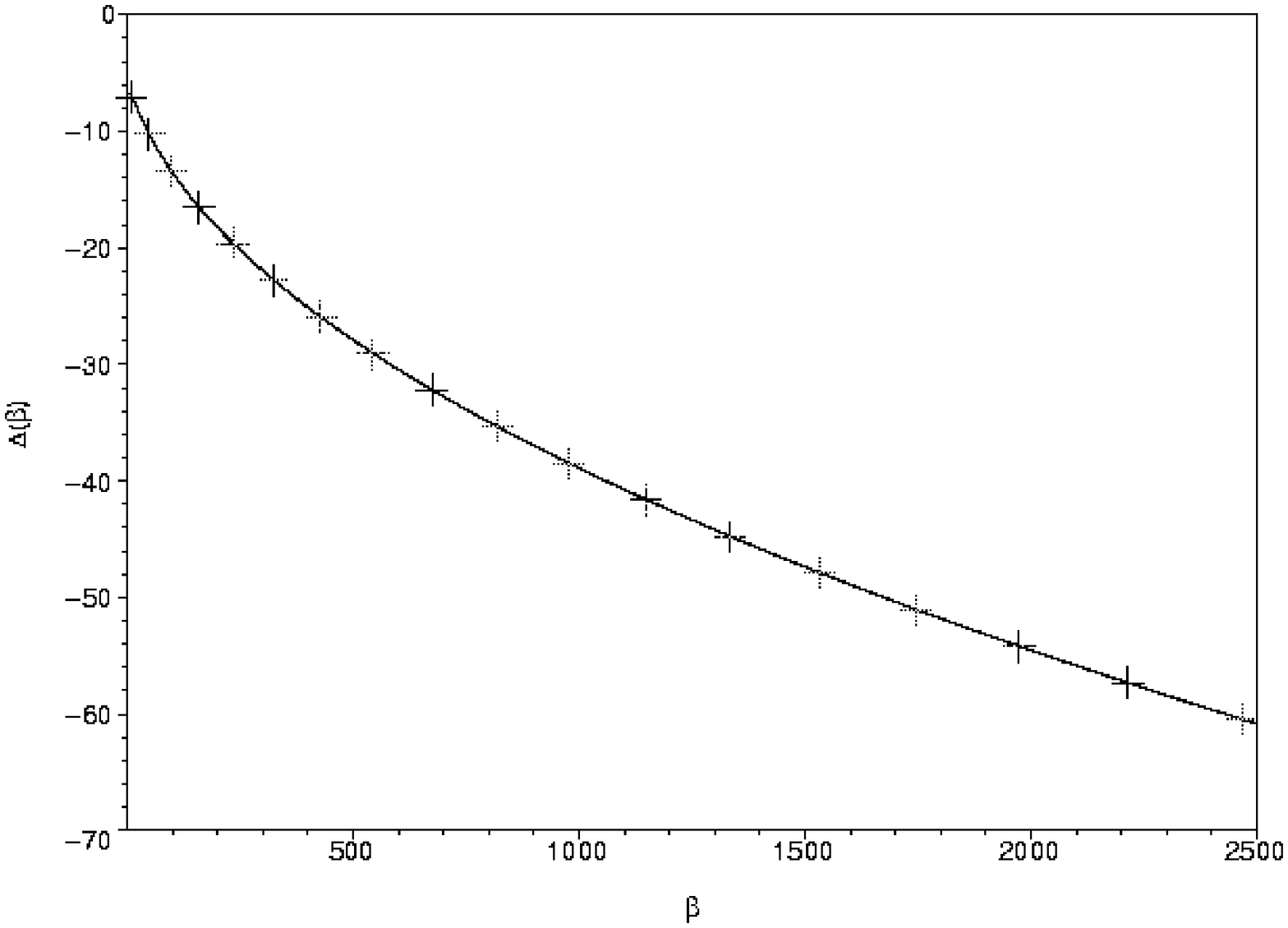}
  \caption{$\Delta^{\mbox{\scriptsize WKB}}_{2,1}$ and $\Delta^{\mbox{\scriptsize Fr}}_{2,1}$}
  \label{WKB_CF_2_1}
 \end{minipage}%
 \begin{minipage}[t]{0.5\linewidth}
  \centering
  \includegraphics[width=4cm,angle=-90]{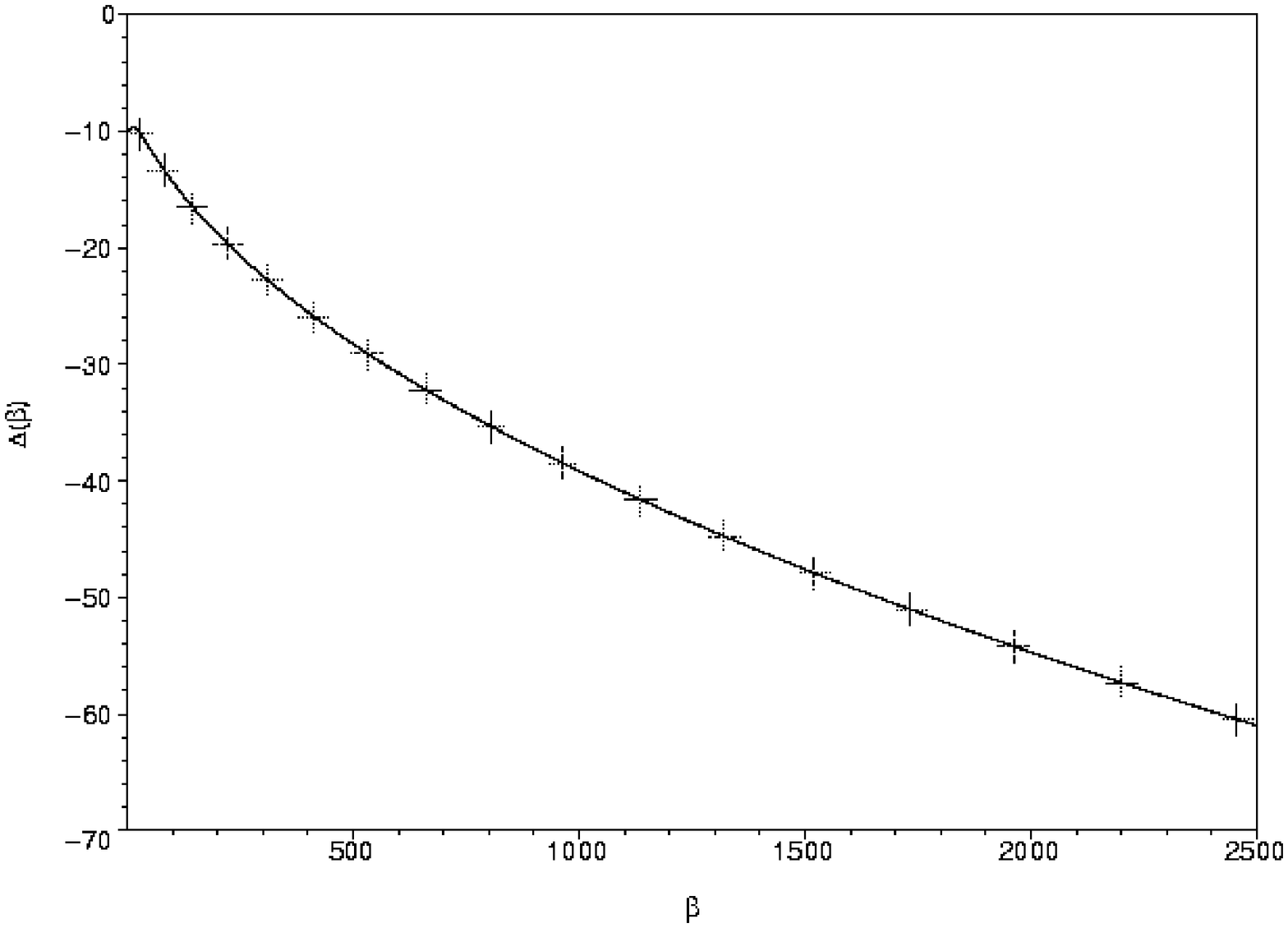}
  \caption{$\Delta^{\mbox{\scriptsize WKB}}_{3,1}$ and $\Delta^{\mbox{\scriptsize Fr}}_{3,1}$}
  \label{WKB_CF_3_1}
 \end{minipage}\\[5pt]
%\end{figure}
%\begin{figure}[ht]
 \begin{minipage}[t]{0.5\linewidth}
  \centering
  \includegraphics[width=4cm,angle=-90]{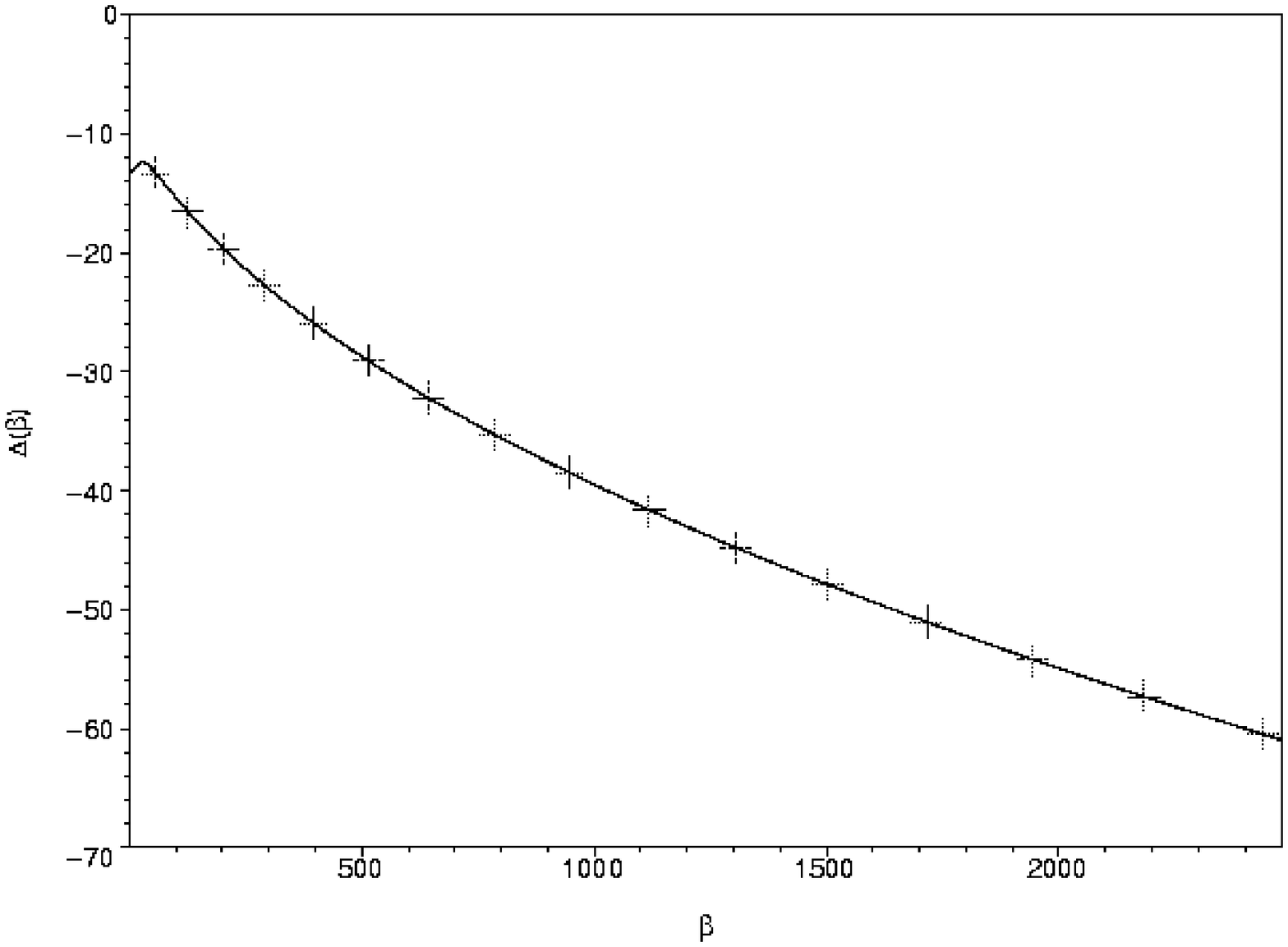}
  \caption{$\Delta^{\mbox{\scriptsize WKB}}_{4,1}$ and $\Delta^{\mbox{\scriptsize Fr}}_{4,1}$}
  \label{WKB_CF_4_1}
 \end{minipage}%
 \begin{minipage}[t]{0.5\linewidth}
  \centering
  \includegraphics[width=4cm,angle=-90]{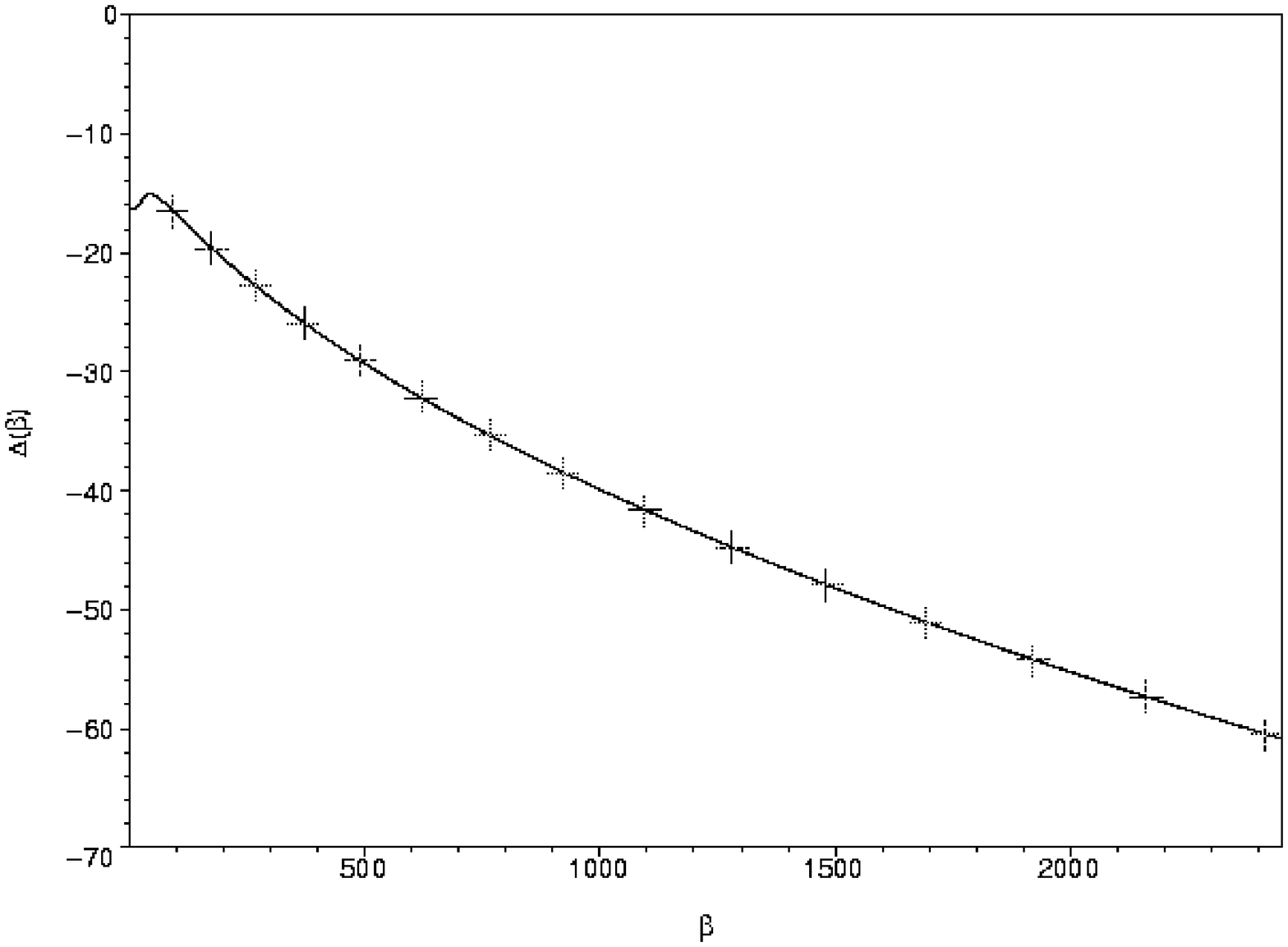}
  \caption{$\Delta^{\mbox{\scriptsize WKB}}_{5,1}$ and $\Delta^{\mbox{\scriptsize Fr}}_{5,1}$}
  \label{WKB_CF_5_1}
 \end{minipage} \\ [5pt]
%\end{figure}
%\begin{figure}[hb]
 \begin{minipage}[t]{0.5\linewidth}
  \centering
  \includegraphics[width=4cm,angle=-90]{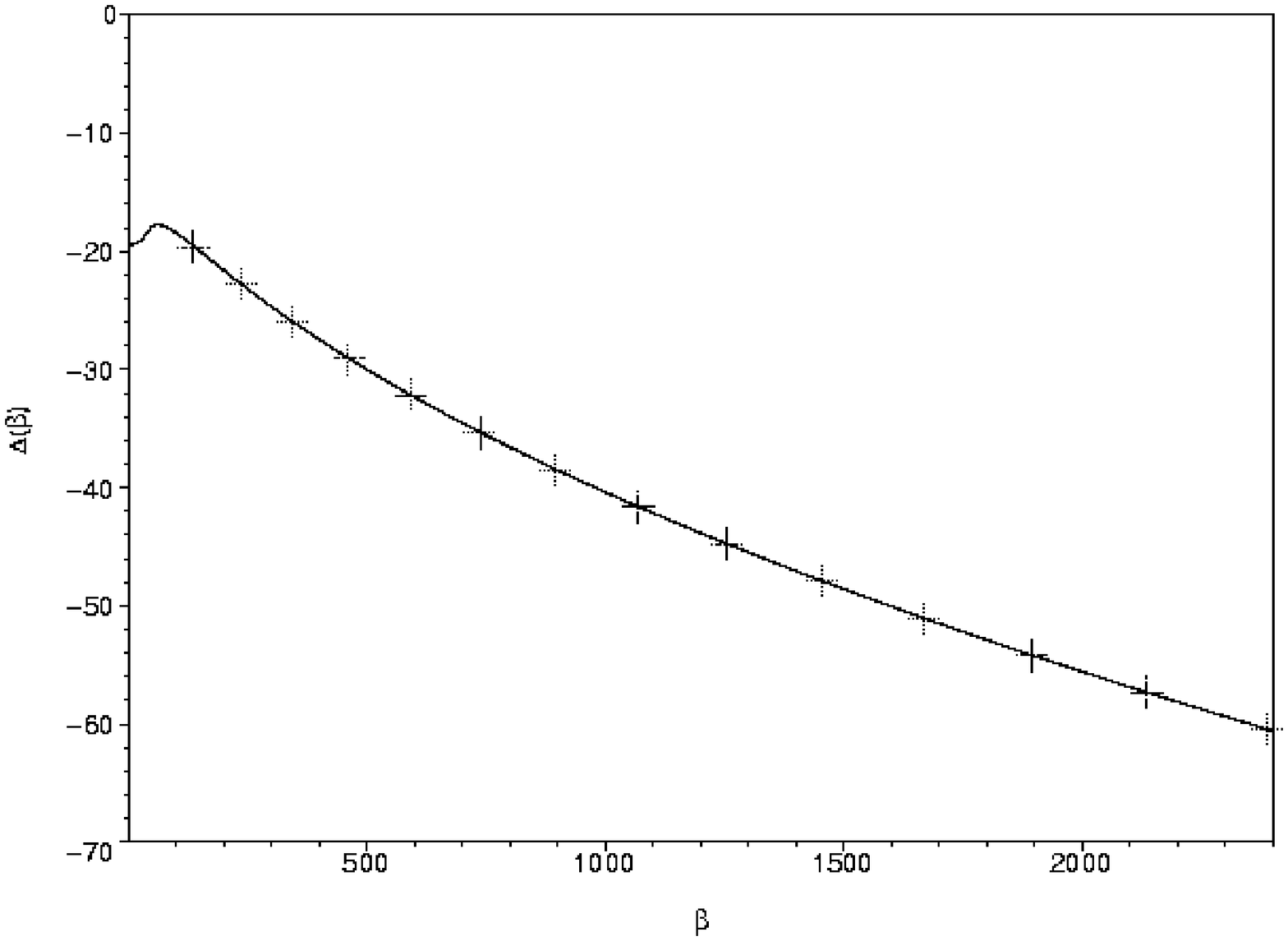}
  \caption{$\Delta^{\mbox{\scriptsize WKB}}_{6,1}$ and $\Delta^{\mbox{\scriptsize Fr}}_{6,1}$}
  \label{WKB_CF_6_1}
 \end{minipage}%
 \begin{minipage}[t]{0.5\linewidth}
  \centering
  \includegraphics[width=4cm,angle=-90]{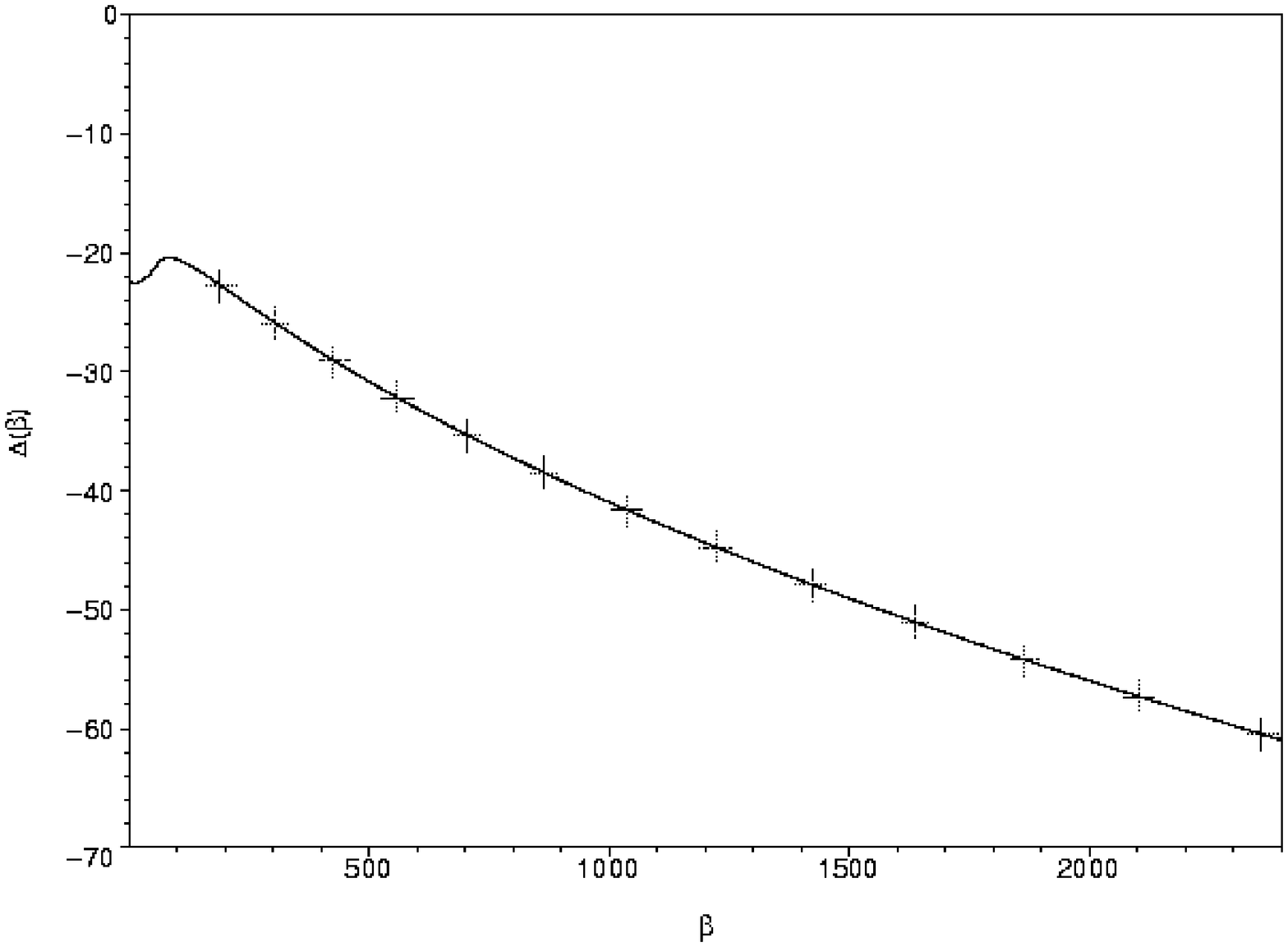}
  \caption{$\Delta^{\mbox{\scriptsize WKB}}_{7,1}$ and $\Delta^{\mbox{\scriptsize Fr}}_{7,1}$}
  \label{WKB_CF_7_1}
 \end{minipage} \\[5pt]
 \begin{minipage}[t]{0.5\linewidth}
  \centering
  \includegraphics[width=4cm,angle=-90]{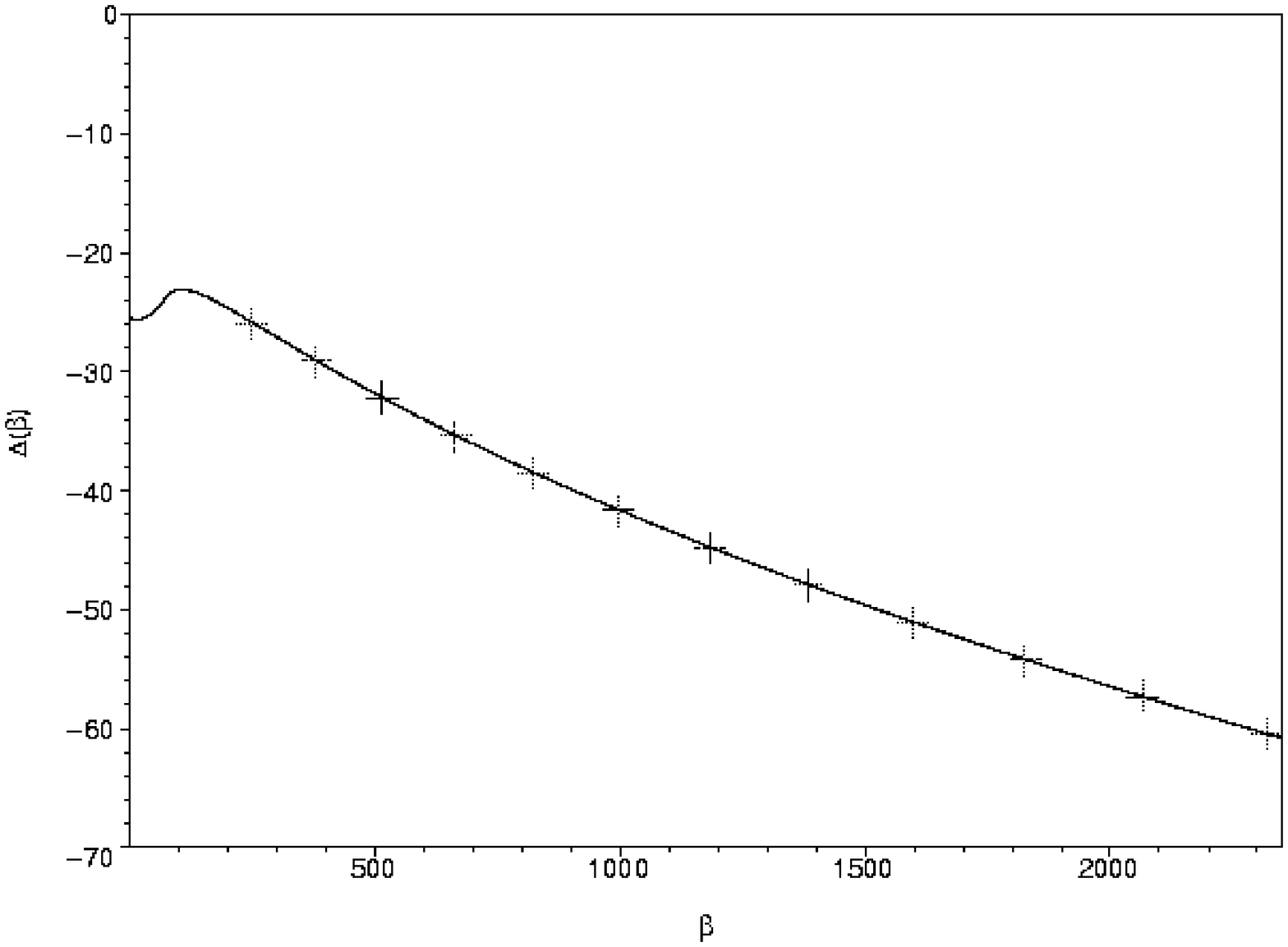}
  \caption{$\Delta^{\mbox{\scriptsize WKB}}_{8,1}$ and $\Delta^{\mbox{\scriptsize Fr}}_{8,1}$}
  \label{WKB_CF_8_1}
 \end{minipage}%
 \begin{minipage}[t]{0.5\linewidth}
  \centering
  \includegraphics[width=4cm,angle=-90]{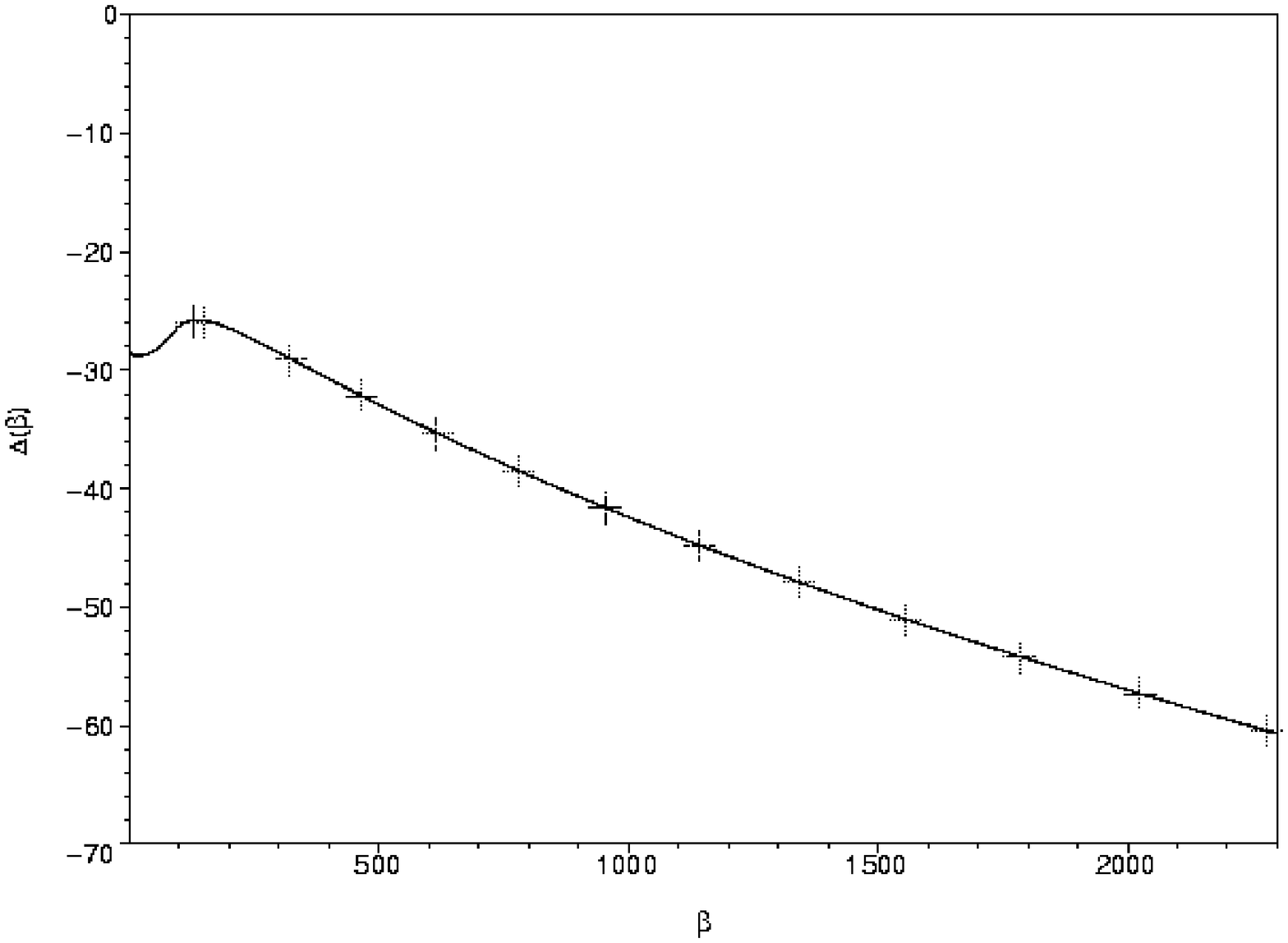}
  \caption{$\Delta^{\mbox{\scriptsize WKB}}_{9,1}$ and $\Delta^{\mbox{\scriptsize Fr}}_{9,1}$}
  \label{WKB_CF_9_1}
 \end{minipage}
\end{figure}

Figs.~\ref{WKB_CF_2_1} to \ref{WKB_CF_9_1} show $\Delta_{j,1}(\beta)$ with varying $\beta$. For $j \le 8$ it is $l(n):=n$. 
$j=9$ is the lowest value where for $q=1$ we have $\beta_1 < \beta_{\max}$, where $\beta_{\max}$ denotes the
local maximum of the WKB approximation: 
We have to relabel the $\beta_n$. This can be done by setting $\forall_{j \ge 9}: \; l(n):=|n-\frac{3}{2}|-\frac{1}{2}$.
Figs.~\ref{WKB_CF_2_1} to \ref{WKB_CF_9_1} show an exact match
of the data sets obtained by the WKB approximation and 
the Frobenius/continued-fraction method.
\medskip

\begin{figure}[ht]
 \begin{minipage}[t]{0.5\linewidth}
  \centering
  \includegraphics[width=4cm,angle=-90]{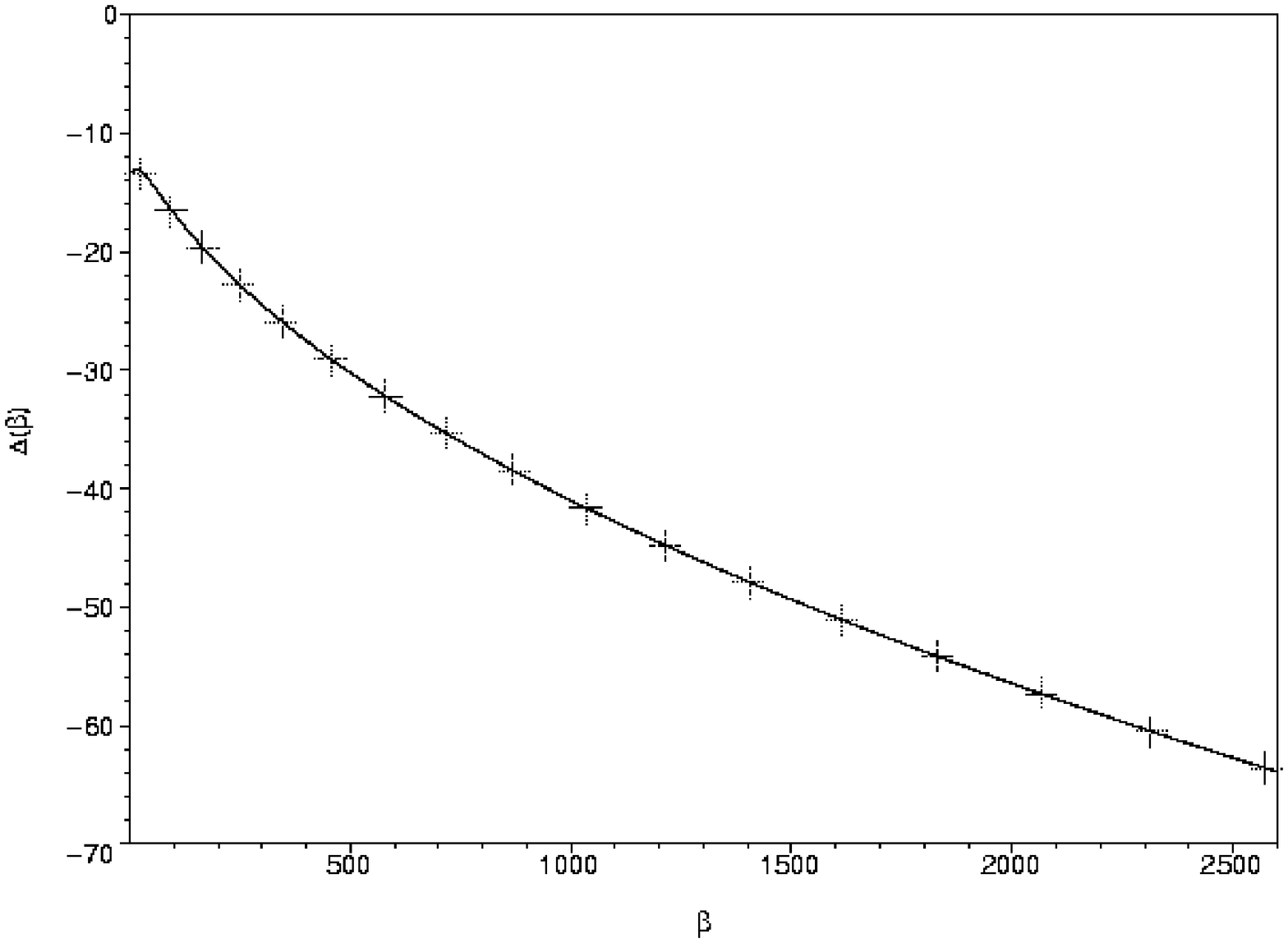}
  \caption{$\Delta^{\mbox{\scriptsize WKB}}_{4,2}$ and $\Delta^{\mbox{\scriptsize Fr}}_{4,2}$}
  \label{WKB_CF_4_2}
 \end{minipage}%
 \begin{minipage}[t]{0.5\linewidth}
  \centering
  \includegraphics[width=4cm,angle=-90]{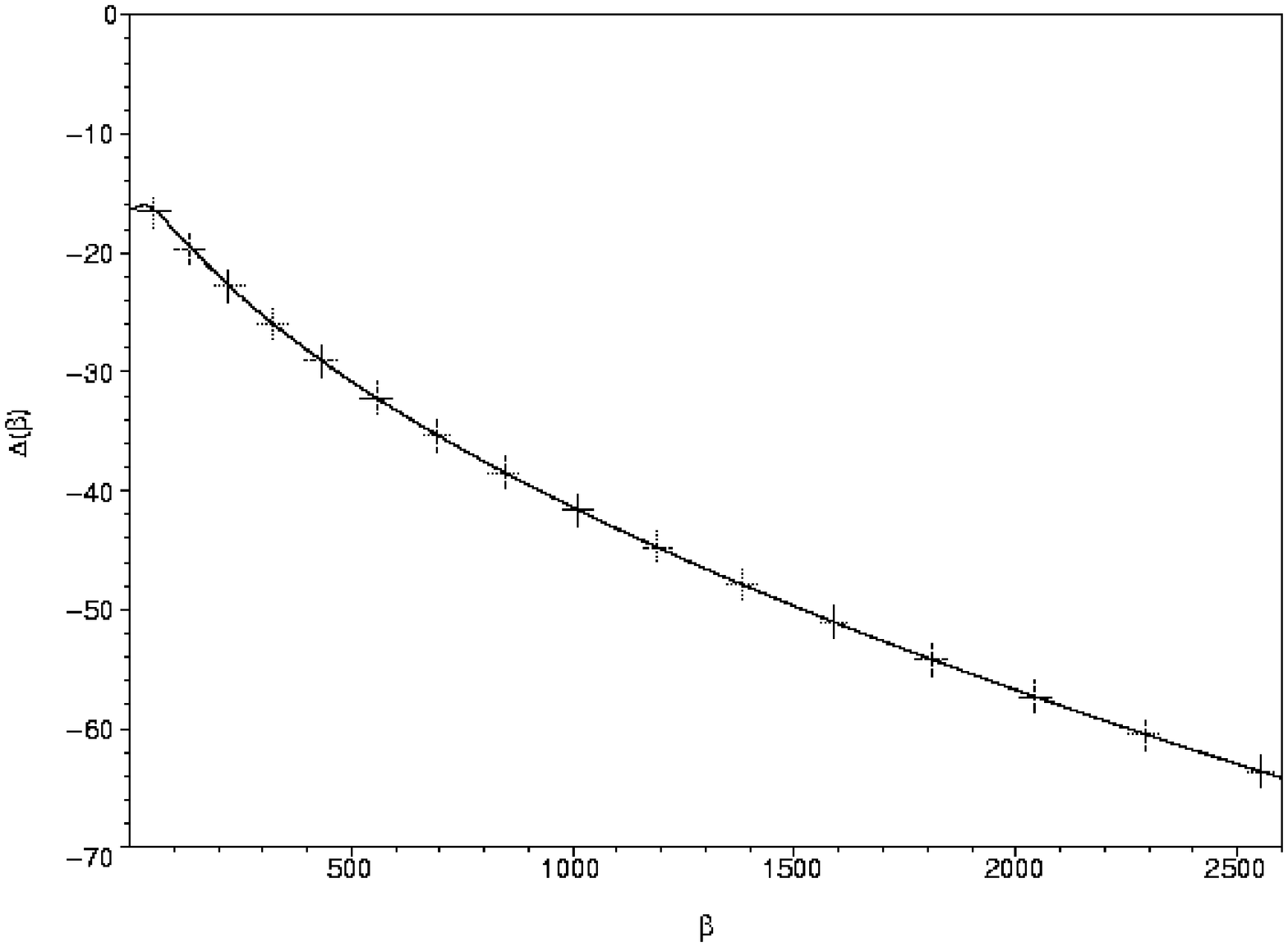}
  \caption{$\Delta^{\mbox{\scriptsize WKB}}_{5,2}$ and $\Delta^{\mbox{\scriptsize Fr}}_{5,2}$}
  \label{WKB_CF_5_2}
 \end{minipage}\\[5pt]
%\end{figure}
%\begin{figure}[ht]
 \begin{minipage}[t]{0.5\linewidth}
  \centering
  \includegraphics[width=4cm,angle=-90]{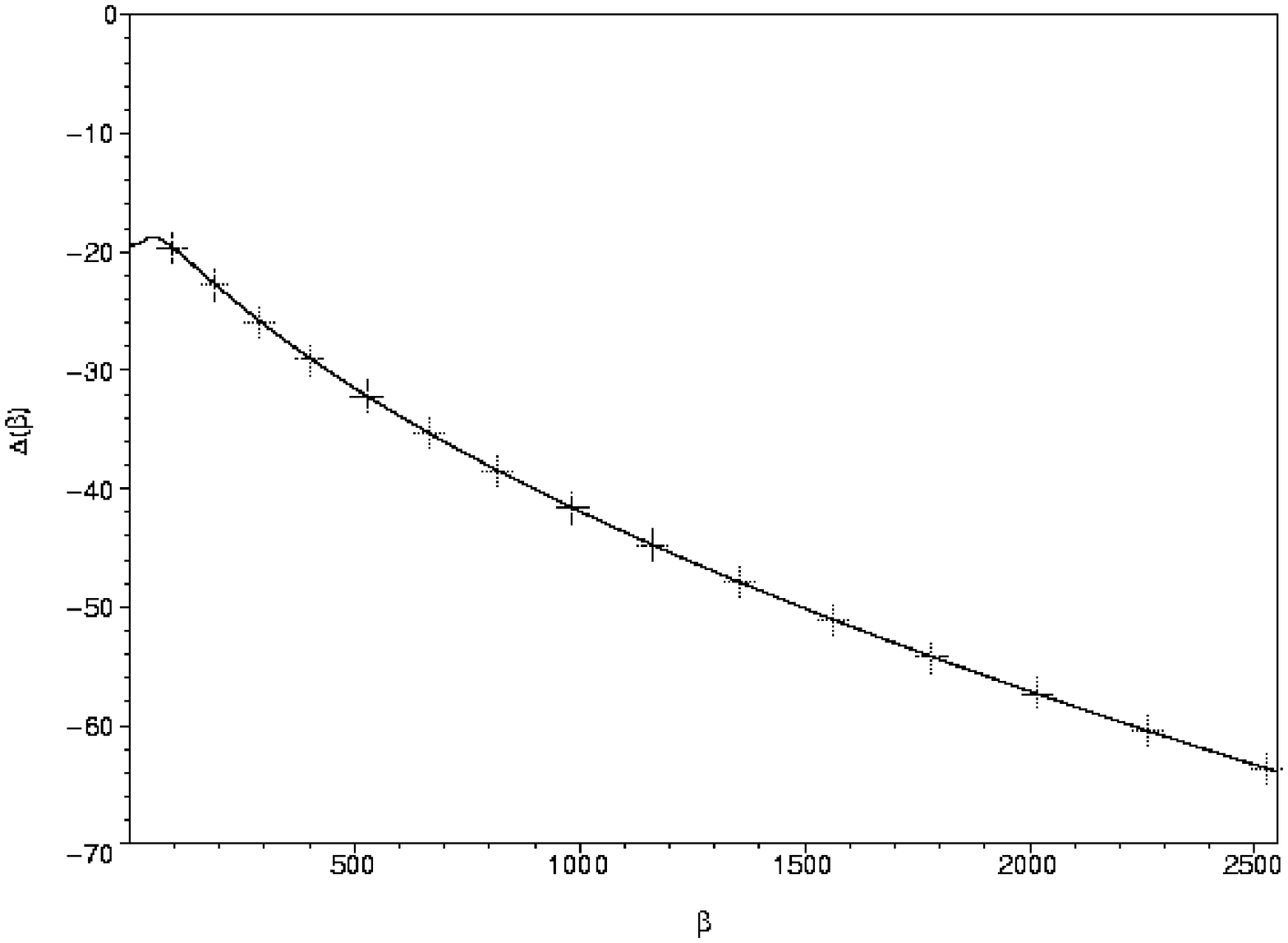}
  \caption{$\Delta^{\mbox{\scriptsize WKB}}_{6,2}$ and $\Delta^{\mbox{\scriptsize Fr}}_{6,2}$}
  \label{WKB_CF_6_2}
 \end{minipage}%
 \begin{minipage}[t]{0.5\linewidth}
  \centering
  \includegraphics[width=4cm,angle=-90]{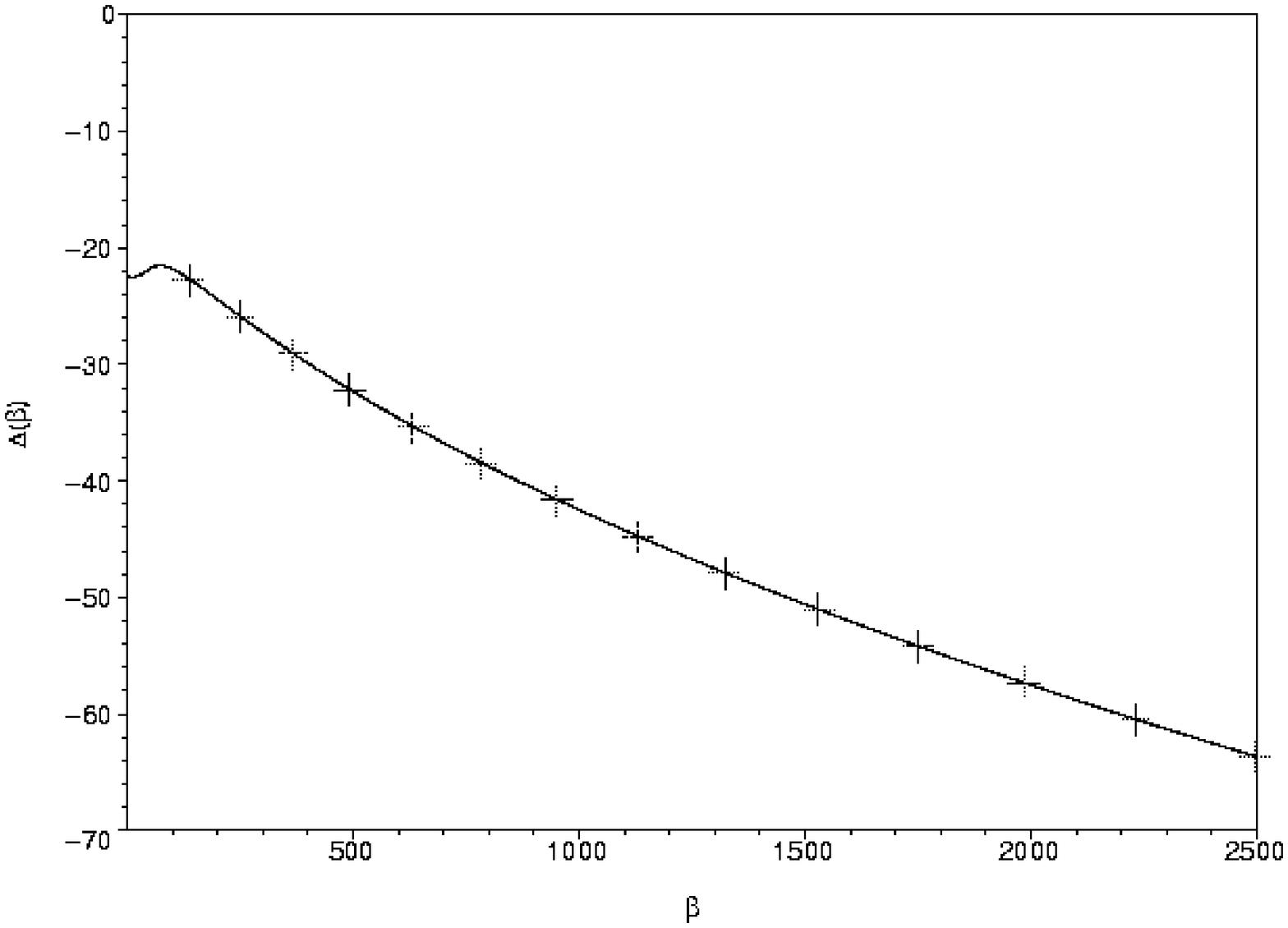}
  \caption{$\Delta^{\mbox{\scriptsize WKB}}_{7,2}$ and $\Delta^{\mbox{\scriptsize Fr}}_{7,2}$}
  \label{WKB_CF_7_2}
 \end{minipage} \\ [5pt]
%\end{figure}
%\begin{figure}[hb]
 \begin{minipage}[t]{0.5\linewidth}
  \centering
  \includegraphics[width=4cm,angle=-90]{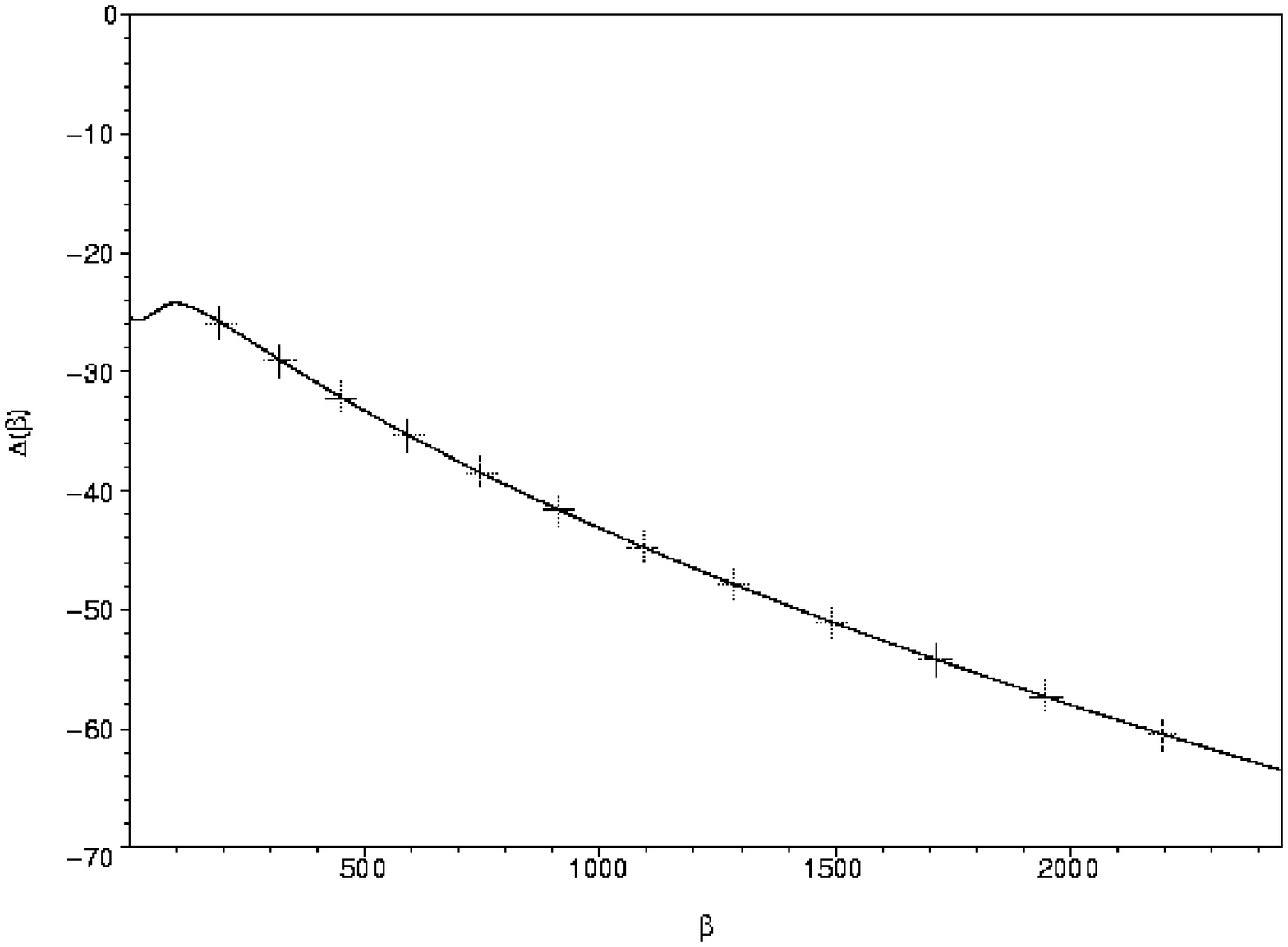}
  \caption{$\Delta^{\mbox{\scriptsize WKB}}_{8,2}$ and $\Delta^{\mbox{\scriptsize Fr}}_{8,2}$}
  \label{WKB_CF_8_2}
 \end{minipage}%
 \begin{minipage}[t]{0.5\linewidth}
  \centering
  \includegraphics[width=4cm,angle=-90]{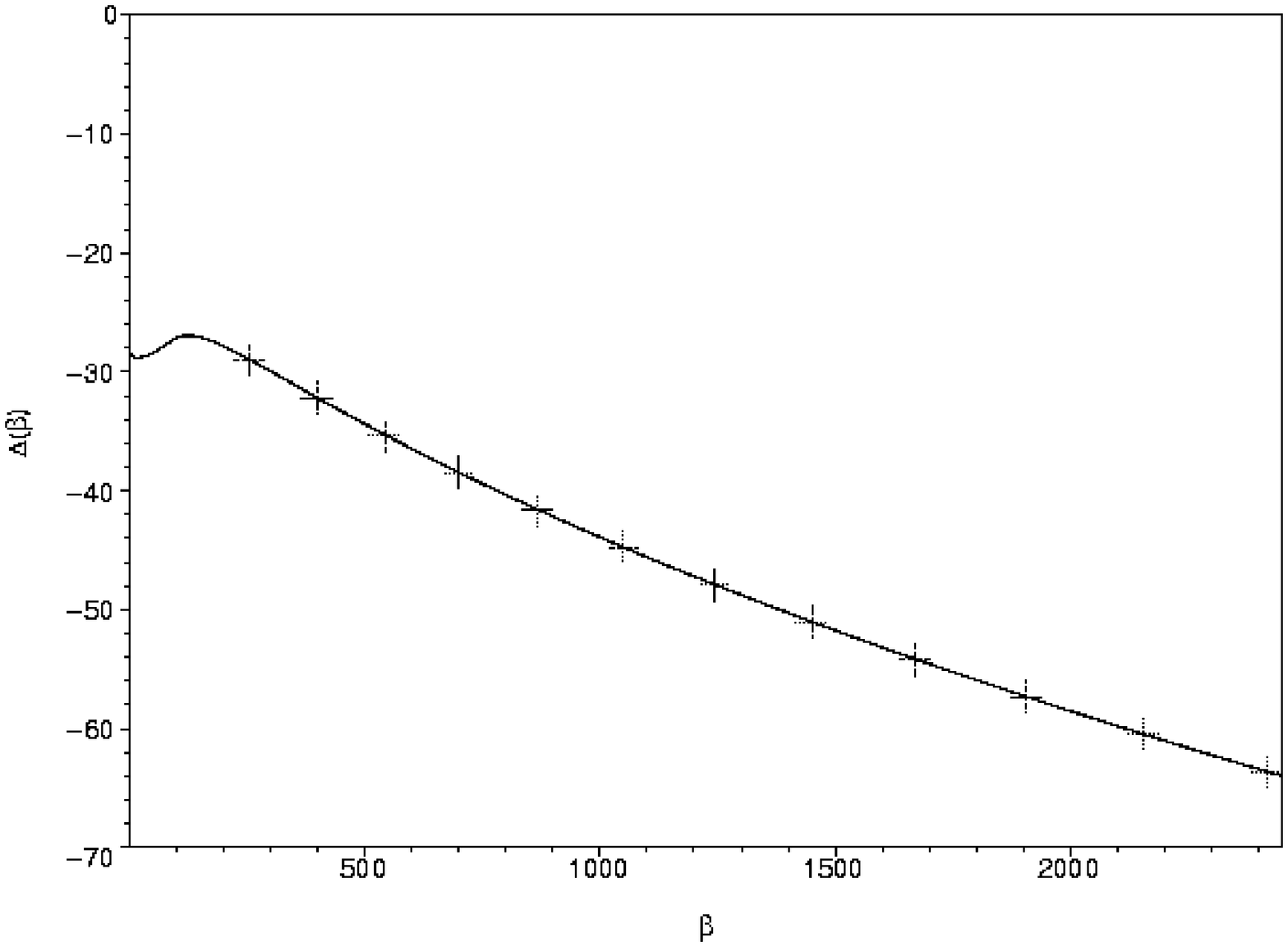}
  \caption{$\Delta^{\mbox{\scriptsize WKB}}_{9,2}$ and $\Delta^{\mbox{\scriptsize Fr}}_{9,2}$}
  \label{WKB_CF_9_2}
 \end{minipage} \\[5pt]
 \begin{minipage}[t]{0.5\linewidth}
  \centering
  \includegraphics[width=4cm,angle=-90]{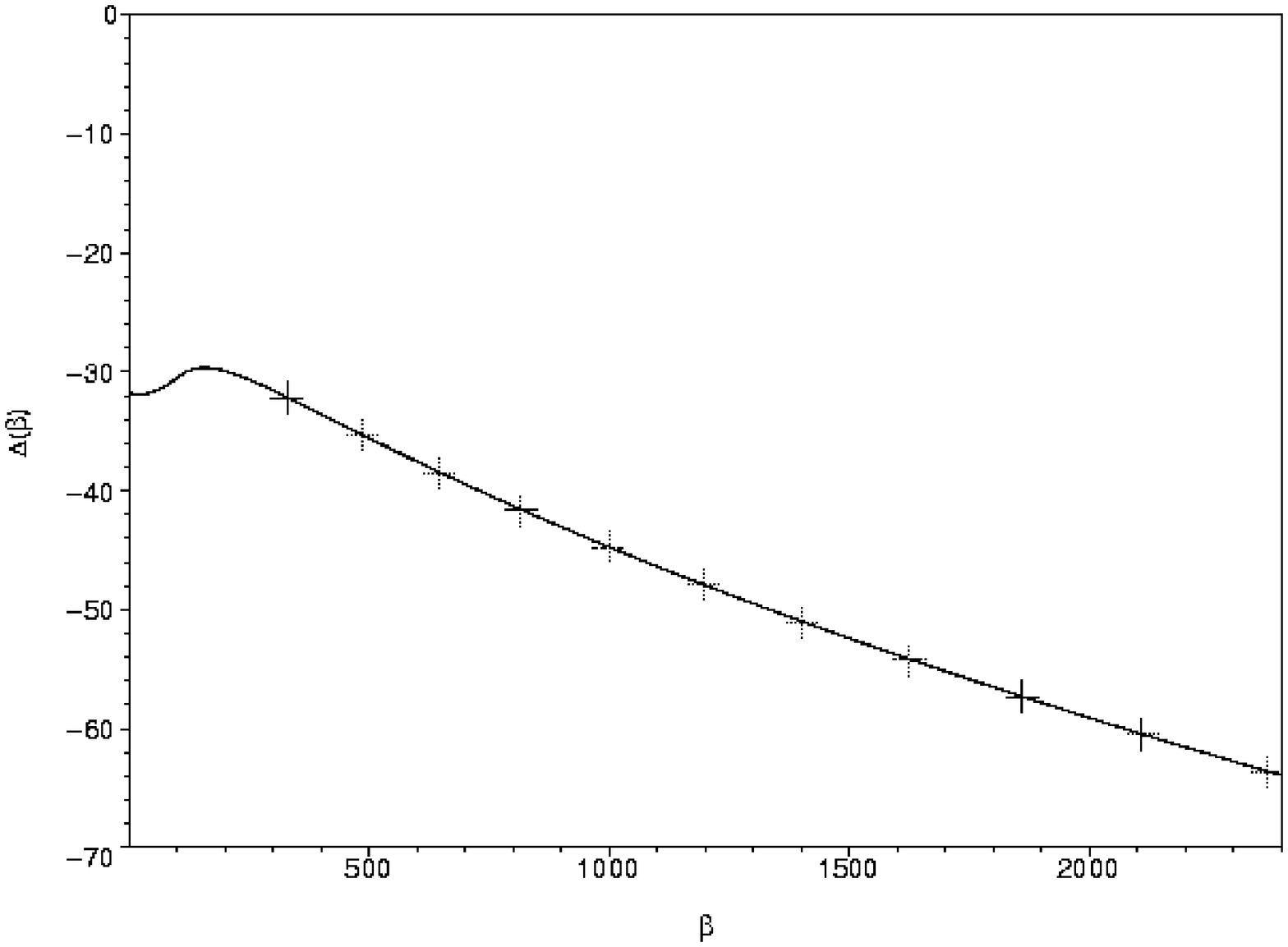}
  \caption{$\Delta^{\mbox{\scriptsize WKB}}_{10,2}$ and $\Delta^{\mbox{\scriptsize Fr}}_{10,2}$}
  \label{WKB_CF_10_2}
 \end{minipage}%
 \begin{minipage}[t]{0.5\linewidth}
  \centering
  \includegraphics[width=4cm,angle=-90]{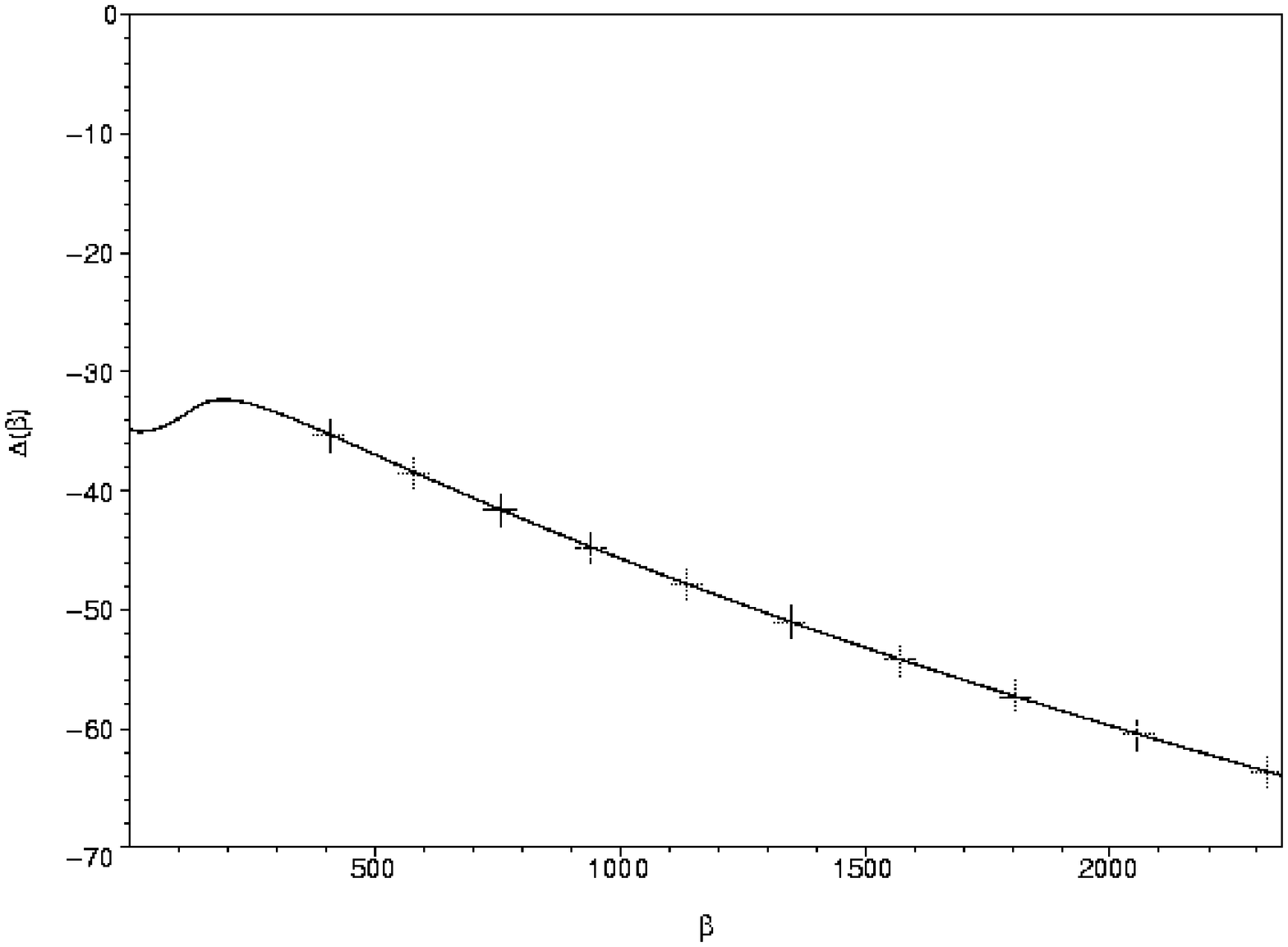}
  \caption{$\Delta^{\mbox{\scriptsize WKB}}_{11,2}$ and $\Delta^{\mbox{\scriptsize Fr}}_{11,2}$}
  \label{WKB_CF_11_2}
 \end{minipage}
\end{figure}

Figs.~\ref{WKB_CF_4_2} to \ref{WKB_CF_11_2} show $\Delta_{j,2}(\beta)$ with varying $\beta$. 
For $j \le 11$ it is $l(n):=n$. 
Figs.~\ref{WKB_CF_4_2} to \ref{WKB_CF_11_2} show an exact match
of the data sets obtained by the WKB approximation and the 
Frobenius/continued-fraction method.
\medskip

The accuracy of the WKB approximation can be explained as follows:
from the error bound in Eq.~(\ref{errorIV}) we can deduce the following numerical result
for the error $\delta_{j,q}(\beta)$ of $\Delta_{j,q}(\beta)$
\beq
 |\delta_{j,q=0}(\beta)| & \le & \frac{\pi}{2} \; \min{\left( 1,
  1.1 \, e^{ \frac{1.2}{\sqrt{\beta}}} -1 \right)} \;.
\eeq
This means that for $\beta>4$ we already have 
\beq
 |\delta_{j,q=0}(\beta)| & \le & \frac{\pi}{2} \left(e^{ \frac{1.2}{\sqrt{\beta}}} -1 \right)  \,,
\eeq
and for $\beta > 100$ we find $|\delta_{j,q=0}(\beta)| \le \frac{\pi}{10}$.

The same computations can be made for $q\ge 1$ with similar results since we also have 
a similar error bound for the general $\delta_{j,q}(\beta)$ in Eq.~(\ref{errorI}).

\section{Conclusions and outlook}
\label{conclusions}
The eigenvalue equation for the Laplace-Beltrami operator acting on
scalar functions on the non-compact Eguchi-Hanson space reduces to a
confluent Heun equation ((\ref{ODE}) or after a suitable substitution
(\ref{ODE2})) with Ince symbol $[0,2,1_2]$.

With the help of the Liouville-Green approximation (WKB) 
we have constructed approximations for the eigenfunctions
by special functions in Sect.~\ref{WKB}. Depending on the quantum numbers that label the
$SU(2)$-representation the approximating functions are either Airy,
Whittaker, or Bessel functions. Furthermore, we have derived the
scattering phases and error bounds in these cases.

Moreover, for specific discrete values of the eigenvalue in Sect.~\ref{Con_Frac} 
we have constructed the exact solutions by the
Frobenius methods. These eigenvalues are given by the zeros of a
meromorphic function defined by the infinite continuous fraction
(\ref{T-frac}).
Together with a monodromy computation this has provided us with the data 
needed for a numerical interpolation of the scattering phases.

Finally, in Sect.~\ref{results} we have shown that these two sets of data (obtained by the WKB approximation and the Frobenius/continued-fraction method) 
agree to a high accuracy. This shows that one can find a 
discrete set of exact values for the spectral density of the Laplace-Beltrami 
operator by the method described in Sect.~\ref{computation}. 
Conversely, it shows that the expressions for the 
eigenfunctions and scattering phases which were derived in Sect.~\ref{WKB} 
and which have the advantage of being given in terms of explicit functions
are very accurate approximations and can be used for all numerical purposes.

It is now interesting to ask whether the meromorphic function defined by the continued fraction
(\ref{T-frac}) can be expressed explicitly as a ratio of special functions.
If so we could obtain the discrete eigenvalues -- for which we have already 
calculated the exact scattering phase and which we have also computed numerically by the continued fraction -- as 
zeros of this ratio of special functions.
A first step towards an explicit representation of the continued fraction
might arise from the method of Pincherle. If we apply \cite[Kap.~21.84, Satz 8]{Perr:1977}
to the continued fraction (\ref{T-frac}) we obtain the following result:
a representation of the meromorphic function defined by the T-fraction (\ref{T-frac}) 
is given by
\beq
 \mathcal{M}(j,q|x) & = & q(q+1)-j(j+1)+x + \frac{\delta_{q0} \, 16\pi x^3 \, 
  \phi\left(\frac{1}{2x}\right)+\int_{\frac{1}{2x}}^{\infty}dz \; \phi(z) z^{-3}}
  {\delta_{q 0} \, 8\pi x^2 \, \phi\left(\frac{1}{2x}\right)
    +\int_{\frac{1}{2x}}^{\infty}dz \; \phi(z) z^{-2}} \;.
\eeq
Here, $\phi$ is the solution of the differential equation
\beq
 \left\lbrace z^3 \, \left(z - \frac{1}{2x} \right) \frac{d^2}{dz^2} + z^2 \left((q+1)z - \frac{q}{x} \right)
 \frac{d}{dz} + \left( \frac{x-q(q+1)+j(j+1)}{2x} z - \frac{1}{2x} \right)\right\rbrace \phi(z)
 =0 
\eeq
that behaves as $(z-\frac{1}{2x})^q$ at the regular singularity $z=\frac{1}{2x}$. An analysis
of the differential equation shows that such a solution always exists. A further analysis of the
behavior at the irregular singularity $z=\infty$  then shows that all the appearing integrals 
also exist.

A closer investigation of the meromorphic function $\mathcal{M}(j,q|x)$ is subject of our
ongoing research.

\section*{Acknowledgments}
I would like to thank Professor W.~Nahm for suggesting the problem and some
useful discussions. Moreover, I wish to thank Professor F.~W.~J.~Olver and
Professor W.~B.~Jones for bringing the articles \cite{Duns:1994} and 
\cite{Thron:1991} to my attention. 

In particular, I would like to thank Katrin Wendland for many helpful 
discussions and a lot of encouragement. I also wish to thank the Physics 
Department of the UNC at Chapel Hill for hospitality, and the \textsl{Studienstiftung 
des deutschen Volkes} and the \textsl{Deutschen Akademischen Austauschdienst} (DAAD) 
for financial support.

\begin{appendix}

\section{The different cases in WKB approximation}
\label{WKB_calculation}
In this appendix we show the explicit construction of the approximate
solutions of the differential equation (\ref{WKB-ODE}) by the 
Liouville-Green approximation (WKB) in the cases I to IV (cf. Sect.~\ref{WKB}).
The construction of an error bound in cases I to III can be found in
the references. The construction of the bound in case IV is given in
App.~\ref{error_bounds} and might explain the construction and philosophy 
behind these error bounds.

To give bounds for the different cases I to IV we need the notion of auxiliary
weight, modulus, and phase function $E,\, M$, and $\theta$: If X,Y are solutions
of the respective differential equation of the first and second kind, then
$ X= M \, E^{-1} \, \sin \theta$ and $Y = E \, M \, \cos\theta$.

In particular, in case I we will need $E,\, M,\, \theta$ for the Airy function
(see \cite[Sect. 11.2]{Olve:1974} for more details):
\beq
\mbox{Ai}(x) =  \frac{M(x)}{E(x)}  \, \sin{\theta(x)} \;, \qquad
\mbox{Bi}(x) =  E(x) \, M(x) \, \cos{\theta(x)} \;.
\eeq         
In the cases III and IV, we will use the functions $E_0, \, M_0, \, \theta_0$ for
the Bessel function (see \cite[Sect. 12.1, 12.3]{Olve:1974} for more details):
\beq
\mbox{J}_0(x) =  \frac{M_0(x)}{E_0(x)}  \, \cos{\theta_0(x)} \;, \qquad
\mbox{Y}_0(x) =  E_0(x) \, M_0(x) \, \sin{\theta_0(x)} \,.
\eeq
For case II, the definition of modulus and weight function is the most complicated. Therefore,
we refer to \cite[Chap. 2, Th. 1]{Duns:1994} for the quite extensive definitions in this case.

\subsection{Case I}
\label{case_I}
In case I it is easy to prove that the function $f$ has a transition point 
(i.e.~simple zero) at $z > 1$ which we denote by $z_0$. 
The idea is now to perform the transformation of the variable $z$
and the function $w(z)$ to $\zeta, \, W(\zeta)$ according to
\beq
 \left\lbrace \begin{array}{lcrcl}
 \forall z \ge z_0: & \; & \frac{2}{3} (-\zeta)^\frac{3}{2} & =  &\int_{z_0}^z \sqrt{f(t)} \, dt \,,\\
 \vspace*{-0.2cm} \\
 \forall z \le z_0: & \; & \frac{2}{3} \zeta^\frac{3}{2}    & =  &\int_z^{z_0} \sqrt{-f(t)} \, dt \,,
 \end{array} \right. 
\eeq
and $w(z) = \sqrt{- \frac{dz}{d\zeta}} \, W(\zeta)$. Note that $\zeta \to - \infty$
corresponds to $z \to \infty$, $\zeta \to 0-$ to $z \to z_0+$, and $\zeta \to \infty$
to $z \to 1$. Eq.~(\ref{WKB-ODE}) becomes
\beqn
 \label{ODE-I}
  \frac{d^2}{d\zeta^2} W(\zeta)   & = & \left[ u^2 \zeta + \psi(\zeta)
\right] W(\zeta) \;,\\
 \nonumber
  \psi(\zeta) & = & \frac{5}{16\zeta^2} - \left[ 4f(z) f''(z) - 5f'(z)^2
\right] \frac{\zeta}{16 f(z)^3}
- \frac{\zeta g(z)}{f(z)} \;.
\eeqn
Approximate solutions of (\ref{WKB-ODE}), i.e. solutions of (\ref{ODE-I})
with $\psi = 0$, that are regular at $z=1$, are given by \cite[Sect. 11.3, Th. 3.1]{Olve:1974}
\beqn
\label{sol-I}
w(z) & = &  \sqrt[4]{- \frac{\zeta}{f(z)}} \left[ \mbox{Ai}(u^\frac{2}{3} \zeta)
+ \varepsilon(u,\zeta) \right] \,,\\
\nonumber
\mbox{with} \quad |\varepsilon(u,\zeta)| & \le & \frac{1}{\lambda} \, \frac{M(u^\frac{2}{3}\zeta)}
{E(u^\frac{2}{3}\zeta)} \, \left[ e^{\frac{2\lambda}{u} \mathcal{V}_{\zeta,\infty}\left(|\zeta|^\frac{1}{2} \mathcal{B}_0(\zeta)\right) } -1 \right] \;. 
\eeqn
The constant $\lambda$ and the function $\mathcal{B}_0$ are defined as follows:
\beq
 \lambda & := & \sup_x{\left\lbrace\pi \, |x|\frac{1}{2} M^2(x)\right\rbrace}\,,\\
 \mathcal{B}_0(\zeta) & := & \frac{1}{2\sqrt{|\zeta|}}
\int_\zeta^\infty \frac{dv}{|v|^\frac{1}{2}} \, \psi(u,v)\,,
\eeq
and $\mathcal{V}$ is the variational operator, i.e.
\beq
 \mathcal{V}_{\zeta,\infty}\left(|\zeta|^\frac{1}{2} \mathcal{B}_0(\zeta)\right) 
 & = & \frac{1}{2} \int_\zeta^\infty \frac{dv}{|v|^\frac{1}{2}} |\psi(u,v)| \;.
\eeq
Following the discussion in \cite[Sect. 13.7.2]{Olve:1974} we can determine
the behavior for $z \to 1$ and $z \to \infty$:
\beq
A(z) & \ub{\longrightarrow}{z\to 1} & \frac{\sqrt[4]{\beta}}{2\sqrt{\pi q}}
(z-1)^\frac{q}{2} \;,\\
A(z) & \ub{\sim}{z\to \infty} & \frac{1}{\sqrt{\pi} \, z^\frac{3}{4}} \,
\sin{\left(u \int_{z_0}^z \sqrt{f(t)} \, dt + \frac{\pi}{4}
+ \delta \right)} \,,
\eeq
where we have used $\lim_{\zeta \to \infty} \varepsilon(u,\zeta) = 0$, and
the phase $\delta$ is determined by 
$\lim_{\zeta \to - \infty} \varepsilon(u,\zeta)$. We know that the argument of 
the $\sin$-function equals $2\sqrt{\beta z}+ \Delta_{j,q}$ for large $z$, i.e.
\beq
A(z) & \ub{\sim}{z\to \infty} &
\frac{1}{\sqrt{\pi} \, z^\frac{3}{4}} \, \sin{\left( 2\sqrt{\beta z}
+ \Delta_{j,q} \right)} \;.
\eeq
Therefore, for the scattering phase we obtain
\beq
\Delta_{j,q} & = & u \lim_{z \to \infty} \left( \int_{z_0}^z \sqrt{f(t)} \, 
dt - 2 \sqrt{z} \right)
+ \frac{\pi}{4} + \delta \,,
\eeq
and with the help of \cite[Sect. 11.2]{Olve:1974} it follows that
\beqn
 \label{errorI}
\frac{2|\delta|}{\pi} & \le &
\min{\left\{ 1, \frac{1}{\lambda}
  \left[ e^{ \frac{2\lambda}{u} \,
\mathcal{V}_{-\infty,\infty}\left(|\zeta|^\frac{1}{2} \mathcal{B}_0(\zeta)\right)
} -1 \right]
\right\}} \;.
\eeqn
To resume, the solution (\ref{sol-I}) can be considered as relation between
the recessive solution
and its asymptotic expansion. In this sense, the solution is a {\em
connection formula} and it is known as
{\em Gans-Jeffreys} formula.

\subsection{Case II}
\label{case_II}
In case II, $f(z)=\frac{z-a}{z^2-1}$, we can apply a result of Olver and Nestor that has been
generalized in \cite{Duns:1994}. The Liouville transformation takes the form
\beq
 \int_\alpha^\zeta \left( \frac{\tau - \alpha}{\tau} \right)^\frac{1}{2} \, d\tau
 & = & \int_a^z dt \sqrt{f(t)} \\
\mbox{with} \qquad \alpha & := & \frac{2}{\pi} \int_1^a \sqrt{-f(t)} \, dt \;.
\eeq
Therefore, we perform the transformation of the variable $z$
and the function $w(z)$ to $\zeta,\, W(\zeta)$ according to
\beq
 \zeta^\frac{1}{2} (\zeta-\alpha)^\frac{1}{2} - \frac{\alpha}{2}
\ln{\left( \frac{2\zeta - \alpha +2
\zeta^\frac{1}{2}(\zeta-\alpha)^\frac{1}{2}}{\alpha} \right)}
 =  \int_a^z \sqrt{f(t)} \, dt \\
\eeq
and $w(z) = \sqrt{ \frac{dz}{d\zeta}} \, W(\zeta)$. The integral can be expressed
through an elliptic integral of the second kind:
\beq
\int_a^z \sqrt{f(t)} \, dt 
& = & 2 \sqrt{\frac{(z-a)(z+1)}{z-1}}
- 2 \sqrt{1+a} \;\; \mbox{E}\left[ \arcsin{\left( \sqrt{ \frac{z-a}{z-1}}
\right)},\sqrt{\frac{2}{1+a}} \right] \,.
\eeq
Note that $\zeta \to \alpha$ corresponds
to $z \to a$, $\zeta \to 0$ to $z \to 1$, where branches of $\zeta$ must be chosen such that
it is an analytic function at both values. Moreover, $\zeta \to \infty$ corresponds
to $z \to \infty$. Eq.~(\ref{WKB-ODE}) becomes
\beq
  \frac{d^2}{d\zeta^2} W(\zeta)   & = & \left[ u^2
\frac{\alpha-\zeta}{\zeta} - \frac{1}{4\zeta^2} + \frac{\psi(\zeta)}{\zeta}
\right] W(\zeta) \;,\\
  \psi(\zeta) & = & \frac{4\zeta^2+\alpha^2}{16\zeta \,(\alpha-\zeta)^2} +
\left[ 4f(z) f''(z) - 5f'(z)^2 \right] \frac{\zeta-\alpha}{16 f(z)^3} \\ & +
& \frac{(\zeta -\alpha) \, g(z)}{f(z)} \;.
\eeq
Approximate solutions of (\ref{WKB-ODE}), i.e. solutions of (\ref{ODE-I})
with $\psi = 0$, that are regular at $z=1$, are given in 
\cite[Chap. 2, case 1, $m=0$]{Duns:1994} in terms of the Whittaker function $\mbox{M}$, i.e.
\beqn
\label{sol-II}
w(z) & = & \sqrt[4]{\frac{\zeta-\alpha}{\zeta \,f(z)}} \left[
e^{-\frac{i\pi}{4}} \, \sqrt{\frac{2\pi}{1+e^{\pi u \alpha}}} \, 
\mbox{M}_{\frac{i u\alpha}{2},0} \left(2i u\zeta
\right)  + \varepsilon(u,\zeta) \right] \\
\nonumber
\mbox{with} \quad |\varepsilon(u,\zeta)| & \le & 
\frac{M^{(1)}_{\frac{u\alpha}{2},0}(2u\zeta)}{E^{(0)}_{\frac{u\alpha}{2},0}(2u\zeta)} \,
\left[ e^{\frac{\kappa_0 \left(\frac{u\alpha}{2}+1\right)^\frac{1}{3}}{2u} \mathcal{V}_
{\mathcal{P}^{(0)}}(F)} - 1 \right] \;.
\eeqn
For the error bound we have adopted the notation of \cite[Chap. 2, Th. 1]{Duns:1994}. Using 
formula (2.19) in \cite{Duns:1994} we determine the behavior of the solution for $z \to 1$
\beq
A(z) & \ub{\longrightarrow}{z\to 1} & 
 \sqrt{\frac{2\pi}{1+e^{\pi u \alpha}}} \; \sqrt[4]{\frac{2\alpha u^2}{a-1}}
 \; \left( \frac{\zeta}{z-1} \right)^\frac{1}{4} \,,
\eeq
where we have used $\lim_{\zeta \to 0} \varepsilon(u,\zeta)=0$.
A close examination of the Liouville transformation then shows that
$\zeta/(z-1) \to (a-1)/(2\alpha)$ for $\zeta \to 0$ and $z \to 1$.
Using \cite[(2.16)-(2.18)]{Duns:1994} we obtain for $z \to \infty$
\beq
A(z) & \ub{\sim}{z \to \infty} & \frac{2}{z^\frac{3}{4}}
\sin{\left( u\zeta - \frac{u\alpha}{2} \ln{(2u\zeta)} + \frac{\pi}{4} + \gamma
+ \delta \right)}\,,
\eeq
where 
\beq
\gamma & = & \mbox{arg} \; \Gamma\left( \frac{1}{2} + \frac{iu\alpha}{2}
\right)
\eeq
and the phase $\delta$ is determined by $\lim_{\zeta \to \infty} \varepsilon(u,\zeta)$.
We know that the argument in the $\sin$-function equals $ 2\sqrt{\beta z} + 
\Delta_{j,0}$
for large $z$, i.e.
\beq
A(z) & \ub{\sim}{z \to \infty} & \frac{1}{z^\frac{3}{4}}
\sin{\left( 2\sqrt{\beta z} + \Delta_{j,0} \mid_{a>1} \right)} \;.
\eeq
A careful examination of the Liouville transformation reveals the following relation
between large $\zeta$ and $z$
\beq
 \zeta - \frac{\alpha}{2} - \frac{\alpha}{2} \ln{(\zeta)} + \frac{\alpha}{2}\ln{(\frac{\alpha}{4})}
 + O(\zeta^{-1})
 = 2 \sqrt{z} - 2 \sqrt{1+a}\; \;
\mbox{E}\left[\sqrt{\frac{2}{1+a}} \right] + O(z^{-\frac{1}{2}})\;.
\eeq 
Therefore, for the scattering phase we obtain
\beqn
\label{errorII}
\nonumber
\Delta_{j,0}\mid_{a>1} & = & u \lim_{z\to \infty} \left( \zeta - 
\frac{\alpha}{2} \ln{(2u\zeta)}
- 2\sqrt{z} \right) + \frac{\pi}{4} + \gamma + \delta\\
\nonumber
& = & - 2u\, \sqrt{1+a} \; \;
\mbox{E}\left[\sqrt{\frac{2}{1+a}} \right]
+ \frac{\alpha u}{2} - \frac{\alpha u}{2} \ln{(\frac{\alpha u}{2})}+  \frac{\pi}{4} + \gamma \\
\nonumber
  \gamma & = & \mbox{arg} \; \Gamma\left( \frac{1}{2} + \frac{iu\alpha}{2}
\right) \\
 \frac{2|\delta|}{\pi} & \le & \min{\left\{ 1,\left[ e^{\frac{\kappa_0 
 \left(\frac{u\alpha}{2}+1\right)^\frac{1}{3}}{2u} \mathcal{V}_{\mathcal{P}^{(0)}}(F)} - 1 
\right]\right\}}\;.
\eeqn
The phase $\Delta_{j,0}\mid_{a>1}$ differs from the one obtained in
\cite{Mign:1991}.

\subsection{Case III}
\label{case_III}
In case III, the function $f$ has  neither a transition point nor a pole
for $z\ge 1$. This case can be understood as the limit $a \to 1$ of case II, 
and applying further results of Olver's \cite{Olve:1977}, in this limit the 
Whittaker function will become the Bessel function $\mbox{J}_0$.

The idea is to perform the transformation of the variable $z$ and the
function $w(z)$ to $\zeta,\, W(\zeta)$ according to
\beq
\zeta = \int_{1}^z \sqrt{f(t)} \, dt = 2 \sqrt{z+1} - 2\sqrt{2}
\eeq
and $w(z) = \sqrt{ \frac{dz}{d\zeta}} \, W(\zeta)$. Note that $\zeta \to \infty$
corresponds to $z \to \infty$ and $\zeta \to 0$ to $z \to 1$. Eq.~(\ref{WKB-ODE}) 
becomes
\beqn
\label{ODE-III}
  \frac{d^2}{d\zeta^2} W(\zeta)   & = & \left[ - u^2 - \frac{1}{4\zeta^2} +
\frac{\psi(\zeta)}{\zeta} \right] W(\zeta) \;,\\ 
\nonumber
  \psi(\zeta) & = & \frac{1}{4\zeta} + \left[ 4f(z) f''(z) - 5f'(z)^2
\right] \frac{\zeta}{16 f(z)^3}
+ \frac{\zeta \, g(z)}{f(z)} \;.
\eeqn
Approximate solutions of (\ref{WKB-ODE}), i.e. solutions of (\ref{ODE-III})
with $\psi=0$, that are regular at $z=1$, are given by
\cite[Chap. 5, Th. 2, case 2, $\nu=0,\, \mu=0$]{Olve:1977}, i.e.
\beqn
\label{sol-III}
w(z) & = & \frac{1}{\sqrt[4]{f(z)}} \left[ \zeta^\frac{1}{2} \, \mbox{J}_0(u
\zeta)  + \varepsilon(u,\zeta) \right] \,,\\
\nonumber
\mbox{with} \quad |\varepsilon(u,\zeta)| & \le & \frac{l_{0,1}}{l_{0,0}} \, \frac{M_0(u\zeta)}{E_0(u\zeta)} \, \left[ e^{l_{0,0} \mathcal{V}_{0,\zeta}(H)} -1 \right] \;. 
\eeqn
The constants $l_{0,0}$, $l_{0,1}$, and the function $H$ are defined as follows:
\beq
 l_{0,0} & := & \sup_x{\left\lbrace\pi \, \Omega_0(x) M_0^2(x)\right\rbrace} \,,\\
 l_{0,1} & := & \sup_x{\left\lbrace\pi \, \Omega_0(x) \, |J_0(x)| \, E_0(x) M_0(x)\right\rbrace} \,,\\
 H(u,\zeta) & := & \frac{1}{2} \int \frac{\psi(u,\zeta)}{\Omega_0(\zeta)} \, d\zeta \,,\\
 \Omega_0(x) & := & \frac{1+x}{\ln{(e + \frac{1}{x})}}\,,
\eeq
and $\mathcal{V}$ is the variational operator, i.e.
\beq
 \mathcal{V}_{0,\zeta}(H) & = & \frac{1}{2} \int_0^\zeta \frac{dv}{\Omega_0(uv)} \, |\psi(u,v)|\;.
\eeq
We determine the behavior of the solutions for $z \to 1$ and $z \to \infty$:
\beq
A(z) & \ub{\longrightarrow}{z\to 1}     & \frac{1}{\sqrt{2}} \,
\mbox{J}_0\left(u \, \frac{z-1}{\sqrt{2}}\right) \to \frac{1}{\sqrt{2}}
\;,\\
A(z) & \ub{\sim}{z\to\infty} & \frac{\sqrt{2}}{\sqrt{\pi u} z^\frac{3}{4}}
\sin{\left(u \int_1^z \sqrt{f(t)} dt + \frac{\pi}{4} +  \delta \right)} \,,
\eeq
where we have used $\lim_{\zeta \to 0+} \varepsilon(u,\zeta) = 0$, and
the phase $\delta$ is determined by $\lim_{\zeta \to \infty} \varepsilon(u,\zeta)$.
We know that the argument of the $\sin$-function equals $2\sqrt{\beta z} + 
\Delta_{j,0}\mid_{a=1}$ for large $z$, i.e.
\beq
A(z) & \ub{\sim}{z\to\infty} &  \frac{\sqrt{2}}{\sqrt{\pi u} z^\frac{3}{4}} 
\sin{\left(
2\sqrt{\beta z} + \Delta_{j,0}\mid_{a=1} \right)} \;.
\eeq
Therefore, for the scattering phase we obtain
\beq
\Delta_{j,0}\mid_{a=1} & = & u \lim_{z\to \infty}  \left( \int_1^z 
\sqrt{f(t)} \, dt - 2\sqrt{z} \right) + \frac{\pi}{4} + \, \delta\\
& = & - \sqrt{8 \beta} + \frac{\pi}{4} + \, \delta
\eeq
and with the help of \cite[Chapt.~5]{Olve:1977} and analogous to 
App.~\ref{error_bounds} the bound for $\varepsilon$ in (\ref{sol-III})
gives a bound for the phase $\delta$. It follows that
\beqn
\label{errorIII}
\frac{2|\delta|}{\pi} & \le &
\min{\left\{ 1, \frac{l_{0,1}}{l_{0,0}}
  \left[ e^{ \frac{l_{0,0}}{u} \, \mathcal{V}_{0,\infty}(H) } -1 \right]
\right\}} \;.
\eeqn
The phase $\Delta_{j,0}$ differs from the one obtained in \cite{Mign:1991}
by $\frac{\pi}{4}$.

\subsection{Case IV}
\label{case_IV}
In case IV the function $f$ has no transition point for $z\ge 1$ and 
takes the form $f(z)= \frac{z-a}{z^2-1}$. The idea is now to perform 
the transformation of
the variable $z$ and the function $w(z)$ to $\zeta,\, W(\zeta)$
according to
\beq
  (-\zeta)^\frac{1}{2} & = & \int_1^z \sqrt{f(t)} \, dt \,\\
\eeq
and $w(z) = \sqrt{- \frac{dz}{d\zeta}} \, W(\zeta)$. 
The integral can be expressed through elliptic integrals of the first and 
second kind: 
\beq
\int_1^z \sqrt{f(t)} \, dt \ & = &                                                               
2 \sqrt{\frac{z^2 -1}{z-a}}
   + (1-a) \, \sqrt{2} \; \mbox{F}\left[ \arcsin{\left( \sqrt{
\frac{z-1}{z-a}} \right)},
   \sqrt{\frac{1+a}{2}} \right] \\ & - & 2 \sqrt{2} \; \; \mbox{E}\left[
\arcsin{\left( \sqrt{ \frac{z-1}{z-a}} \right)},
   \sqrt{\frac{1+a}{2}} \right] \,.
\eeq
Note that $\zeta \to -\infty$ 
corresponds to $z \to \infty$ and $\zeta \to 0-$ to $z \to 1+$. Eq.~(\ref{WKB-ODE}) 
becomes
\beqn
\label{ODE-IV}
  \frac{d^2}{d\zeta^2} W(\zeta)   & = & \left[ \frac{u^2}{4\zeta} -
\frac{1}{4\zeta^2}
+ \frac{\psi(\zeta)}{\zeta} \right] W(\zeta) \;,\\
\nonumber
  \psi(\zeta) & = & \frac{1}{16\zeta} - \frac{g(z)}{4f(z)} - \frac{ 4f(z)
f''(z) - 5f'(z)^2}{64 f(z)^3} \;.
\eeqn
Approximate solutions of (\ref{WKB-ODE}), i.e. solutions of (\ref{ODE-IV})
with $\psi=0$, that are regular at $z=1$, are then given by
\cite[Sect. 12.4, Th. 4.1, $\nu=0, \, n=0$]{Olve:1974}, i.e.
\beqn
\label{sol-IV}
w(z) & = & \frac{1}{\sqrt[4]{4|\zeta|f(z)}} \left[ |\zeta|^\frac{1}{2} \,
\mbox{J}_0(u |\zeta|^\frac{1}{2}) + \varepsilon(u,\zeta) \right] \\
\nonumber
\mbox{with} \quad |\varepsilon(u,\zeta)| & \le & \frac{\lambda_{0,1}}{\lambda_{0,0}} \, |\zeta|^\frac{1}{2} \frac{M_0(u\zeta)}{E_0(u\zeta)} \, \left[ e^{\frac{\lambda_{0,0}}{u} \mathcal{V}_{\zeta,0}\left(|\zeta|^\frac{1}{2} \mathcal{B}_0(\zeta)\right)} -1 \right] \;. 
\eeqn
The bound for the error $\varepsilon$ is derived in App.~\ref{error_bounds}.
The constants $\lambda_{0,0}$,$\lambda_{0,1}$, as well as the function $\mathcal{B}_0$
are defined as follows:
\beq   
\lambda_{0,0} & := & \sup_{x\ge0}\left\lbrace \pi x M_0^2(x)\right\rbrace \,,\\
\lambda_{0,1} & := & \sup_{x\ge0}\left\lbrace \pi x M_0^2(x) \,,
     \cos{\theta_0(x)}\right\rbrace \,,\\       
\mathcal{B}_0(\zeta) & := & \frac{1}{|\zeta|^\frac{1}{2}} \int_\zeta^0
     \frac{dv}{|v|^\frac{1}{2}} \, \psi(v) \,,
\eeq
and $\mathcal{V}$ is the variational operator, i.e.
\beq
\mathcal{V}_{\zeta,0}\left(|\zeta|^\frac{1}{2} \mathcal{B}_0(\zeta)\right)
& = & \int_\zeta^0 \frac{dv}{|v|^\frac{1}{2}} |\psi(v)| \;.
\eeq 
We determine the behavior of the solutions for 
$z \to 1$ and $z \to \infty$:
\beq
A(z) & \ub{\longrightarrow}{z\to 1}     & \frac{1}{\sqrt{2}} \,
\mbox{J}_0\left(u\sqrt{2(1-a)} \, \sqrt{z-1}\right) \to \frac{1}{\sqrt{2}}
\,,\\
A(z) & \ub{\sim}{z\to\infty} & \frac{1}{\sqrt{\pi u} z^\frac{3}{4}}
\sin{\left(u \int_1^z \sqrt{f(t)} dt + \frac{\pi}{4} +  \delta \right)}\,,
\eeq
where we have used that $\lim_{\zeta \to 0-} \varepsilon(u,\zeta) = 0$ and
the phase $\delta$ is determined by $\lim_{\zeta \to -\infty} \varepsilon(u,\zeta)$.
We know that the argument of the $\sin$-function equals $2\sqrt{\beta z} + 
\Delta_{j,0}$ for large $z$, i.e.
\beq 
A(z) & 
\sim & \frac{1}{\sqrt{\pi u} z^\frac{3}{4}} \sin{\left( 2\sqrt{\beta z} +
\Delta_{j,0} \right)} \;.
\eeq
Therefore, we obtain for the scattering phase                                
\beq
\Delta_{j,0}\mid_{a<1} & = & u \lim_{z\to \infty}  \left( \int_1^z 
\sqrt{f(t)} \, dt - 2\sqrt{z}  \right) + \frac{\pi}{4} + \, \delta \\
& = & (1-a) \, \sqrt{2 \beta} \;\;
\mbox{K}\left[\sqrt{\frac{1+a}{2}} \right] -
2 \sqrt{2 \beta} \; \; \mbox{E}\left[\sqrt{\frac{1+a}{2}} \right]
+ \frac{\pi}{4} + \delta \;.
\eeq
In addition, we will derive in App.~\ref{error_bounds}
\beqn
\label{errorIV}
\frac{2|\delta|}{\pi} & \le &
\min{\left\{ 1, \frac{\lambda_{0,1}}{\lambda_{0,0}}
  \left[ e^{ \frac{\lambda_{0,0}}{u} \,
\mathcal{V}_{-\infty,0}\left(|\zeta|^\frac{1}{2} \mathcal{B}_0(\zeta)\right)
} -1 \right]
\right\}} \;.
\eeqn
The phase $\Delta_{j,0}\mid_{a<1}$ differs from the one obtained in
\cite{Mign:1991} by $\frac{\pi}{4}$ and the factor in front 
of the elliptic integral of the first kind.

\section{Error Bounds}
\label{error_bounds}
In this chapter we apply the construction of error bounds presented in 
\cite[Sect. 11.3]{Olve:1974} to case IV (cf. \cite[Sect. 12.4, Ex. 4.4]{Olve:1974}). 
In case IV the Bessel functions
$ |\zeta|^\frac{1}{2} \, \mbox{J}_0(u |\zeta|^\frac{1}{2})$ and
$|\zeta|^\frac{1}{2}\, \mbox{Y}_0(u |\zeta|^\frac{1}{2})$ of first and second
type, are first approximations for $W$. They are exact solutions if $\psi \equiv 0$.
Substituting
$W(\zeta) = |\zeta|^\frac{1}{2}\, \mbox{J}_0(u |\zeta|^\frac{1}{2}) +
\varepsilon(u,\zeta)$, the differential equation becomes
\beq
\frac{d^2}{d\zeta^2}\varepsilon(u,\zeta)  - \left[ \frac{u^2}{4\zeta} -
\frac{1}{4\zeta^2} \right] \varepsilon(u,\zeta)
& = & \frac{\psi(\zeta)}{\zeta}  \left[ \varepsilon(u,\zeta) +
|\zeta|^\frac{1}{2}\, \mbox{J}_0(u |\zeta|^\frac{1}{2}) \right] \;.
\eeq
If we rewrite it as an integral equation we obtain a Volterra integral equation
\beq
\varepsilon(u,\zeta) & = & \int_\zeta^0 \mathsf{K}(\zeta,v) \,
\frac{\psi(v)}{|v|^\frac{1}{2}}
\left[ \varepsilon(v) + |v|^\frac{1}{2} \, \mbox{J}_0(u |v|^\frac{1}{2}) 
\right]
dv\,,\\
\mbox{where} \quad \mathsf{K}(\zeta,v) & = &  \pi |v|^\frac{1}{2} \left[
\mbox{J}_0(u |\zeta|^\frac{1}{2}) \, \mbox{Y}_0(u |v|^\frac{1}{2})
- \mbox{Y}_0(u |\zeta|^\frac{1}{2}) \, \mbox{J}_0(u |v|^\frac{1}{2}) \right]
\;.
\eeq
Using the fact that $E_0$ is a non-increasing function we obtain
\beq
 \forall \zeta \le v \le 0 \,: \quad
| \mathsf{K}(\zeta,v) | & \le & \underbrace{\frac{M_0(u
|\zeta|^\frac{1}{2})}{E_0(u |\zeta|^\frac{1}{2})}}_{=: \, P_0(\zeta)}
\; \underbrace{\pi |v|^\frac{1}{2} M_0(u |v|^\frac{1}{2}) \, E_0(u
|v|^\frac{1}{2})}_{=: \, Q(v)} \;.
\eeq
We introduce the parameters $\kappa_0$, $\kappa$ by
\beq
\kappa_0 & := & \sup_{v\in [\zeta,0]}\left\lbrace P_0(v) \, Q(v) \right\rbrace 
\le \frac{\lambda_{0,0}}{u} \,,\\
\kappa & := & \sup_{v\in [\zeta,0]}\left\lbrace Q(v) \,  \left||v|^\frac{1}{2} 
\,\mbox{J}_0(u |v|^\frac{1}{2})\right| \right\rbrace 
 \le \frac{|\zeta|^\frac{1}{2}\, \lambda_{0,1}}{u} \;.
\eeq
With the help of \cite[Sect. 6.10, Th. 10.2]{Olve:1974} which 
establishes a bound for the solution of a Volterra integral equation we end up with a
bound for $\varepsilon$,
\beqn
\label{error}
\nonumber
\frac{|\varepsilon(u,\zeta)|}{P_0(\zeta)} & \le &  \frac{\kappa}{\kappa_0}
\left[ \exp{\left\{\kappa_0 \,
\mathcal{V}_{\zeta,0}\left(|\zeta|^\frac{1}{2} \mathcal{B}_0(\zeta)\right)
\right\}} -1 \right] \\
& \le & |\zeta|^\frac{1}{2} \, \frac{\lambda_{0,1}}{\lambda_{0,0}}
\left[ \exp{\left\{\frac{\lambda_{0,0}}{u} \,
\mathcal{V}_{\zeta,0}\left(|\zeta|^\frac{1}{2} \mathcal{B}_0(\zeta)\right)
\right\}} -1 \right]\,.
\eeqn
We are interested in an error bound for the phase in the asymptotic
expansion of the solution. Assume that the solution $W$ takes the following
form for large arguments $\zeta$: 
\beq
W(u,\zeta) & \ub{\sim}{\zeta \to \infty} & \left( \frac{2|\zeta|^\frac{1}{2}}{\pi u } 
\right)^\frac{1}{2} \Big((1+\rho) \, 
\sin{(u|\zeta|^\frac{1}{2} + \frac{\pi}{4} +\delta)} + o(1) \Big)\;.
\eeq
As in \cite[Sect. 6.7]{Olve:1974} we rewrite the 
difference of the exact and the approximate solution
as a trigonometric function
\beq
\varepsilon(u,\zeta) & = & W(u,\zeta) - |\zeta|^\frac{1}{2} \, \mbox{J}_0(u
|\zeta|^\frac{1}{2}) \\
& \ub{\sim}{\zeta \to \infty}  & \left( \frac{2|\zeta|^\frac{1}{2}}{\pi u } 
\right)^\frac{1}{2}
\left\{ [1 + \rho ] \sin{\left(u |\zeta|^\frac{1}{2}+\frac{\pi}{4}
+\delta\right)}
- \sin{\left(u |\zeta|^\frac{1}{2}+\frac{\pi}{4}\right)} \right\}\\
& = & \left( \frac{2|\zeta|^\frac{1}{2}}{\pi u} \right)^\frac{1}{2}
        \sigma \, \sin{\left(u |\zeta|^\frac{1}{2}+\frac{\pi}{4} +
\eta\right)} \,,
\eeq
where the new parameters $\sigma$ (always considered to be positive) and $\eta$
are related to $\rho$ and $\delta$ by 
\beq
 (1+\rho) \, e^{i \delta} = 1 + \sigma \, e^{i \eta} \;.
\eeq
By elementary geometry it follows that $\frac{2 |\delta|}{\pi} \le \sigma$.
Then, the bound for $\varepsilon$ in (\ref{error}) gives a bound for 
$\sigma$. Now, we choose a sequence $(\zeta_n)$ with $\lim_{n\to \infty}\zeta_n = - \infty$
for which $u |\zeta|^\frac{1}{2}+\frac{\pi}{4} +\eta$ is an odd integer multiple of
$\frac{\pi}{2}$. This shows that 
\beq
 \sigma_{\zeta \to -\infty} = \lim_{\zeta \to - \infty} \left(\frac{\pi u}{2
|\zeta|^\frac{1}{2}} \right)^\frac{1}{2} \, |\varepsilon(u,\zeta)| \;.
\eeq
Finally, since
\beq
 M_0(x) & \ub{\sim}{x \to \infty} & \left( \frac{2}{\pi x} \right)^\frac{1}{2} \,,\\
 E_0(x) & \ub{=}{x \to \infty} & 1\,,
\eeq
we obtain a bound for $\delta$:
\beq
\frac{2|\delta|}{\pi} 
  & \le & \min{\left\{ 1, \frac{\lambda_{0,1}}{\lambda_{0,0}}
  \left[ e^{ \frac{\lambda_{0,0}}{u} \,
\mathcal{V}_{-\infty,0}\left(|\zeta|^\frac{1}{2} \mathcal{B}_0(\zeta)\right)
} -1 \right]
\right\}}\;.
\eeq

\end{appendix}
\small
\addcontentsline{toc}{section}{References}
\bibliographystyle{alpha}
\bibliography{article}

\end{document}